\documentclass[11pt,twoside]{article}
\usepackage[hmargin=0.8in,vmargin=0.9in]{geometry}
\usepackage{fancyhdr}
\usepackage{graphicx}
\usepackage{subfigure}
\usepackage{amssymb}
\usepackage{amsmath}
\usepackage{amsfonts}
\usepackage{theorem}
\usepackage{mathrsfs}
\usepackage{mathtools}
\usepackage{bm}
\usepackage{color}
\usepackage{setspace}
\usepackage{exscale}
\usepackage{relsize}
\usepackage{epstopdf}
\usepackage{float}
\usepackage{cite}

\usepackage{color}

\usepackage{booktabs,multirow} % for much better looking tables
\usepackage{array} % for better arrays (eg matrices) in maths
\usepackage{paralist} % very flexible & customisable lists (eg. enumerate/itemize, etc.)
\usepackage{verbatim} % adds environment for commenting out blocks of text & for better verbatim
\usepackage{subfigure} % make it possible to include more than one captioned figure/table in a single float
\theoremstyle{plain}			% use "default" font
\newtheorem{thm}{Theorem}[section]

{\theorembodyfont{\rmfamily}\newtheorem{remark}{Remark}[section]}

\allowdisplaybreaks[1]

\numberwithin{equation}{section}
\numberwithin{figure}{section}
\numberwithin{table}{section}

\newcommand\eref[1]{(\ref{#1})}

\newcommand*\xbar[1]{%
  \hbox{%
    \vbox{%
      \hrule height 0.5pt % The actual bar
      \kern0.4ex%         % Distance between bar and symbol
      \hbox{%
        \kern-0.05em%      % Shortening on the left side
        \ensuremath{#1}%
        \kern-0.00em%      % Shortening on the right side
      }%
    }%
  }%
}

\newcommand{\mF}{\bm{F}}

\newcommand{\mG}{\bm{G}}

\newcommand{\mH}{\bm{H}}
\newcommand{\mU}{\bm{U}}

\newcommand{\dt}{\Delta t}
\newcommand{\dx}{\Delta x}
\newcommand{\dy}{\Delta y}

\newcommand{\hf}{{\frac{1}{2}}}

\newcommand{\jph}{{j+\frac{1}{2}}}
\newcommand{\jmh}{{j-\frac{1}{2}}}
\newcommand{\kph}{{k+\frac{1}{2}}}
\newcommand{\kmh}{{k-\frac{1}{2}}}

\newtheorem{rmk}[thm]{Remark}
\title{High-Order Flux Splitting Schemes for the Euler Equations of Gas Dynamics}
\author{Shaoshuai Chu\thanks{Department of Mathematics, RWTH Aachen University, 52056, Aachen, Germany;
{\tt chu@igpm.rwth-aachen.de}}, ~Michael Herty\thanks{Department of Mathematics, RWTH Aachen
University, 52056, Aachen, Germany; Department of Mathematics and Applied Mathematics, University of Pretoria, Private Bag X20, Hatfield 0028, South Africa; {\tt herty@igpm.rwth-aachen.de}},~ and Eleuterio F. Toro\thanks{Laboratory of Applied Mathematics,
DICAM, University of Trento, 14-38122,  Trento,  Italy; {\tt  toro@ing.unitn.it}}}
\date{}
\begin{document}
\date{}
\maketitle
\begin{abstract}
We develop high-order flux splitting schemes for the one- and two-dimensional Euler equations of gas dynamics. The proposed schemes are high-order extensions of the existing first-order flux splitting schemes introduced in [{\sc E. F. Toro, M. E. V\'azquez-Cend\'on}, Comput. \& Fluids, 70 (2012), pp. 1--12], where the Euler equations of gas dynamics are split into two subsystems: the advection and pressure systems. In this paper, we formulate the TV splitting within the semi-discrete framework to extend it to higher orders of accuracy for the first time. The second-order extension is obtained by using piecewise linear interpolant to reconstruct the one-sided point values of the unknowns. The third- and fifth-order schemes are developed using the finite-difference alternative weighted essentially non-oscillatory (A-WENO) framework, which is particularly effective in handling multidimensional problems and provides a more straightforward approach to constructing higher-order WENO schemes. These extensions significantly improve the resolution of discontinuities and the accuracy of numerical solutions, as demonstrated by a series of numerical experiments of both the one- and two-dimensional Euler equations of gas dynamics.
\end{abstract}

\noindent
{\bf Key words:} Flux splitting schemes, Euler equations of gas dynamics, advection system, pressure system, A-WENO schemes  

\smallskip
\noindent
{\bf AMS subject classification:} 35L65,  65M06, 65M08, 76M12, 76M20, 76L05, 76N15 

\section{Introduction}
This paper focuses on numerical solutions of the Euler equations of gas dynamics, which in the one-dimensional
(1-D) case, read as 
\begin{equation}
\mU_t+\mF(\mU)_x=\bm0,
\label{1.1}
\end{equation}
with 
\begin{equation}\label{1.2}
 \bm U:=(\rho,\rho u,E)^\top, \quad {\rm and} \quad \bm F(\bm U)=(\rho u,\rho u^2+p,u(E+p))^\top,
\end{equation}
where $t$ is the time, $x$ is the spatial variable, $\rho$, $u$, $p$, and $E$ are the density, velocity, pressure, and total energy, respectively. The system is completed through the following equations of state (EOS): 
\begin{equation}\label{1.3}
p=(\gamma-1)\Big[E-\hf\rho u^2\Big],
\end{equation}
where the parameter $\gamma$ represents the specific heat ratio.

It is well-known that equations \eref{1.1}--\eref{1.3} model a wide range of physical phenomena such as shock waves, turbulence, and compressible flows. These systems often involve complex wave patterns, including discontinuities like shocks and rarefactions, even when the initial data are very smooth, which presents significant challenges for numerical methods. Traditional approaches, particularly first-order methods, may struggle to resolve these discontinuities accurately, often requiring very fine meshes to achieve satisfactory results.

Since the pioneering works of \cite{Fri,Lax,Godunov59}, numerous methods have been developed to solve hyperbolic systems like the one in \eref{1.1}; see, e.g., the monographs and review papers \cite{Hesthaven18,Leveque02,KLR20,Shu20,Toro2009,BAF} and the references therein. Here, we focus on flux splitting methods, which are fundamental in computational fluid dynamics due to their ability to decompose fluxes into components corresponding to different wave families, such as shock, contact, and expansion waves. This decomposition enhances the resolution of shock waves and discontinuities, which is crucial for accurate simulations of compressible flows; see, e.g., \cite{ATvL1986, ATR1989, SW1981, VL1982, VL1982a, Lious1993, Lious1996, Lious1998, Lious2006,DT2011, LLN2012, GVM1999, KSFW2011, ZB1993}. While classical flux splitting schemes (e.g., \cite{ATvL1986, ATR1989, SW1981, VL1982, VL1982a}) are effective, they often struggle to resolve intermediate characteristic fields, leading to excessive dissipation or numerical artifacts. To overcome these limitations, more advanced methods, such as the advection upstream splitting method (AUSM), have been developed (see, e.g.,  \cite{Lious1993}) offering better handling of contact waves and improved resolution of wave patterns. AUSM has since garnered considerable attention and undergone refinement, with further developments \cite{Lious1996, Lious1998, Lious2006}. Additionally, a flux splitting approach similar to AUSM was proposed in \cite{ZB1993}, and subsequent advancements have been recorded in works such as \cite{DT2011, LLN2012, GVM1999, KSFW2011}, further enriching the field of computational fluid dynamics.

Recently, a new flux splitting method, known as the TV splitting scheme, was introduced for the 1-D Euler equations \eref{1.1}--\eref{1.3} in \cite{TV2012} and later extended to higher dimensions in \cite{TCL2015}, where the system \eref{1.1}--\eref{1.3} has been divided into two subsystems: the advection and pressure systems. These subsystems are analyzed, and corresponding Godunov-type discretization schemes are formulated. The proposed schemes are characterized by their simplicity, robustness, and accuracy, offering significant improvements over existing flux splitting methods. In particular, they effectively capture contact and shear waves while precisely preserving isolated stationary contacts. Beyond its application to high-dimensional Euler equations, it has also been extended to other systems, including magnetohydrodynamics and shallow water equations; see, e.g., \cite{TCL2015,BMT2016,DBTF2019,TCVS2022,TT2017}.

The TV splitting schemes introduced in \cite{TV2012,TCL2015} offer a robust and accurate first-order approach for solving the Euler equations of gas dynamics. However, their accuracy is inherently limited. These methods typically require very fine meshes to achieve acceptable resolution for shock waves and other important wave structures. As the complexity increases, higher-order methods are necessary to improve the resolution and reduce the computational cost associated with fine grids. In this paper, we extend these schemes to higher orders of accuracy—second, third, and fifth—in the semi-discrete framework for the first time. The second-order extension is achieved through applying piecewise linear interpolant to reconstruct the one-sided point values of the unknowns, while third- and fifth-order schemes are developed within the finite-difference (FD) alternative essentially non-oscillatory (A-WENO) framework. This framework has been proven to be a powerful tool for generalizing low-order finite-volume (FV) schemes to higher-order FD ones, particularly in multidimensional cases. Its ``dimension-by-dimension" reconstruction process simplifies the development of higher-order WENO schemes; see, e.g., \cite{liu17,wang18,WG2024,Wang20,Jiang13,Jiang96}.

The proposed TV splitting schemes are applied to both 1-D and two-dimensional (2-D) Euler equations of gas dynamics. We test the first-, second-, third-, and fifth-order versions of these schemes on various numerical examples. As expected, the resolution improves with higher-order schemes used, especially when transiting from the first-order scheme to the second-order one. Additionally, we compare the studied TV splitting schemes with the Central-Upwind (CU) schemes (see, e.g.,\cite{Kurganov02,Kurganov01,Kurganov07}), or HLL schemes (see, e.g., \cite{HLL1983}), and the HLLC schemes (see, e.g., \cite{TSS1994}), demonstrating the high efficiency of the proposed methods.

The rest of this paper is organized as follows. In \S\ref{sec2}, beginning with a brief overview of the existing first-order TV splitting scheme for the 1-D Euler equations of gas dynamics, we extend it to second-, third-, and fifth-order accuracy. In \S\ref{sec3}, we first introduce the 2-D first-order TV splitting scheme and similarly extend it to high orders. Finally, in \S\ref{sec4}, we present a number of the 1-D and 2-D numerical results to demonstrate their performances.

\section{One-Dimensional TV Splitting Schemes}\label{sec2}
In this section, we first describe the first-order TV splitting scheme for the 1-D Euler equations \eref{1.1}--\eref{1.3} from \cite{TV2012} and then extend it to second-, third-, and fifth-order accuracy. 

\subsection{1-D First-Order TV Splitting Schemes: A Brief Overview}\label{sec2.1}
Supposing the computational domain is covered with uniform cells $C_j:=[x_\jmh,x_\jph]$ with $x_\jph-x_\jmh\equiv\dx$ centered at $x_j=(x_\jmh+x_\jph)/2$, $\,j=1,\ldots,N$, we assume that the cell average values 
\begin{equation*}
  \xbar\mU_j(t):\approx\frac{1}{\dx}\int\limits_{C_j}\mU(x,t)\,{\rm d}x
\end{equation*}
are available at a certain time level $t$. Note that all of the indexed quantities are time-dependent, but from here on, we suppress the time-dependence of all of the indexed quantities for the sake of brevity.

The computed cell averages $\xbar \mU_j$ of the 1-D system \eref{1.1}--\eref{1.3} are evolved in time by solving the following semi-discrete system of ordinary differential equations (ODEs):
\begin{equation}
\frac{{\rm d}\xbar \mU_j}{{\rm d}t}=-\frac{\bm{{\cal F}}^{\rm FV}_\jph-\bm{{\cal F}}^{\rm FV}_\jmh}{\dx},
\label{2.1}
\end{equation}
where $\bm{{\cal F}}^{\rm FV}_\jph=\bm{{\cal F}}^{\rm FV}_\jph\big(\bm U_\jph^-,\bm U_\jph^+\big)$ is the numerical flux, defined by
\begin{equation}\label{2.1a}
\bm{{\cal F}}^{\rm FV}_\jph\big(\bm U_\jph^-,\bm U_\jph^+\big)= \bm{{\cal F}}^A_\jph\big(\bm U_\jph^-,\bm U_\jph^+\big) + \bm{{\cal F}}^P_\jph \big(\bm U_\jph^-,\bm U_\jph^+\big).
\end{equation}
with $\mU^\pm_\jph$ being the left/right-sided point values of $\mU$ at the cell interfaces $x_\jph$. For the first order scheme, we take $\mU^+_\jph=\mU_{j+1}$ and $\mU^-_\jph=\mU_{j}$. 

Here,  $\bm{{\cal F}}^A_\jph\big(\bm U_\jph^-,\bm U_\jph^+\big)$ is the advection flux given by 
\begin{equation*}
  \bm{{\cal F}}^A_\jph\big(\bm U_\jph^-,\bm U_\jph^+\big)=\begin{cases}
                                                          u^*_\jph\begin{pmatrix} \rho^-_\jph \\[0.8ex] 
                                                          (\rho u)^-_\jph\\[0.8ex]\dfrac{1}{2} \rho^-_\jph (u^-_\jph)^2\end{pmatrix}, & \mbox{if } u^*_\jph \ge 0, \\[8.ex]
                                                          u^*_\jph \begin{pmatrix} \rho^+_\jph \\[0.8ex] (\rho u)^+_\jph\\[0.8ex]\dfrac{1}{2} \rho^+_\jph (u^+_\jph)^2\end{pmatrix}, & \mbox{otherwise},
                                                  \end{cases}
\end{equation*}
and  $\bm{{\cal F}}^P_\jph \big(\bm U_\jph^-,\bm U_\jph^+\big)$ is the pressure flux given by 
\begin{equation*}
  \bm{{\cal F}}^P_\jph \big(\bm U_\jph^-,\bm U_\jph^+\big)= \begin{pmatrix} 0 \\[0.8ex] p^*_\jph\\[1.8ex] \dfrac{\gamma  u^*_\jph  p^*_\jph}{\gamma-1} \end{pmatrix},
\end{equation*}
where 
\begin{equation}\label{2.4}
 \begin{aligned}
  u^*_\jph &= \frac{C^+_\jph u^+_\jph - C^-_\jph u^-_\jph} { C^+_\jph-C^-_\jph}-\frac{2}{ C^+_\jph-C^-_\jph} (p^+_\jph - p^-_\jph), \\
  p^*_\jph &= \frac{C^+_\jph p^-_\jph - C^-_\jph p^+_\jph} { C^+_\jph-C^-_\jph}+\frac{C^+_\jph C^-_\jph}{ 2(C^+_\jph-C^-_\jph)} (u^+_\jph - u^-_\jph),\\
  C^\pm_\jph &=\rho^\pm_\jph \bigg(u^\pm_\jph \pm \sqrt{\big(u^\pm_\jph \big)^2+4 \big(c^+_\jph\big)^2}\bigg),
 \end{aligned}
\end{equation}
with  
$$
u^\pm_\jph = \dfrac{(\rho u)^\pm_\jph}{\rho^\pm_\jph}, \,\,\, p^\pm_\jph=(\gamma-1)\Big(E^\pm_\jph-\frac{1}{2}\rho^\pm_\jph (u^\pm_\jph)^2\Big),\,\,\, {\rm and} \,\,\, c^\pm_\jph=\sqrt{{\gamma p^\pm_\jph}/{\rho^\pm_\jph}}.
$$

\subsection{1-D Second-Order TV Splitting Scheme}\label{sec2.2}
We now extend the first-order TV splitting scheme introduced in \S \ref{sec2.1} to the second order of accuracy. The resulting scheme \eref{2.1}--\eref{2.4} achieves second-order accuracy provided that the one-sided point values $\mU^\pm_\jph$, used to compute the numerical flux \(\bm{{\cal F}}^{\rm FV}_\jph\) in \eref{2.1a}, are second-order of accurate. To this end, we approximate $\mU^\pm_\jph$ using a piecewise linear interpolant
\begin{equation}\label{2.5}
\widetilde\mU(x)=\,\xbar\mU_j+(\mU_x)_j(x-x_j),\quad x\in C_j,
\end{equation}
which leads to
\begin{equation}\label{2.5a}
\mU^-_\jph=\,\xbar\mU_j+\frac{\dx}{2}(\mU_x)_j,\quad\mU^+_\jph=\,\xbar\mU_{j+1}-\frac{\dx}{2}(\mU_x)_{j+1}.
\end{equation}
In order to ensure the reconstruction \eref{2.5}--\eref{2.5a} is non-oscillatory, one needs to compute the slopes $(\mU_x)_j$ in \eref{2.5} with the help of a 
nonlinear limiter. In all of the numerical experiments reported in \S\ref{sec4}, we have used a generalized minmod limiter
\cite{lie03,Nessyahu90,Sweby84}:
\begin{equation}
(\mU_x)_j={\rm minmod}\left(\theta\frac{\,\xbar\mU_j-\,\xbar\mU_{j-1}}{\dx},\,\frac{\,\xbar\mU_{j+1}-\,\xbar\mU_{j-1}}{2\dx},\,
\theta\frac{\,\xbar\mU_{j+1}-\,\xbar\mU_j}{\dx}\right),\quad\theta\in[1,2],
\label{2.7}
\end{equation}
applied in a component-wise manner. Here, the minmod function is defined as
\begin{equation}\label{2.7a}
{\rm minmod}(z_1,z_2,\ldots):=\begin{cases}
\min_j\{z_j\}&\mbox{if}~z_j>0\quad\forall\,j,\\
\max_j\{z_j\}&\mbox{if}~z_j<0\quad\forall\,j,\\
0            &\mbox{otherwise}.
\end{cases}
\end{equation}
The parameter $\theta$ in \eref{2.7} is used to control the amount of numerical viscosity present in the resulting scheme, and larger values of $\theta$ correspond to sharper but, in general, more oscillatory reconstructions. In this paper, we use $\theta= 1.3$.

\subsection{1-D Third-Order TV Splitting Schemes}
In this section, we extend the second-order TV splitting scheme introduced in \S \ref{sec2.2} to the third-order of accuracy in the framework of the finite-difference A-WENO scheme introduced in \cite{Jiang13} (see also \cite{liu17,wang18,Wang20,WG2024}), which has been proven to be a powerful tool for generalizing low-order FV schemes to higher-order FD ones.

Following \cite{Jiang13}, the point values $\mU_j$ are evolved in time by solving the following system of ODEs:
\begin{equation}
\frac{{\rm d}\mU_j}{{\rm d}t}=-\frac{{\mH_\jph}-{\mH_\jmh}}{\dx},
\label{2.8}
\end{equation}
where ${\mH_\jph}$ is the (third-order accurate) numerical flux defined by
\begin{equation*}
{\mH_\jph}=\bm{{\cal F}}^{\rm FV}_\jph-\frac{1}{24}(\dx)^2(\mF_{xx})_\jph.
\end{equation*}
Here, $\bm{{\cal F}}^{\rm FV}_\jph$ is the FV numerical flux as in \eref{2.1a} and $(\mF_{xx})_\jph$ is the higher-order correction term used to increase the order of the numerical flux. The correction term $(\mF_{xx})_\jph$ can be approximated using the standard central FDs:
\begin{equation*}
(\mF_{xx})_\jph=\frac{1}{2(\dx)^2}\big[{\mF}_{j-1}-{\mF}_{j}-{\mF}_{j+1}+{\mF}_{j+2}\big].
\end{equation*}
We stress that the resulting scheme is third-order once the one-sided point values $\mU^\pm_\jph$ employed to compute the numerical flux $\bm{{\cal F}}^{\rm FV}_\jph$ are third-order accurate. This can be done by implementing a certain nonlinear limiting procedure like the third-order WENO-type interpolation (see, e.g, \cite{WG2024,wang18,CH_third}) applied to the local characteristic variables; see Appendix \ref{appa} for a detailed explanation.

\subsection{1-D Fifth-Order TV Splitting Schemes}
According to \cite{Jiang13}, to achieve fifth-order accuracy, the point values $\mU_j$ are evolved in time by solving the system \eref{2.8} with the (fifth-order accurate) numerical flux 
\begin{equation*}
{\mH_\jph}=\bm{{\cal F}}^{\rm FV}_\jph-\frac{1}{24}(\dx)^2(\mF_{xx})_\jph+\frac{7}{5760}(\dx)^4(\mF_{xxxx})_\jph,
\end{equation*}
where $\bm{{\cal F}}^{\rm FV}_\jph$ is the FV numerical flux as in \eref{2.1a}, $(\mF_{xx})_\jph$ and $(\mF_{xxxx})_\jph$ are the higher-order correction terms computed by the fourth- and second-order accurate FDs, respectively; see, e.g., \cite{wang18,Wang20,CKX23,CCK23_Adaptive,Chu21}. Here, we have used the following
higher-order correction terms:
\begin{equation*}
\begin{aligned}
&(\mF_{xx})_\jph=\frac{1}{48(\dx)^2}\Big[-5\mF_{j-2}+39\mF_{j-1}-34\mF_j-34\mF_{j+1}+39\mF_{j+2}-5\mF_{j+3}\Big],\\
&(\mF_{xxxx})_\jph=\frac{1}{2(\dx)^4}\Big[\mF_{j-2}-3\mF_{j-1}+2\mF_j+2\mF_{j+1}-3\mF_{j+2}+\mF_{j+3}\Big],
\end{aligned}
\end{equation*}
where $\mF_{j}:=\mF(\mU_{j})$. In order to ensure the resulting scheme is fifth order, the one-sided point values $\mU^\pm_\jph$ employed to compute the numerical flux $\bm{{\cal F}}^{\rm FV}_\jph$ need to be at least fifth-order accurate. This can be done by using a certain nonlinear limiting procedure like the fifth-order WENO-Z interpolation from \cite{DLGW,Gao20,Jiang13,liu17,wang18} applied to the local characteristic variables; see Appendix \ref{appb} for details.

\begin{remark}
Note that we extend the first-order TV splitting schemes to third- and fifth-order accuracy in the framework of FD A-WENO methods. It is also easy to extend it to even higher-order accuracy in this framework; see, e.g. \cite{Gao20}. 
\end{remark}

\section{Two-Dimensional Schemes}\label{sec3}
In this section, we first briefly overview the first-order  TV splitting scheme for the 2-D Euler equations of gas dynamics and then extend it to the higher order of accuracy.

The 2-D Euler equations of gas dynamics read as 
\begin{equation}
\mU_t+\mF(\mU)_x+\mG(\mU)_y=\bm0,
\label{3.1}
\end{equation}
with $\bm U:=(\rho,\rho u,\rho v,E)^\top$, $\bm F(\bm U)=(\rho u,\rho u^2+p,\rho uv,u(E+p))^\top$, and $\bm G(\bm U)=(\rho v,\rho uv,\rho v^2+p,v(E+p))^\top$. Here $v$ is the $y$-velocity and the rest of the notations are the same as in the 1-D case \eref{1.2}--\eref{1.3}. The system is completed through the following
EOS:
\begin{equation}\label{3.2}
p=(\gamma-1)\Big[E-\frac{\rho}{2}(u^2+v^2)\Big].
\end{equation}

\subsection{2-D First-Order TV Splitting Schemes}
Supposing that the computational domain is covered with uniform cells $C_{j,\,k}:=[x_\jmh,x_\jph]\times[y_\kmh,y_\kph]$ centered at
$(x_j,y_k)=\big((x_\jmh+x_\jph)/2$, $(y_\kph+y_\kmh)/2\big)$ with $x_\jph-x_\jmh\equiv\dx$ and $y_\kph-y_\kmh\equiv\dy$ for all $j,k$, we 
assume that the cell averages 
\begin{equation*}
  \xbar\mU_{j,k}(t):\approx\frac{1}{\dx \dy}\int\limits_{C_{j,k}}\mU(x,y,t)\,{\rm d}x{\rm d}y
\end{equation*}
are available at a certain time level $t$. The cell averages $\xbar \mU_{j,k}$ are then evolved in time by numerically solving the following system of ODEs:
\begin{equation}\label{3.2a}
\begin{aligned}
\frac{{\rm d}\xbar \mU_{j,k}}{{\rm d}t}=-\frac{\bm{{\cal F}}^{\rm FV}_{\jph,k}-\bm{{\cal F}}^{\rm FV}_{\jmh,k}}{\dx}-\frac{\bm{{\cal G}}^{\rm FV}_{j,\kph}-\bm{{\cal G}}^{\rm FV}_{j,\kmh}}{\dy}.
\end{aligned}
\end{equation}
Here, $\bm{{\cal F}}^{\rm FV}_{\jph,k}\big(\bm U_{\jph,k}^-,\bm U_{\jph,k}^+\big)$ and $\bm{{\cal G}}^{\rm FV}_{j,\kph}\big(\bm U_{j,\kph}^-,\bm U_{j,\kph}^+\big)$ are the numerical fluxes, defined by
\begin{equation}\label{3.3}
\begin{aligned}
&\bm{{\cal F}}^{\rm FV}_{\jph,k}\big(\bm U_{\jph,k}^-,\bm U_{\jph,k}^+\big)= \bm{{\cal F}}^A_{\jph,k}\big(\bm U_{\jph,k}^-,\bm U_{\jph,k}^+\big) + \bm{{\cal F}}^P_{\jph,k} \big(\bm U_{\jph,k}^-,\bm U_{\jph,k}^+\big),\\
&\bm{{\cal G}}^{\rm FV}_{j,\kph}\big(\bm U_{j,\kph}^-,\bm U_{j,\kph}^+\big)= \bm{{\cal G}}^A_{j,\kph}\big(\bm U_{j,\kph}^-,\bm U_{j,\kph}^+\big) + \bm{{\cal G}}^P_{j,\kph} \big(\bm U_{j,\kph}^-,\bm U_{j,\kph}^+\big),
\end{aligned}
\end{equation}
where  $\mU^\pm_{\jph,k}$ and $\mU^\pm_{j,\kph}$ are the left/right-sided point values of $\mU$ at the cell interfaces $(x_\jph,y_k)$ and $(x_j,y_\kph)$, respectively. In the first-order scheme,  we take $\mU^+_{\jph,k}=\mU_{j+1,k}$, $\mU^+_{j,\kph}=\mU_{j,k+1}$, and $\mU^-_{\jph,k}=\mU^-_{j,\kph}=\mU_{j,k}$. 

Here, $\bm{{\cal F}}^A_{\jph,k}\big(\bm U_{\jph,k}^-,\bm U_{\jph,k}^+\big)$ is the $x$-direction advection flux given by 
\begin{equation*}
  \bm{{\cal F}}^A_{\jph,k}\big(\bm U_{\jph,k}^-,\bm U_{\jph,k}^+\big)=\begin{cases}
                                                          u^*_{\jph,k}\begin{pmatrix} \rho^-_{\jph,k} \\[0.8ex]
                                                           (\rho u)^-_{\jph,k}\\[0.8ex]
                                                           (\rho v)^-_{\jph,k}\\[0.8ex]
                                                          \dfrac{1}{2} \rho^-_\jph \Big[ (u^-_{\jph,k})^2+(v^-_{\jph,k})^2 \Big]\end{pmatrix}, & \mbox{if } u^*_{\jph,k} \ge 0, \\[12.ex]
                                                          u^*_{\jph,k} \begin{pmatrix} \rho^+_{\jph,k} \\[0.8ex]
                                                           (\rho u)^+_{\jph,k}\\[0.8ex]
                                                           (\rho v)^+_{\jph,k}\\[0.8ex]
                                                          \dfrac{1}{2} \rho^+_\jph \Big[ (u^+_{\jph,k})^2+(v^+_{\jph,k})^2 \Big]\end{pmatrix}, & \mbox{otherwise},
                                                  \end{cases}
\end{equation*}
and $\bm{{\cal F}}^P_{\jph,k} \big(\bm U_{\jph,k}^-,\bm U_{\jph,k}^+\big)$ is the $x$-direction pressure flux given by 
\begin{equation*}
  \bm{{\cal F}}^P_{\jph,k} \big(\bm U_{\jph,k}^-,\bm U_{\jph,k}^+\big)= \begin{pmatrix} 0 \\[0.8ex] p^*_{\jph,k}\\[1.8ex] 0 \\\dfrac{\gamma  u^*_{\jph,k}\,  p^*_{\jph,k}}{\gamma-1} \end{pmatrix},
\end{equation*}
where 
\begin{equation}\label{3.5aa}
 \begin{aligned}
  u^*_{\jph,k} &= \frac{C^+_{\jph,k} u^+_{\jph,k} - C^-_{\jph,k} u^-_{\jph,k}} {C^+_{\jph,k}-C^-_{\jph,k}}-\frac{2}{C^+_{\jph,k}-C^-_{\jph,k}} \big(p^+_{\jph,k} - p^-_{\jph,k}\big), \\
  p^*_{\jph,k} &= \frac{C^+_{\jph,k} p^-_{\jph,k} - C^-_{\jph,k} p^+_{\jph,k}} {C^+_{\jph,k}-C^-_{\jph,k}}+\frac{C^+_{\jph,k} C^-_{\jph,k}}{ 2(C^+_{\jph,k}-C^-_{\jph,k})} \big(u^+_{\jph,k} - u^-_{\jph,k}\big),\\
   C^\pm_{\jph,k}&=\rho^\pm_{\jph,k} \bigg(u^\pm_{\jph,k} \pm \sqrt{\big(u^\pm_{\jph,k} \big)^2+4 \big(c^\pm_{\jph,k}\big)^2}\bigg),
 \end{aligned}
\end{equation}
with 
$$
u^\pm_{\jph,k}=\dfrac{(\rho u)^\pm_{\jph,k}}{\rho^\pm_{\jph,k}}, \,\,\, v^\pm_{\jph,k}=\dfrac{(\rho v)^\pm_{\jph,k}}{\rho^\pm_{\jph,k}}, \,\,\, p^\pm_{\jph,k}=(\gamma-1)\big(E^\pm_{\jph,k}-\frac{1}{2}\rho^\pm_{\jph,k} [(u^\pm_{\jph,k})^2+(v^\pm_{\jph,k})^2]\big),
$$
and 
$$
c^\pm_{\jph,k}=\sqrt{{\gamma p^\pm_{\jph,k}}/{\rho^\pm_{\jph,k}}}.
$$
Similarly, 
$\bm{{\cal G}}^A_{j,\kph}\big(\bm U_{j,\kph}^-,\bm U_{j,\kph}^+\big)$ is the $y$-direction advection flux given by 
\begin{equation*}
  \bm{{\cal G}}^A_{j,\kph}\big(\bm U_{j,\kph}^-,\bm U_{\jph,k}^+\big)=\begin{cases}
                                                          v^*_{j,\kph}\begin{pmatrix} \rho^-_{j,\kph} \\[0.8ex]
                                                           (\rho u)^-_{j,\kph}\\[0.8ex]
                                                           (\rho v)^-_{j,\kph}\\[0.8ex]
                                                          \dfrac{1}{2} \rho^-_\jph \Big[ (u^-_{j,\kph})^2+(v^-_{j,\kph})^2 \Big]\end{pmatrix}, & \mbox{if } v^*_{j,\kph} \ge 0, \\[12.ex]
                                                          v^*_{j,\kph} \begin{pmatrix} \rho^+_{j,\kph} \\[0.8ex]
                                                           (\rho u)^+_{j,\kph}\\[0.8ex]
                                                           (\rho v)^+_{j,\kph}\\[0.8ex]
                                                          \dfrac{1}{2} \rho^+_{j,\kph} \Big[ (u^+_{j,\kph})^2+(v^+_{j,\kph})^2 \Big]\end{pmatrix}, & \mbox{otherwise},
                                                  \end{cases}
\end{equation*}
and $\bm{{\cal G}}^P_{j,\kph} \big(\bm U_{j,\kph}^-,\bm U_{j,\kph}^+\big)$ is the $y$-direction pressure flux given by 
\begin{equation*}
  \bm{{\cal G}}^P_{j,\kph} \big(\bm U_{j,\kph}^-,\bm U_{j,\kph}^+\big)= \begin{pmatrix} 0 \\ 0 \\[0.8ex] p^*_{j,\kph}\\[1.8ex] \dfrac{\gamma  v^*_{j,\kph} \, p^*_{j,\kph}}{\gamma-1} \end{pmatrix},
\end{equation*}
where 
\begin{equation}\label{3.5bb}
 \begin{aligned}
  v^*_{j,\kph} &= \frac{C^+_{j,\kph} v^+_{j,\kph} - C^-_{j,\kph} v^-_{j,\kph}} { C^+_{j,\kph}-C^-_{j,\kph}}-\frac{2}{C^+_{j,\kph}-C^-_{j,\kph}} \big(p^+_{j,\kph} - p^-_{j,\kph}\big), \\[0.8ex]
  p^*_{j,\kph} &= \frac{C^+_{j,\kph} p^-_{j,\kph} - C^-_{j,\kph} p^+_{j,\kph}} { C^+_{j,\kph}-C^-_{j,\kph}}+\frac{C^+_{j,\kph} C^-_{j,\kph}}{ 2(C^+_{j,\kph}-C^-_{j,\kph})} \big(v^+_{j,\kph} - v^-_{j,\kph}\big),\\[1.8ex]
  C^\pm_{j,\kph}&=\rho^\pm_{j,\kph} \bigg(v^-_{j,\kph} \pm \sqrt{\big(v^\pm_{j,\kph} \big)^2+4 \big(c^\pm_{j,\kph}\big)^2}\bigg),
 \end{aligned}
\end{equation}
with 
$$
u^\pm_{j,\kph}=\dfrac{(\rho u)^\pm_{j,\kph}}{\rho^\pm_{j,\kph}},\,\,\, v^\pm_{j,\kph}=\dfrac{(\rho v)^\pm_{j,\kph}}{\rho^\pm_{j,\kph}}, \,\,\, p^\pm_{j,\kph}=(\gamma-1)\big(E^\pm_{j,\kph}-\frac{1}{2}\rho^\pm_{j,\kph} [(u^\pm_{j,\kph})^2+(v^\pm_{j,\kph})^2]\big),$$ 
and 
$$c^\pm_{j,\kph}=\sqrt{{\gamma p^\pm_{j,\kph}}/{\rho^\pm_{j,\kph}}}.$$

\subsection{2-D Second-Order TV Splitting Scheme}
As in the 1-D case, the resulting scheme \eref{3.2a}--\eref{3.3} is second-order accurate once the one-sided point values $\mU^\pm_{\jph,k}$ and $\mU^\pm_{j,\kph}$ employed to compute the numerical fluxes \eref{3.3} are second order. To this end, we approximate $\mU^\pm_{\jph,k}$ and $\mU^\pm_{j,\kph}$ by
\begin{equation*}
\widetilde\mU(x,y)=\,\xbar\mU_{j,k}+(\mU_x)_{j,k}(x-x_j)+(\mU_y)_{j,k}(y-y_k),\quad x\in C_{j,k},
\end{equation*}
which leads to
$$
\begin{aligned}
& \mU^-_{\jph,k}=\,\xbar\mU_{j,k}+\frac{\dx}{2}(\mU_x)_{j,k},\quad\mU^+_{\jph,k}=\,\xbar\mU_{j+1,k}-\frac{\dx}{2}(\mU_x)_{j+1,k},\\
& \mU^-_{j,\kph}=\,\xbar\mU_{j,k}+\frac{\dy}{2}(\mU_y)_{j,k},\quad\mU^+_{j,\kph}=\,\xbar\mU_{j+1,k}-\frac{\dy}{2}(\mU_y)_{j+1,k},\\
\end{aligned}
$$
where
\begin{equation*}
\begin{aligned}
(\mU_x)_{j,k}={\rm minmod}\left(\theta\frac{\,\xbar\mU_{j,k}-\,\xbar\mU_{j-1,k}}{\dx},\,\frac{\,\xbar\mU_{j+1,k}-\,\xbar\mU_{j-1,k}}{2\dx},\,
\theta\frac{\,\xbar\mU_{j+1,k}-\,\xbar\mU_{j,k}}{\dx}\right),\\[1.ex]
(\mU_y)_{j,k}={\rm minmod}\left(\theta\frac{\,\xbar\mU_{j,k}-\,\xbar\mU_{j,k-1}}{\dy},\,\frac{\,\xbar\mU_{j,k+1}-\,\xbar\mU_{j,k-1}}{2\dy},\,
\theta\frac{\,\xbar\mU_{j,k+1}-\,\xbar\mU_{j,k}}{\dy}\right).
\end{aligned}
\end{equation*}
Here, the minmod function is defined by \eref{2.7a}.

\subsection{2-D Third-Order TV Splitting Schemes}
Following \cite{Jiang13}, the point values $\mU_{j,k}$ are evolved in time by solving the following system of ODEs:
\begin{equation}\label{3.5}
\frac{{\rm d}\mU_{j,k}}{{\rm d}t}=-\frac{{\mH_{\jph,k}}-{\mH_{\jmh,k}}}{\dx}-\frac{{\mH_{j,\kph}}-{\mH_{j,\kph}}}{\dy},
\end{equation}
where the numerical fluxes $\mH_{\jph,k}$ and $\mH_{j,\kph}$ are defined by
\begin{equation*}
{\mH_{\jph,k}}=\bm{{\cal F}}^{\rm FV}_{\jph,k}-\frac{1}{24}(\dx)^2(\mF_{xx})_{\jph,k}, \quad  {\mH_{j,\kph}}=\bm{{\cal G}}^{\rm FV}_{j,\kph}-\frac{1}{24}(\dy)^2(\mG_{yy})_{j,\kph}.
\end{equation*}
Here, $\bm{{\cal F}}^{\rm FV}_{\jph,k}$ and ${\mG^{\rm FV}_{j,\kph}}$ are the FV numerical fluxes as in \eref{3.3}, $(\mF_{xx})_{\jph,k}$ and $(\mG_{yy})_{j,\kph}$ are the higher-order correction terms computed by the standard central FDs:
\begin{equation*}
\begin{aligned}
&(\mF_{xx})_{\jph,k}=\frac{1}{12(\dx)^2}\big[\mF_{j-1,k}-\mF_{j,k}-\mF_{j+1,k}+\mF_{j+2,k}\big],\\
&(\mG_{yy})_{j,\kph}=\frac{1}{12(\dy)^2}\big[\mG_{j,k-1}-\mG_{j,k}-\mG_{j,k+1}+\mG_{j,k+2}\big].\\
\end{aligned}
\end{equation*}
To ensure the resulting scheme is third-order accuracy, the one-sided point values $\mU^\pm_{\jph,k}$ and $\mU^\pm_{j,\kph}$ are also computed using third-order WENO-type interpolation applied to the local characteristic variables. Note that this can be done in a ``dimension-by-dimension" manner as in the 1-D case, we therefore omit the details for the sake of brevity.  

\subsection{2-D Fifth-Order TV Splitting Schemes}
According to \cite{Jiang13}, the point values $\mU_j$ are evolved in time by solving the system of ODEs \eref{3.5} with the following numerical fluxes ${\mH_{\jph,k}}$ and ${\mH_{j,\kph}}$:
\begin{equation*}
\begin{aligned}
& {\mH_{\jph,k}}=\bm{{\cal F}}^{\rm FV}_{\jph,k}-\frac{1}{24}(\dx)^2(\mF_{xx})_{\jph,k}+\frac{7}{5760}(\dx)^4(\mF_{xxxx})_{\jph,k},\\
& {\mG_{j,\kph}}=\bm{{\cal G}}^{\rm FV}_{j,\kph}-\frac{1}{24}(\dy)^2(\mG_{yy})_{j,\kph}+\frac{7}{5760}(\dy)^4(\mG_{yyyy})_{j,\kph}.
\end{aligned}
\end{equation*}
Here, $\bm{{\cal F}}^{\rm FV}_{\jph,k}$ and ${\mG^{\rm FV}_{j,\kph}}$ are the FV numerical fluxes as in \eref{3.3}, $(\mF_{xx})_{\jph,k}$,
$(\mF_{xxxx})_{\jph,k}$, $(\mG_{yy})_{j,\kph}$, $(\mG_{yyyy})_{j,\kph}$ are approximations of the second- and fourth-order spatial derivatives of $\mF$ at $(x,y)=(x_\jph,y_k)$ and $\mG$ at $(x,y)=(x_j,y_\kph)$, respectively. In this paper, we have used the following higher-order correction terms from \cite{CCK23_Adaptive}:
\begin{equation*}
\begin{aligned}
&(\mF_{xx})_{\jph,k}=\frac{1}{48(\dx)^2}\left(-5\mF_{j-2,k}+39\mF_{j-1,k}-34\mF_{j,k}-34\mF_{j+1,k}+39\mF_{j+2,k}-5\mF_{j+3,k}\right),\\
&(\mF_{xxxx})_{\jph,k}=\frac{1}{2(\dx)^4}\left(\mF_{j-2,k}-3\mF_{j-1,k}+2\mF_{j,k}+2\mF_{j+1,k}-3\mF_{j+2,k}+\mF_{j+3,k}\right),\\
&(\mG_{yy})_{j,\kph}=\frac{1}{48(\dy)^2}\left(-5\mG_{j,k-2}+39\mG_{j,k-1}-34\mG_{j,k}-34\mG_{j,k+1}+39\mG_{j,k+2}-5\mG_{j,k+3}\right),\\
&(\mG_{yyyy})_{j,\kph}=\frac{1}{2(\dy)^4}\left(\mG_{j,k-2}-3\mG_{j,k-1}+2\mG_{j,k}+2\mG_{j,k+1}-3\mG_{j,k+2}+\mG_{j,k+3}\right),
\end{aligned}
\end{equation*}
where $\mF_{j,k}:=\mF(\mU_{j,k})$ and $\mG_{j,k}:=\mG(\mU_{j,k})$.
To achieve fifth-order accuracy, the one-sided point values $\mU^\pm_{\jph,k}$ and  $\mU^\pm_{j,\kph}$ employed to compute the numerical flux $\bm{{\cal F}}^{\rm FV}_{\jph,k}$ and $\bm{{\cal G}}^{\rm FV}_{j,\kph}$ need to be at least fifth-order accurate. This can also be done in a ``dimension-by-dimension" manner as in the 1-D case, we therefore omit the details for the sake of brevity.

\section{Numerical Examples} \label{sec4}
In this section, we test the studied first-, second-, third-, and fifth-order schemes on several numerical examples and compare their performances. For the sake of brevity, these schemes will be referred to as the 1-, 2-, 3-, and 5-Order schemes, respectively.

We numerically integrate the ODE systems \eref{2.1}, \eref{2.8}, \eref{3.2a}, and \eref{3.5} by the three-stage third-order strong stability preserving (SSP) Runge-Kutta method (see, e.g., \cite{Gottlieb11,Gottlieb12}) and use the CFL number 0.45.

\subsection{One-Dimensional Examples}
We begin with the 1-D Euler equations of gas dynamics \eref{1.1}--\eref{1.3}. In all of the Examples 1--5, we take the specific heat ratio $\gamma=1.4$.

\subsubsection*{Example 1---1-D Accuracy Test}
In the first example, we consider the system \eref{1.1}--\eref{1.3} subject to the following periodic initial
conditions,
\begin{equation*}
\rho(x,0)=1+\frac{1}{10}\sin(2\pi x),\quad u(x,0)\equiv1,\quad p(x,0)\equiv1.
\end{equation*}
The exact solution of this initial value problem is given by
$$
\rho(x,t)=1+\frac{1}{10}\sin\left[2\pi(x-t)\right],\quad u(x,t)\equiv1,\quad p(x,0)\equiv1.
$$

We first compute the numerical solution on the computational domain $[-1,1]$ until the final time $t=0.1$ by the 1-, 2-, 3-, and 5-Order schemes on a sequence of uniform meshes: 100, 200, 400, and 800, measure the $L^1$-errors, and then compute the corresponding experimental convergence rates for the density. The obtained results are presented in Table \ref{tab1},  where one can clearly see that the expected order of accuracy is achieved for the studied schemes.  
\begin{table}[ht!]
\centering
\begin{tabular}{|c|cc|cc|cc|cc|cc|cc|}
\hline
\multirow{2}{2em}{Mesh}&\multicolumn{2}{c|}{1-Order}&\multicolumn{2}{c|}{2-Order}&\multicolumn{2}{c|}{3-Order}&\multicolumn{2}{c|}{5-Order}\\
\cline{2-9}&Error&Rate&Error&Rate&Error&Rate&Error&Rate\\
\hline
$100$&4.93e-03 &---   &4.70e-04 &---  &1.02e-05&--- &1.33e-07 &---\\
$200$&2.49e-03 &0.985 &1.12e-04 &2.07 &1.24e-06&3.04&4.40e-09 &4.92\\
$400$&1.25e-03 &0.993 &2.76e-05 &2.03 &1.55e-07&3.00&1.42e-10&4.95\\
$800$&6.27e-04 &0.996 &6.30e-06 &2.13 &1.94e-08&3.00&4.55e-12&5.00\\
\hline
\end{tabular}
\caption{\sf Example 1: The $L^1$-errors and experimental convergence rates for the density $\rho$ computed by the 1-, 2-, 3-, and 5-Order schemes.\label{tab1}}
\end{table}

\begin{rmk}
We stress that in order to achieve the fifth order of accuracy for the 5-Order scheme, we use smaller time steps with $\dt \sim (\dx)^{\frac{5}{3}}$ to balance the spatial and temporal errors.
\end{rmk}

\subsubsection*{Example 2---Advection of Smooth Density}
In the second example taken from \cite{Toro2009}, we evaluate the efficiency of the four studied schemes. We consider the following initial conditions:
\begin{equation*}
(\rho, u,p)(x,0)=(2+\sin(\pi x))^4,1,1),
\end{equation*}
subject to the period boundary conditions at both ends on the computational interval $[-1,1]$. The exact solution of this initial value problem can be easily
obtained and is given by
$$
(\rho, u,p)(x,t)=(2+\sin(\pi (x-t)))^4,1,1).
$$

\begin{figure}[ht!]
\centerline{\includegraphics[trim=2.8cm .2cm 4.8cm 1.3cm, clip, width=13.4cm]{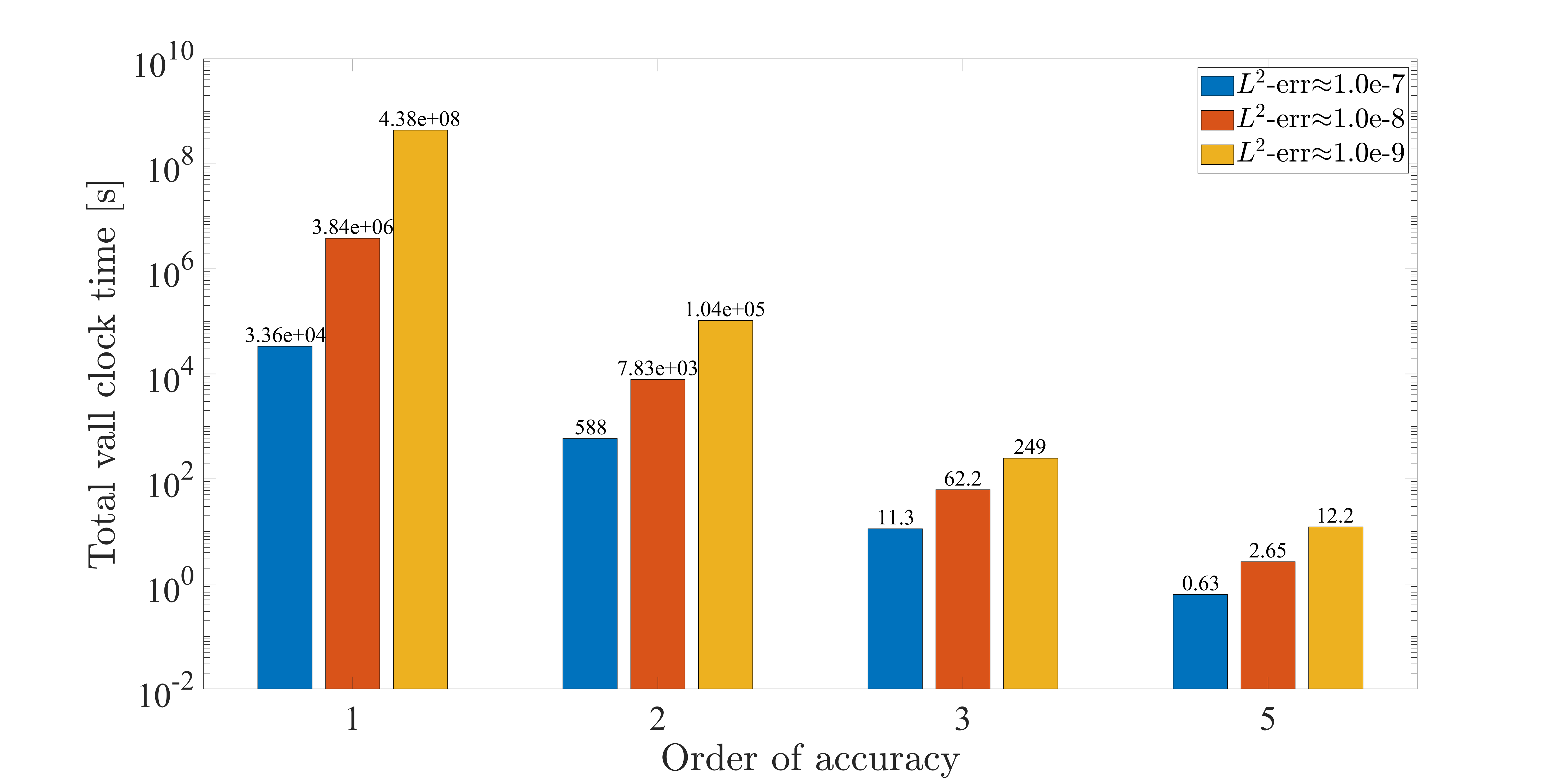}\hspace{1cm}}
\caption{\sf Example 2: Computational cost (bars) against order of accuracy in space and time for three prescribed errors: $10^{-7}$, $10^{-8}$, and $10^{-9}$.
\label{fig2a}}
\end{figure}

In Figure \ref{fig2a}, we present the computational cost (bars) versus the order of accuracy for three cases, each corresponding to prescribed $L^2$-errors ($10^{-7}$, $10^{-8}$, and $10^{-9}$). As one can see, the computational cost decreases significantly with increasing order of accuracy, demonstrating the superior efficiency of high-order methods. In fact, high-order methods can achieve the same accuracy as low-order methods while requiring orders of magnitude less CPU time, making them substantially more efficient

\subsubsection*{Example 3---Shock-Density Wave Interaction Problem.} In this example taken from \cite{SO89}, we consider the shock-density
wave interaction problem with the following initial data,
\begin{equation*}
(\rho,u,p)\Big|_{(x,0)}=\begin{cases}\bigg(\dfrac{27}{7},\dfrac{4\sqrt{35}}{9},\dfrac{31}{3}\bigg),&x<-4,\\[0.8ex]
(1+0.2\sin(5x),0,1),&x>-4,
\end{cases}
\end{equation*}
prescribed in the computational domain $[-5,5]$ subject to the free boundary conditions.

\begin{figure}[ht!]
\centerline{\includegraphics[trim=0.8cm 0.3cm 1.cm 0.5cm, clip, width=6.cm]{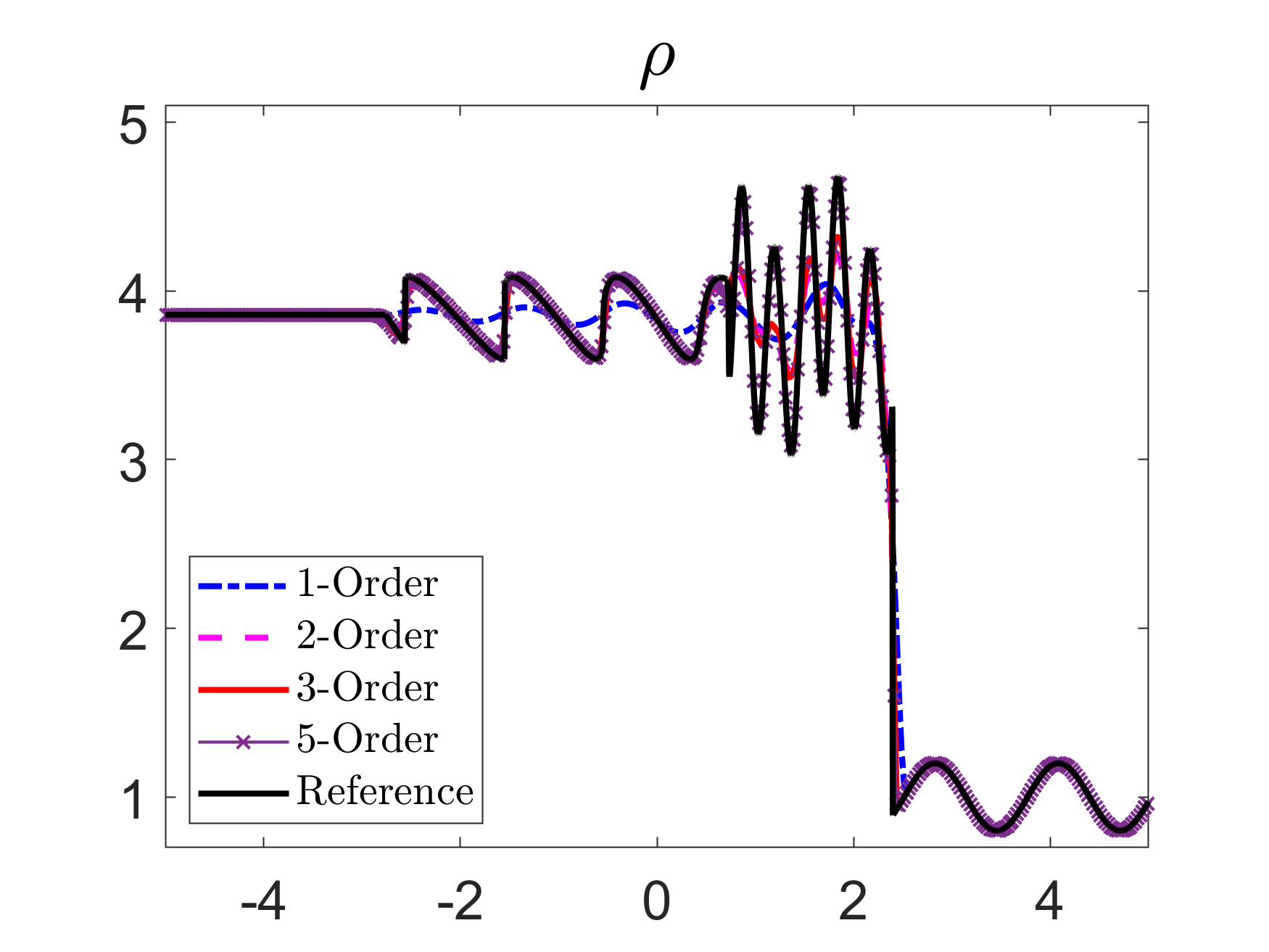}\hspace{1cm}
            \includegraphics[trim=0.8cm 0.3cm 1.cm 0.5cm, clip, width=6.cm]{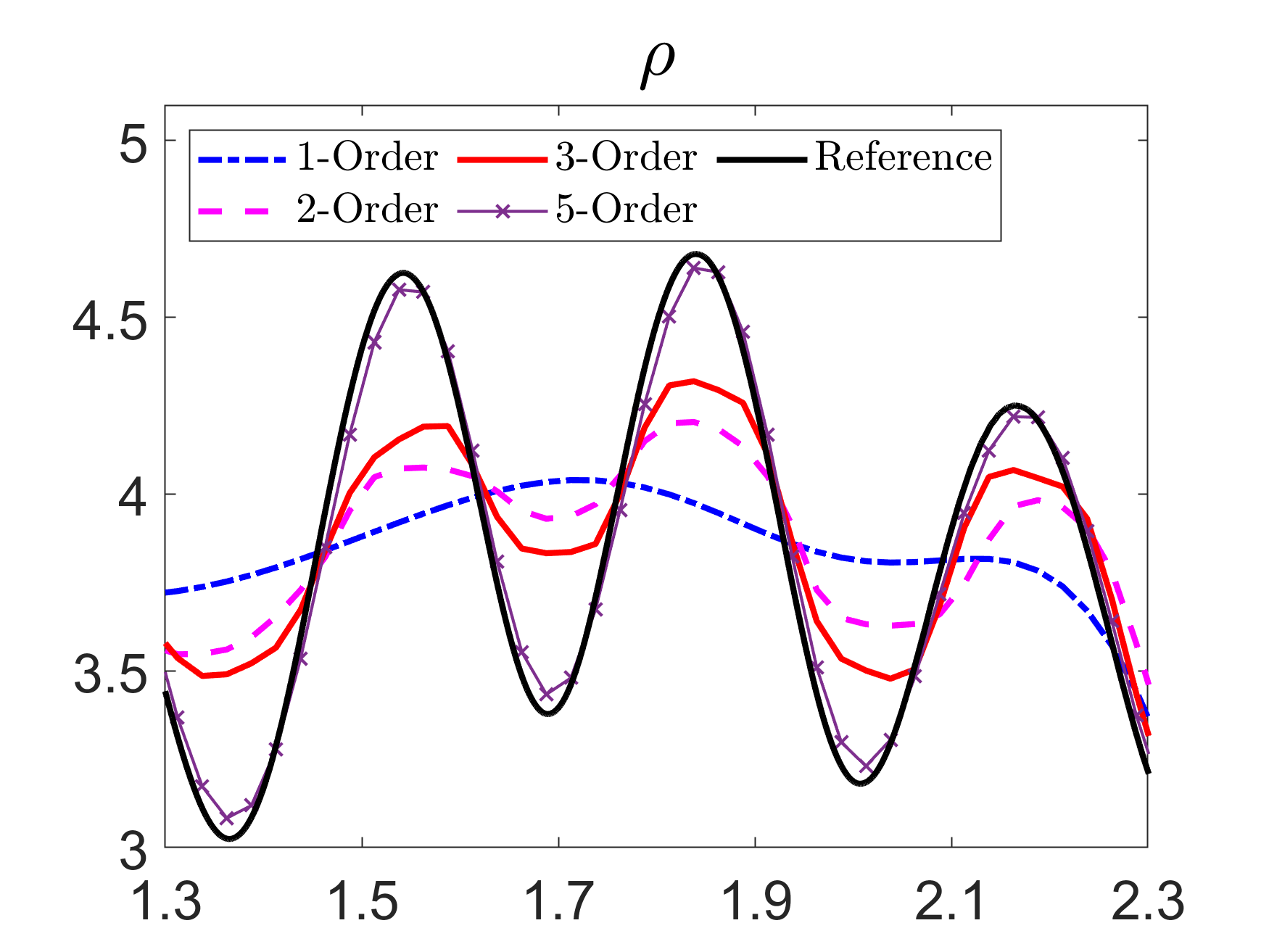}}
\caption{\sf Example 3: Density $\rho$ computed by the 1-, 2-, 3-, and 5-Order schemes (left) and zoom at $[1.3,2.3]$ (right).
\label{fig2}}
\end{figure}

We compute the numerical solutions until the final time $t=5$ by the 1-, 2-, 3-, and 5-Order schemes on a uniform mesh of 400 cells, and present the obtained numerical results in Figure \ref{fig2} together with the reference solution computed by the 5-Order scheme on a much finer mesh of 8000 cells. One can clearly see that the resolution of the computed density improves significantly when high-order schemes are used.

{\color{red}
\subsubsection*{Example 4--- Titarev-Toro Problem}
In this example, we consider the Titarev-Toro problem taken from \cite{Toro2005}; see also \cite{Shu88,Toro2005a}.} The initial conditions,
\begin{equation*}
(\rho,u,p)(x,0)=\begin{cases}
(1.51695,0.523346,1.805),&x<-4.5,\\
(1+0.1\sin(20x),0,1),&x>-4.5,
\end{cases}
\end{equation*}
correspond to a forward-facing shock wave of Mach 1.1 interacting with high-frequency density perturbations, that is, as the shock wave
moves, the perturbations spread ahead. In this example, we set free boundary conditions at both ends of the computational domain $[-10,5]$.

\begin{figure}[ht!]
   \centerline{\includegraphics[trim=0.8cm 0.3cm 0.7cm 0.4cm, clip, width=6.cm]{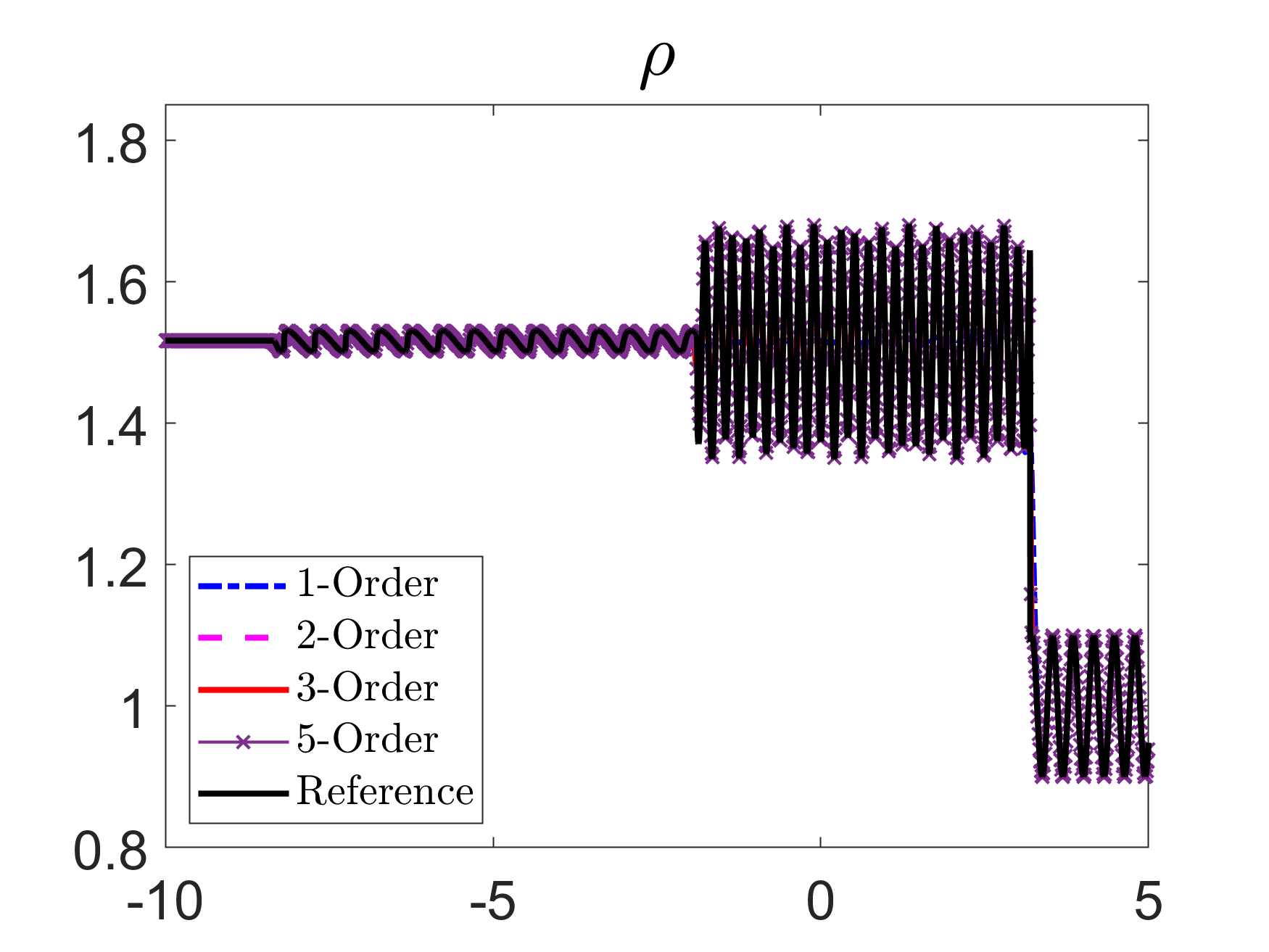}\hspace{0.5cm}
               \includegraphics[trim=0.8cm 0.3cm 0.7cm 0.4cm, clip, width=6.cm]{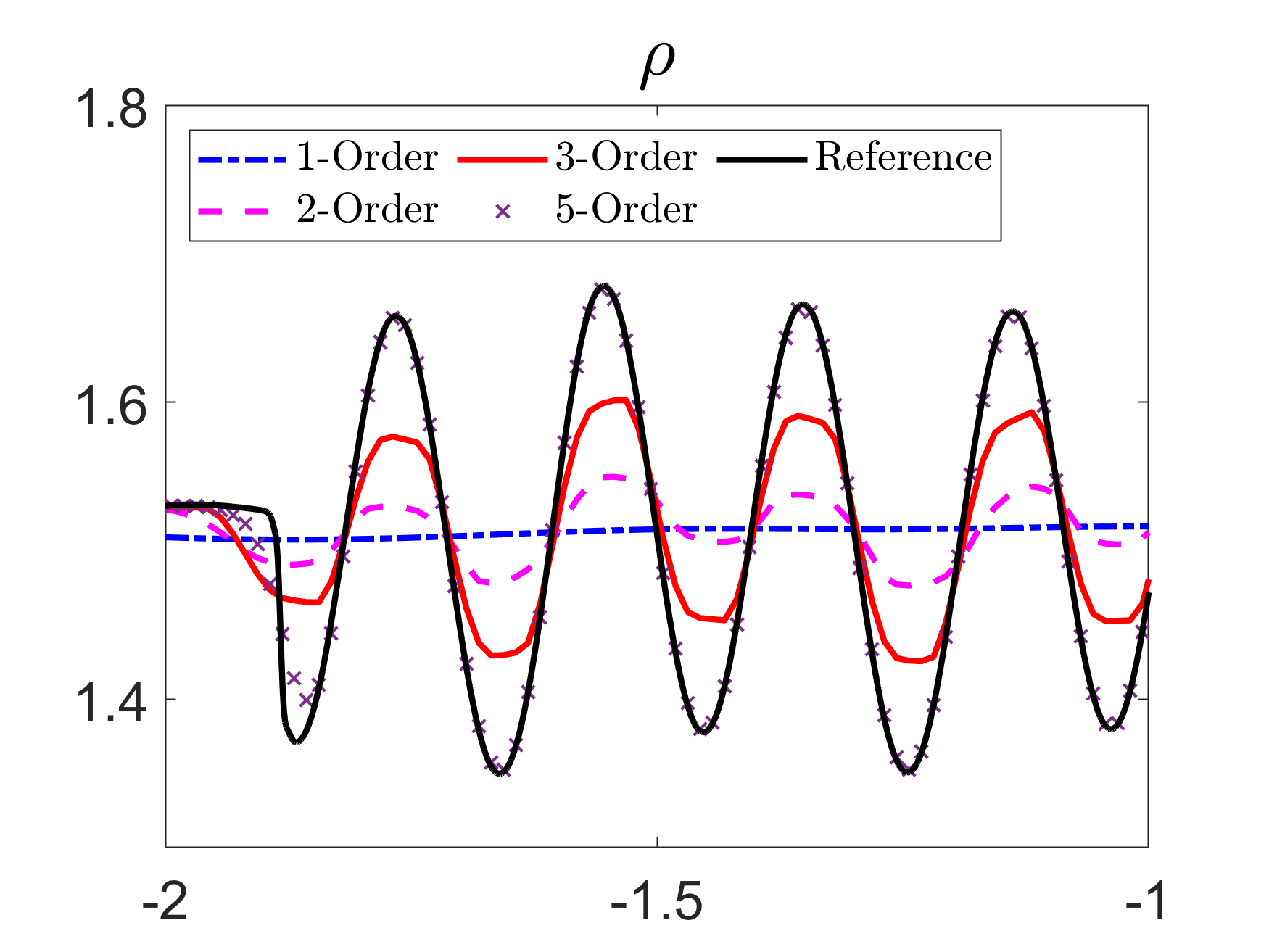}}
   \caption{\sf Example 4: Density $\rho$ computed by the 1-, 2-, 3-, and 5-Order schemes (left) and zoom at $[-2, -1]$ (right).}\label{fig3}
\end{figure}

We compute the solutions until the final time $t=5$ by the 1-, 2-, 3-, and 5-Order schemes on a uniform mesh of 1200 cells. The numerical results are
shown in Figure \ref{fig3} along with the reference solution computed by the 5-Order scheme on a much finer mesh of 12000 cells. The obtained results clearly demonstrate 
a substantial difference in the resolution computed by schemes with different orders. 

\begin{figure}[htb!]
   \centerline{\includegraphics[trim=0.8cm 0.3cm 0.7cm 0.4cm, clip, width=6.cm]{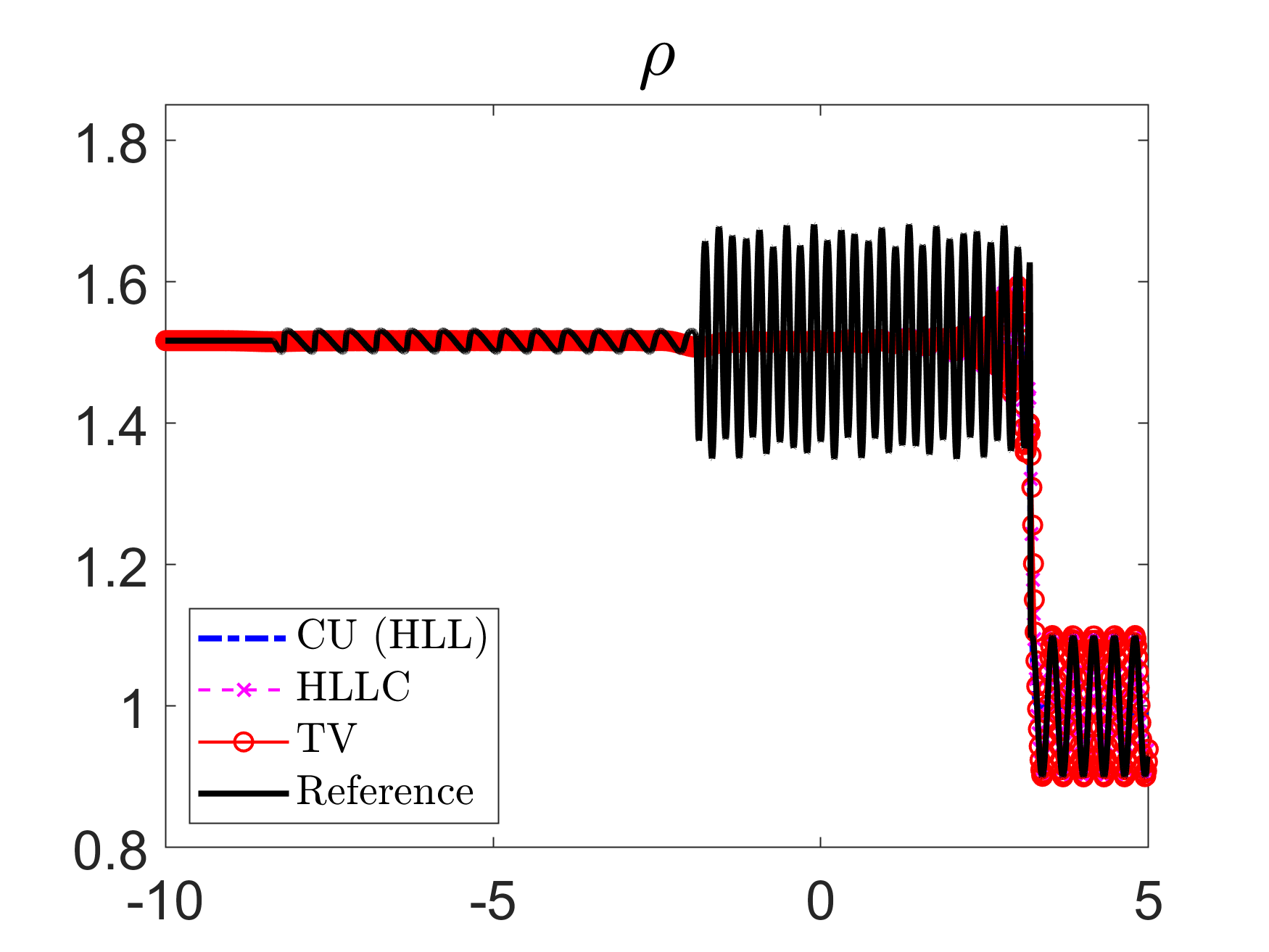}\hspace{0.5cm}
               \includegraphics[trim=0.8cm 0.3cm 0.7cm 0.4cm, clip, width=6.cm]{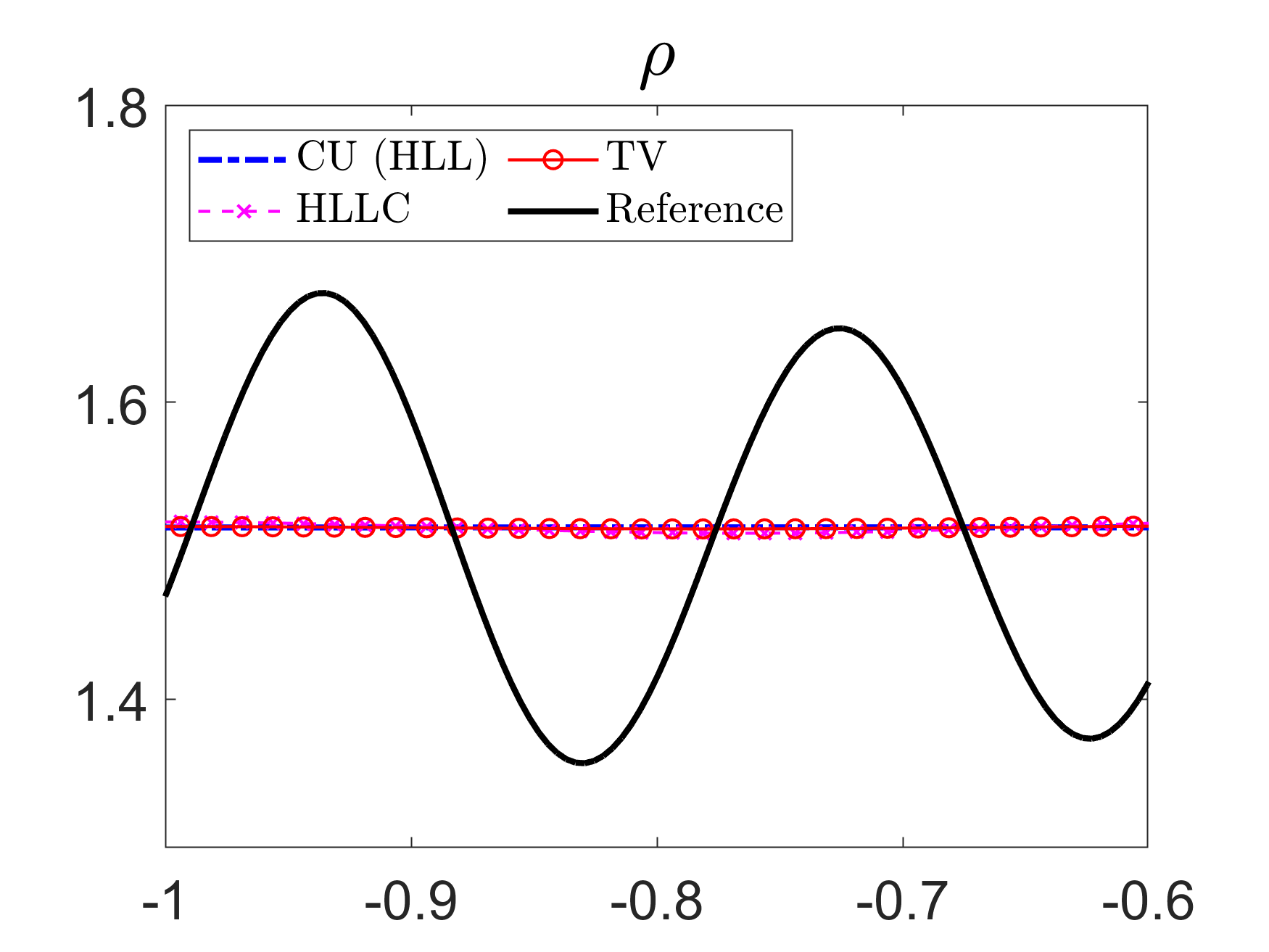}}
   \vskip 8pt
   \centerline{\includegraphics[trim=0.8cm 0.3cm 0.7cm 0.4cm, clip, width=6.cm]{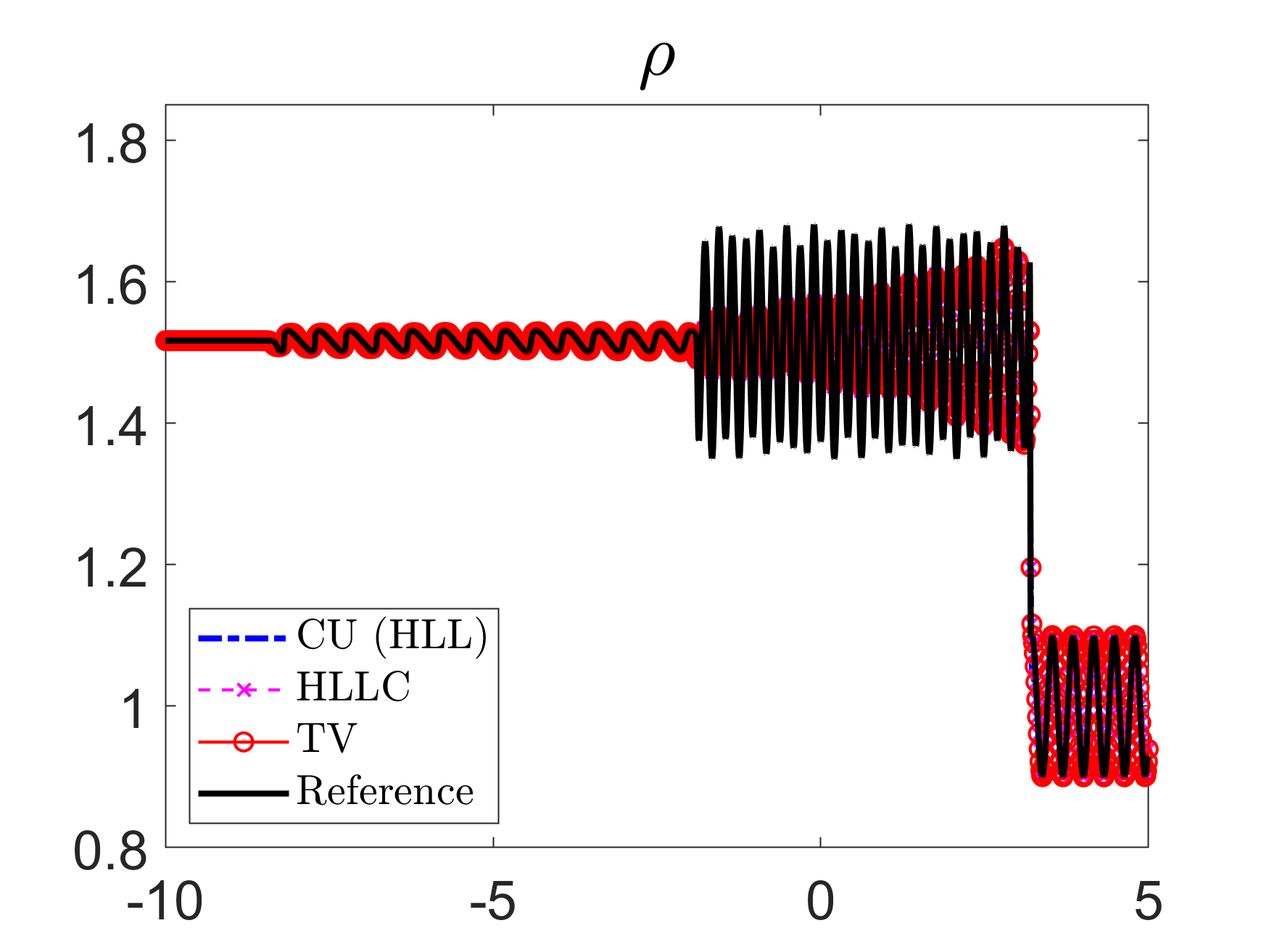}\hspace{0.5cm}
               \includegraphics[trim=0.8cm 0.3cm 0.7cm 0.4cm, clip, width=6.cm]{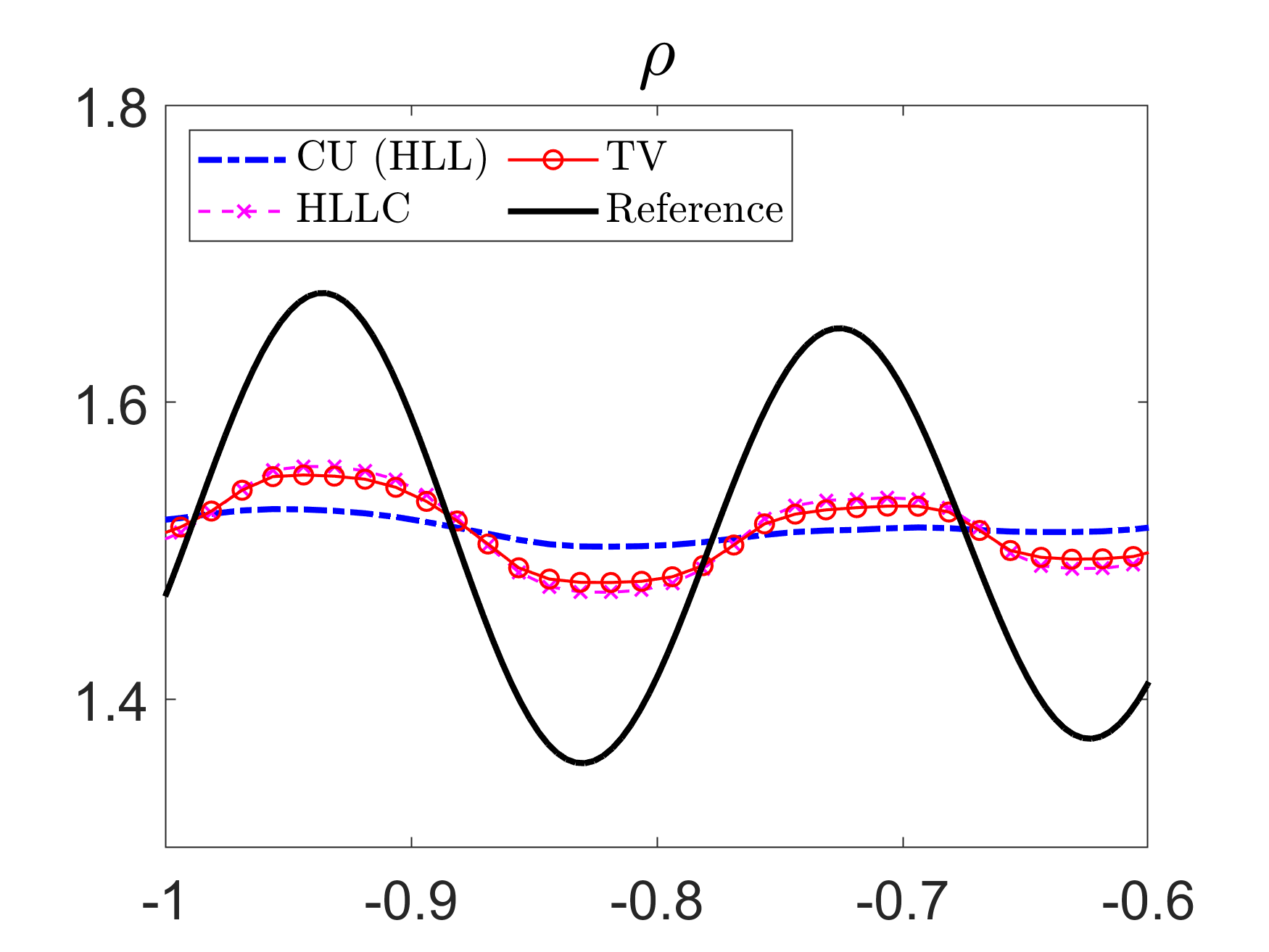}}
    \vskip 8pt
   \centerline{\includegraphics[trim=0.8cm 0.3cm 0.7cm 0.4cm, clip, width=6.cm]{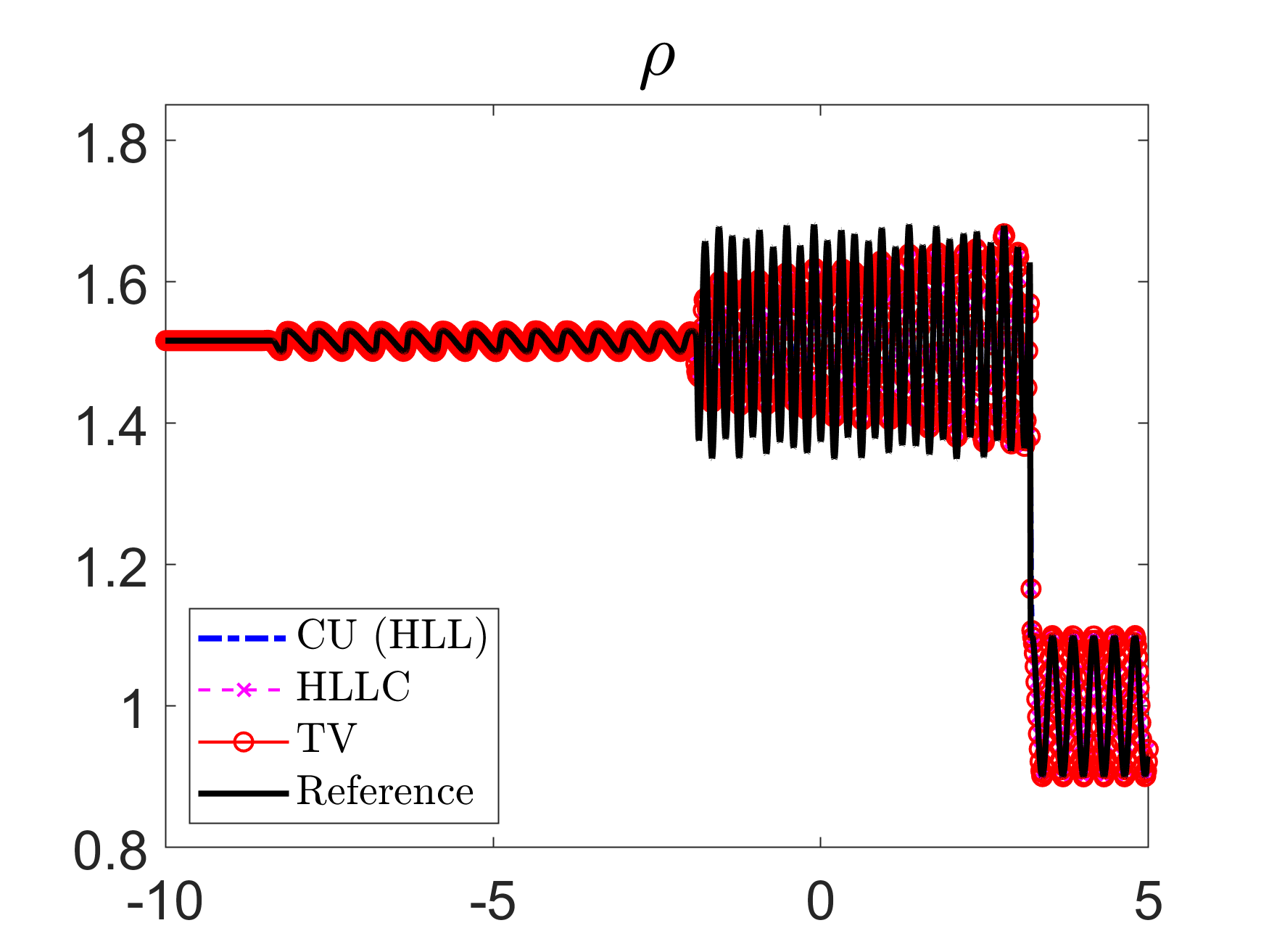}\hspace{0.5cm}
               \includegraphics[trim=0.8cm 0.3cm 0.7cm 0.4cm, clip, width=6.cm]{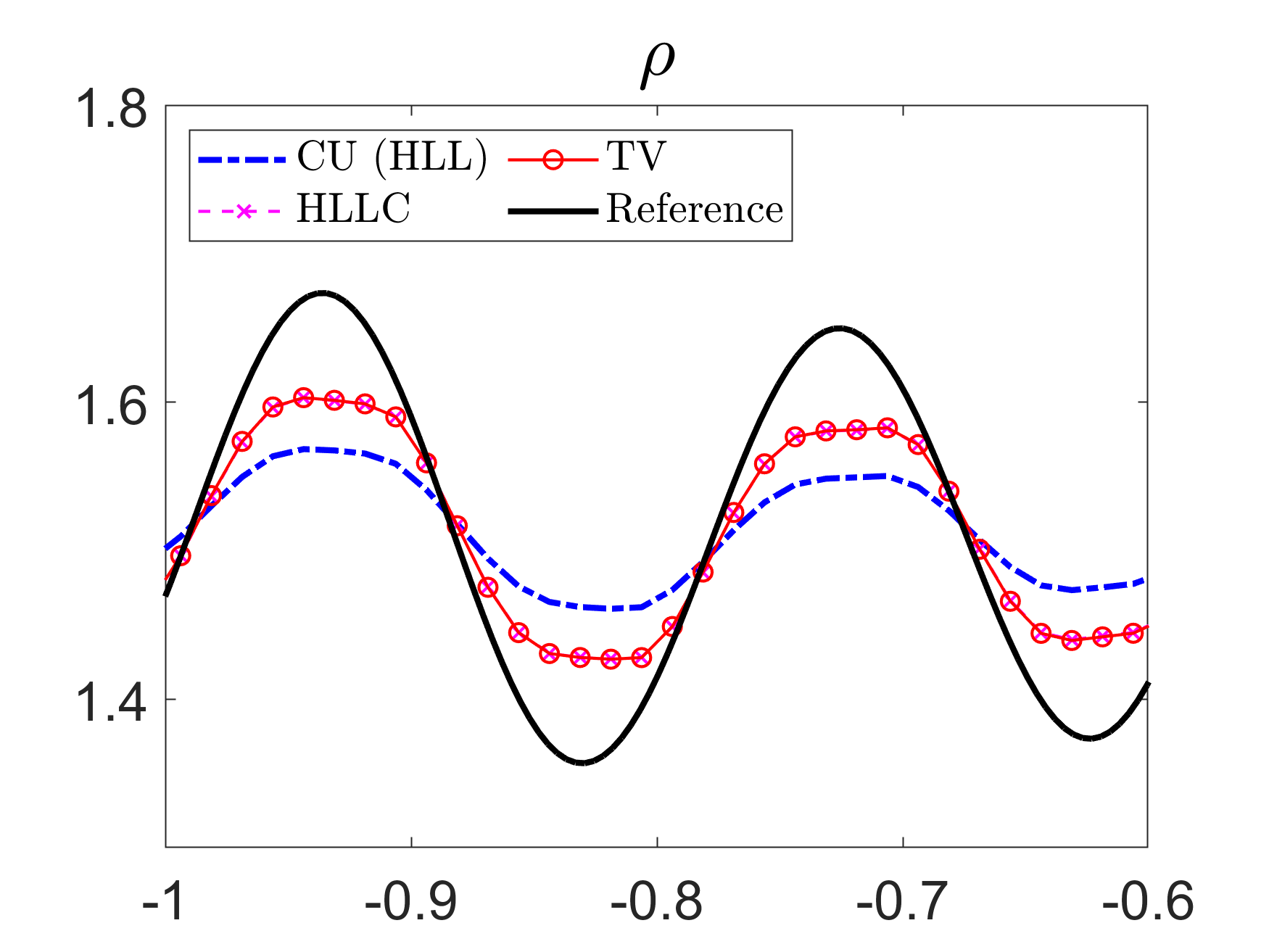}}
    \vskip 8pt
   \centerline{\includegraphics[trim=0.8cm 0.3cm 0.7cm 0.4cm, clip, width=6.cm]{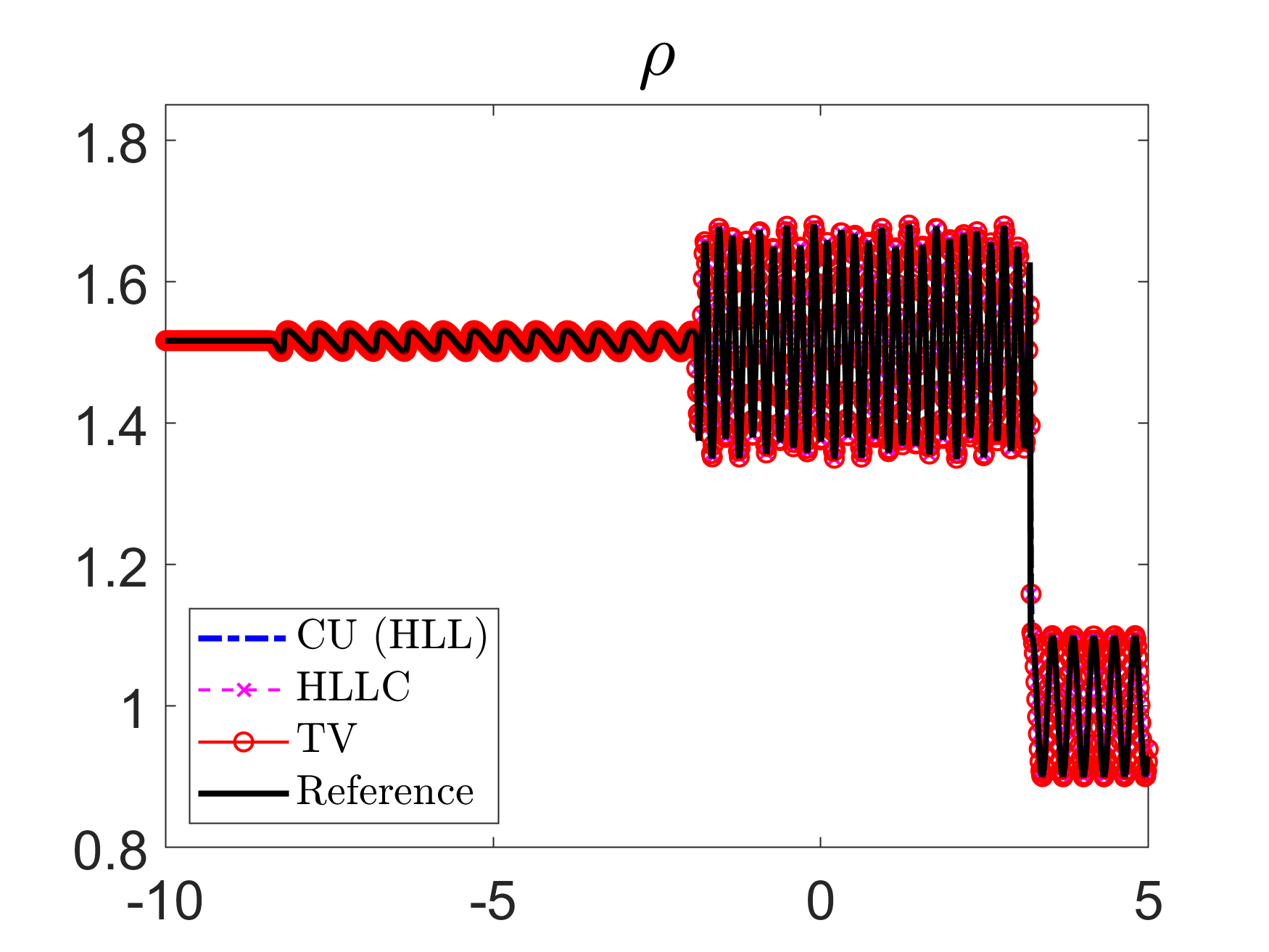}\hspace{0.5cm}
               \includegraphics[trim=0.8cm 0.3cm 0.7cm 0.4cm, clip, width=6.cm]{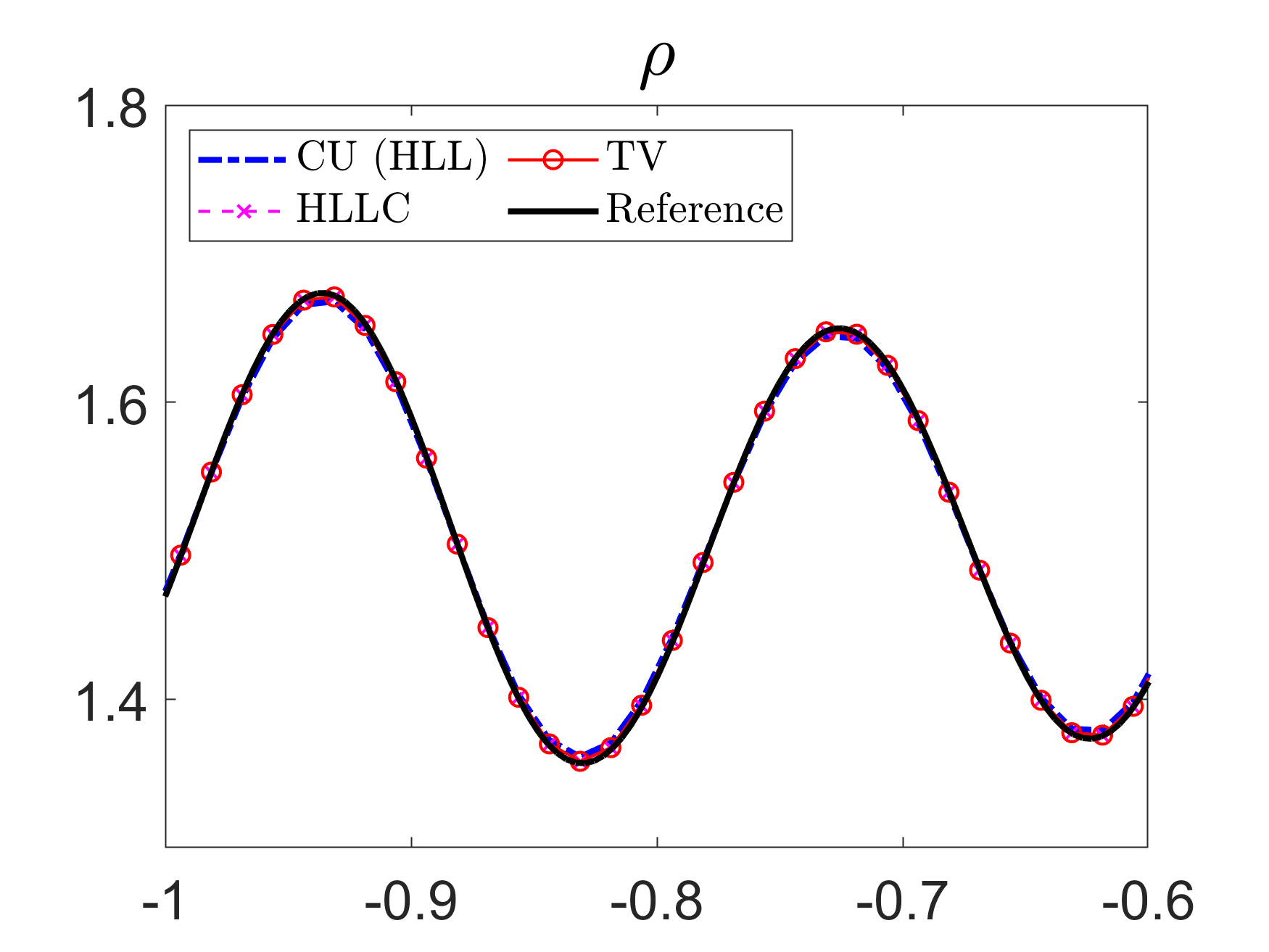}}
   \caption{\sf Example 4: Density $\rho$ computed by the 1- (top row), 2- (second row), 3- (third-row) and 5-Order (bottom row) TV splitting, CU (HLL), and HLLC schemes (left) and zoom at $[-1, -0.6]$ (right).}\label{fig3a}
\end{figure}

In this example, we compare the performances of the studied TV splitting scheme with the CU scheme (see, e.g.,\cite{Kurganov01,Kurganov07}) or  HLL scheme  (see, e.g.,\cite{HLL1983}), and HLLC scheme (see e.g.,\cite{TSS1994}. We compute the solutions by the corresponding CU (HLL) and HLLC schemes and present the obtained numerical densities in Figure \ref{fig3a}, where one can both TV splitting and HLLC schemes achieve better resolutions than the CU (HLL) scheme. At the same time, HLLC scheme contains slightly less numerical dissipations than the studied TV splitting scheme; see the first two rows in Figure \ref{fig3a}, and the results computed by TV splitting and HLLC schemes are almost identical after extended to higher orders; see the last two rows in Figure \ref{fig3a}.

\subsubsection*{Example 5---Blast Wave Problem}
In the last 1-D example,  we consider the strong shocks interaction problem from \cite{Woodward88} with the following initial conditions:
\begin{equation*}
(\rho, u,p)(x,0)=\begin{cases}
(1,0,1000),&x<0.1,\\
(1,0,0.01),&0.1\le x\le 0.9,\\
(1,0,100),&x>0.9,
\end{cases}
\end{equation*}
prescribed in the computational domain $[0,1]$ subject to the solid wall boundary conditions at both ends.

\begin{figure}[ht!]
	\centerline{\includegraphics[trim=0.9cm 0.3cm 0.9cm 0.4cm, clip, width=6.cm]{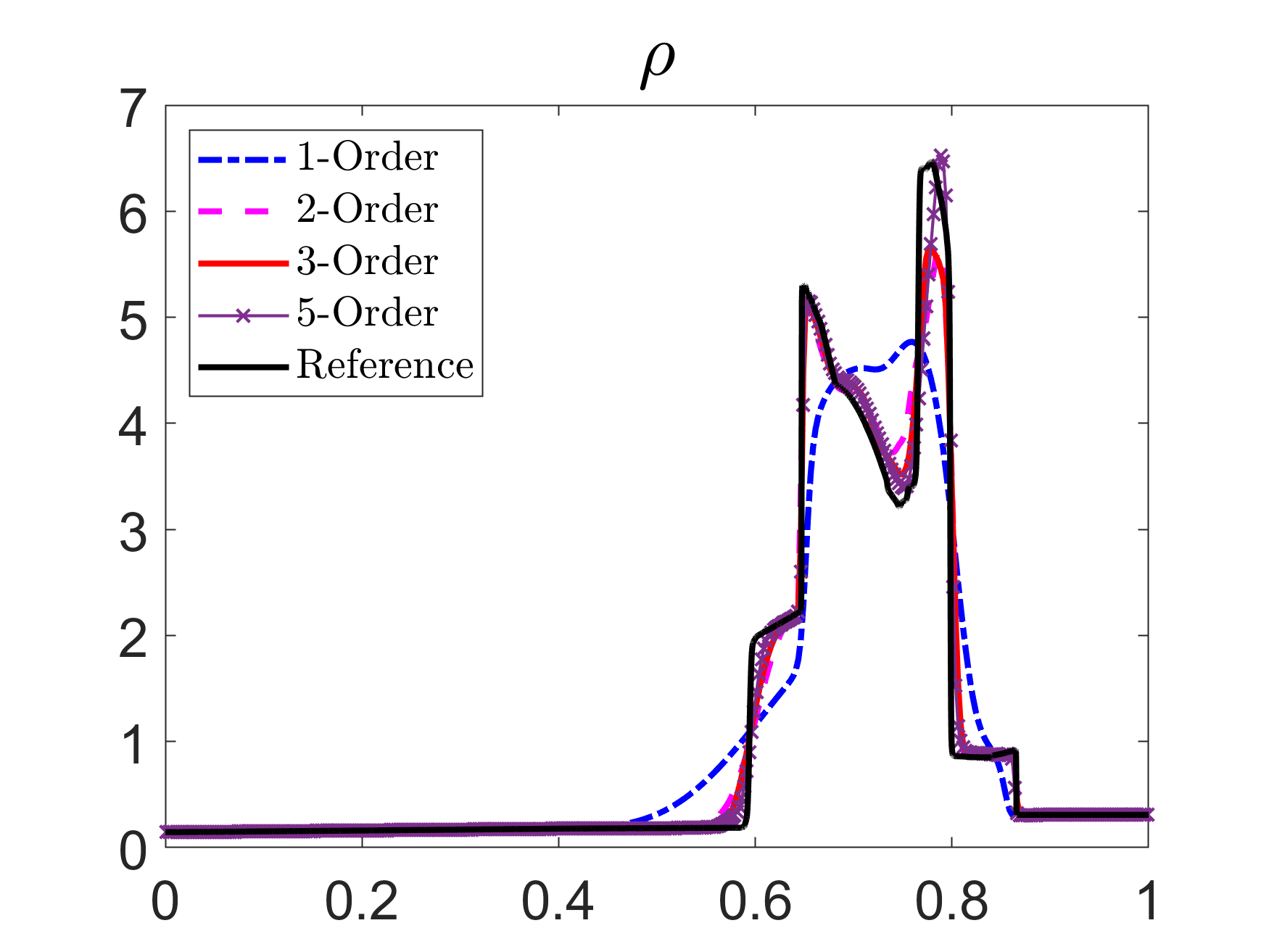}\hspace{0.5cm}
		        \includegraphics[trim=0.9cm 0.3cm 0.9cm 0.4cm, clip, width=6.cm]{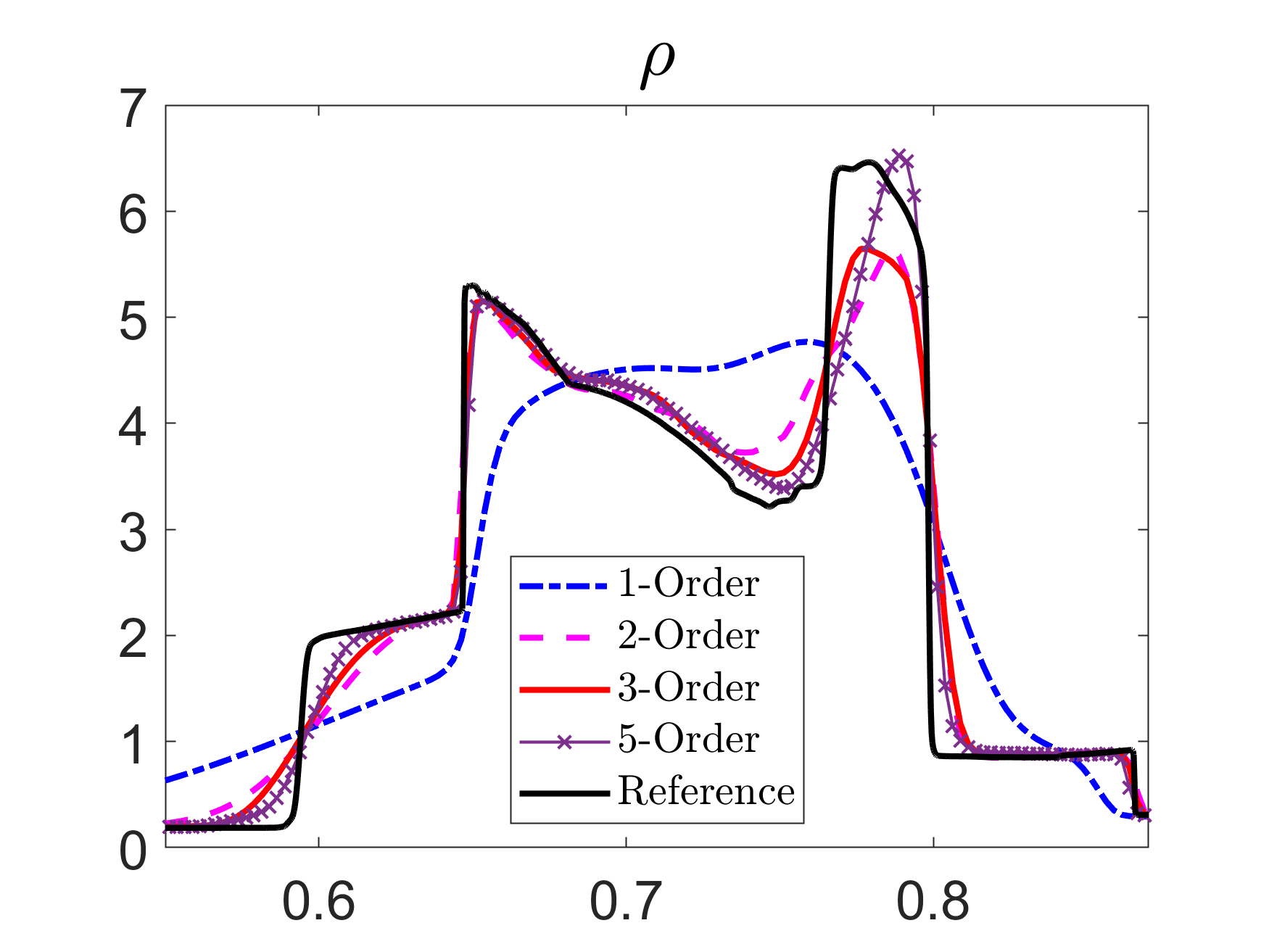}}
	\caption{\sf Example 5: Density $\rho$ computed by the 1-, 2-, 3-, and 5-Order schemes (left) and zoom at $[0.55, 0.87]$ (right).}\label{fig5}
\end{figure}

We compute the numerical solutions until the final time $t=0.038$ by the 1-, 2-, 3-, and 5-Order schemes on a uniform mesh of 400 cells. The obtained results are presented in Figure \ref{fig5} together with the reference solution computed by the 5-Order scheme on a much finer mesh of 4000 cells, demonstrating that the resolution of the computed density improves significantly with the use of high-order schemes, especially when transiting from the 1-Order scheme to the 2-Order one.

\subsection{Two-Dimensional Examples}
In this section, we consider the 2-D Euler equations of gas dynamics \eref{3.1}--\eref{3.2}. In Examples 6--10, we take the specific heat ratio $\gamma=1.4$, while in Example 11, we take $\gamma=5/3$.

\subsubsection*{Example 6---2-D Accuracy Test}
In the first 2-D example taken from \cite{Kurganov17,CCK23_Adaptive,CCHKL_22}, we consider the 2-D Euler equations of gas dynamics subject to the periodic initial
conditions,
\begin{equation*}
\rho(x,y,0)=1+\frac{1}{5}\sin(\pi(x+y)),\quad u(x,y,0)\equiv1,\quad v(x,y,0)\equiv-0.7,\quad p(x,y,0)\equiv1,
\end{equation*}
prescribed on $[-1,1]\times[-1,1]$. The exact solution of this initial value problem can be easily obtained and is given by
$$
\rho(x,y,t)=1+\frac{1}{5}\sin\left[\pi(x+y-0.3t)\right],\quad u(x,y,t)\equiv1,\quad v(x,y,t)\equiv-0.7,\quad p(x,y,0)\equiv1.
$$

We first compute the numerical solution until the final time $t=0.1$ using the 1-, 2-, 3-, and 5-Order schemes on a sequence of uniform meshes: $50\times 50$, $100\times 100$, $200\times 200$, and $400\times 400$. We then measure the $L^1$-errors and the corresponding experimental convergence rates for the density. The obtained results are presented in Table \ref{tab52}, where one can see that the studied 1-, 2-, 3-, and 5-Order schemes achieve the expected order of accuracy. Similar to Example 1, we had to use smaller time steps with $\dt \sim \min\big\{(\dx)^{\frac{5}{3}},(\dy)^{\frac{5}{3}}\big \}$ to achieve the fifth order of accuracy.
\begin{table}[ht!]
\centering
\begin{tabular}{|c|cc|cc|cc|cc|cc|cc|}
\hline
\multirow{2}{2em}{Mesh}&\multicolumn{2}{c|}{1-Order}&\multicolumn{2}{c|}{2-Order}&\multicolumn{2}{c|}{3-Order}&\multicolumn{2}{c|}{5-Order}\\
\cline{2-9}&Error&Rate&Error&Rate&Error&Rate&Error&Rate\\
\hline
$50 \times 50$ &1.68e-02 &---   &1.08e-03 &---  &3.36e-05&--- &2.49e-07 &---\\
$100\times 100$&8.47e-03 &0.986 &2.64e-04 &2.03 &4.21e-06&2.99&7.80e-09 &5.00\\
$200\times 200$&4.25e-03 &0.994 &6.16e-05 &2.10 &5.27e-07&3.00&2.44e-10&5.00\\
$400\times 400$&2.13e-03 &0.997 &1.47e-05 &2.06 &6.59e-08&3.00&7.83e-12&4.96\\
\hline
\end{tabular}
\caption{\sf Example 6: The $L^1$-errors and experimental convergence rates for the density $\rho$ computed by the 1-, 2-, 3-, and 5-Order schemes.\label{tab52}}
\end{table}

\begin{figure}[ht!]
	\centerline{\includegraphics[trim=1.cm 0.9cm 1.4cm 0.2cm, clip, width=6cm]{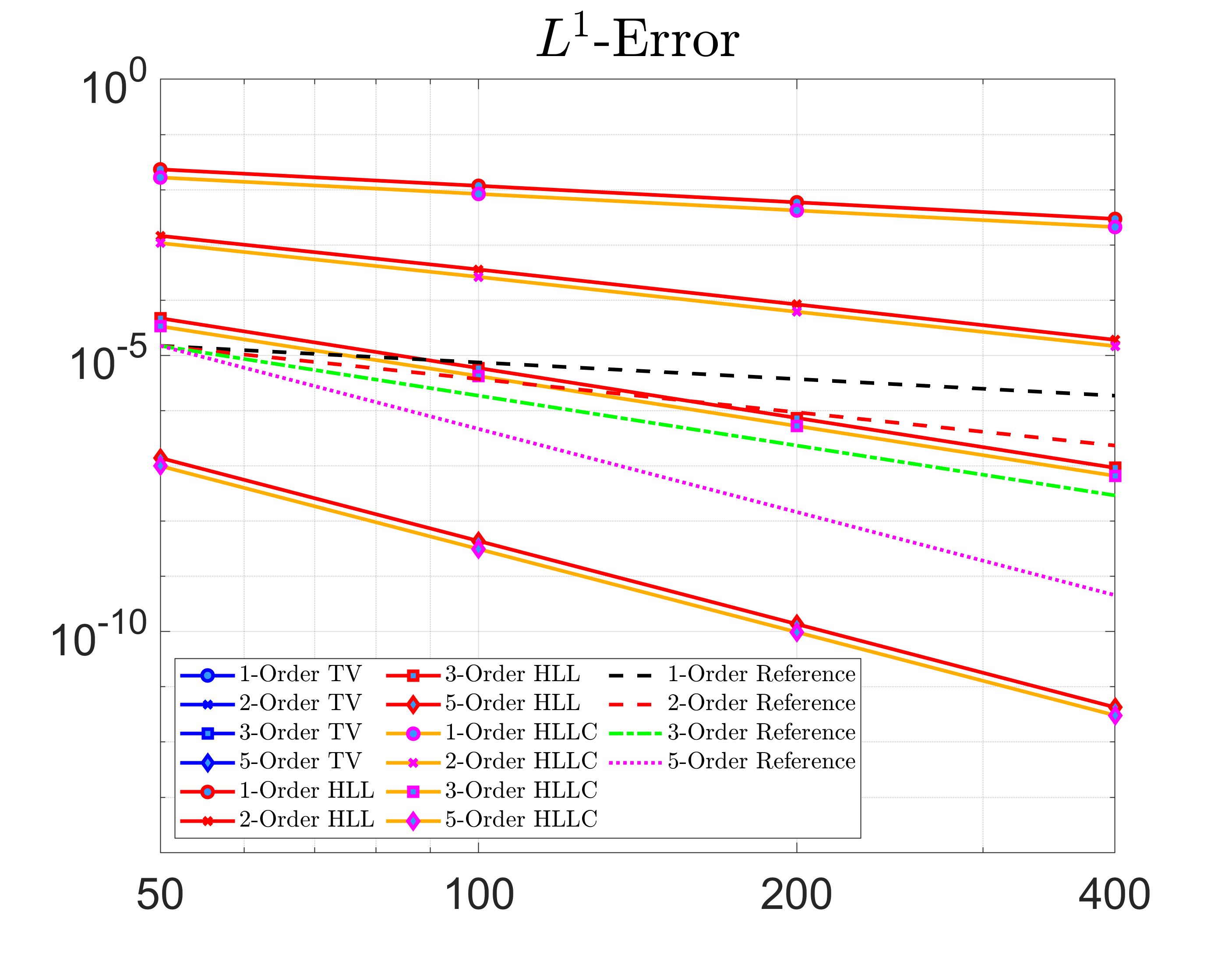}\hspace{1.cm}
		        \includegraphics[trim=1.cm 0.9cm 1.4cm 0.2cm, clip, width=6cm]{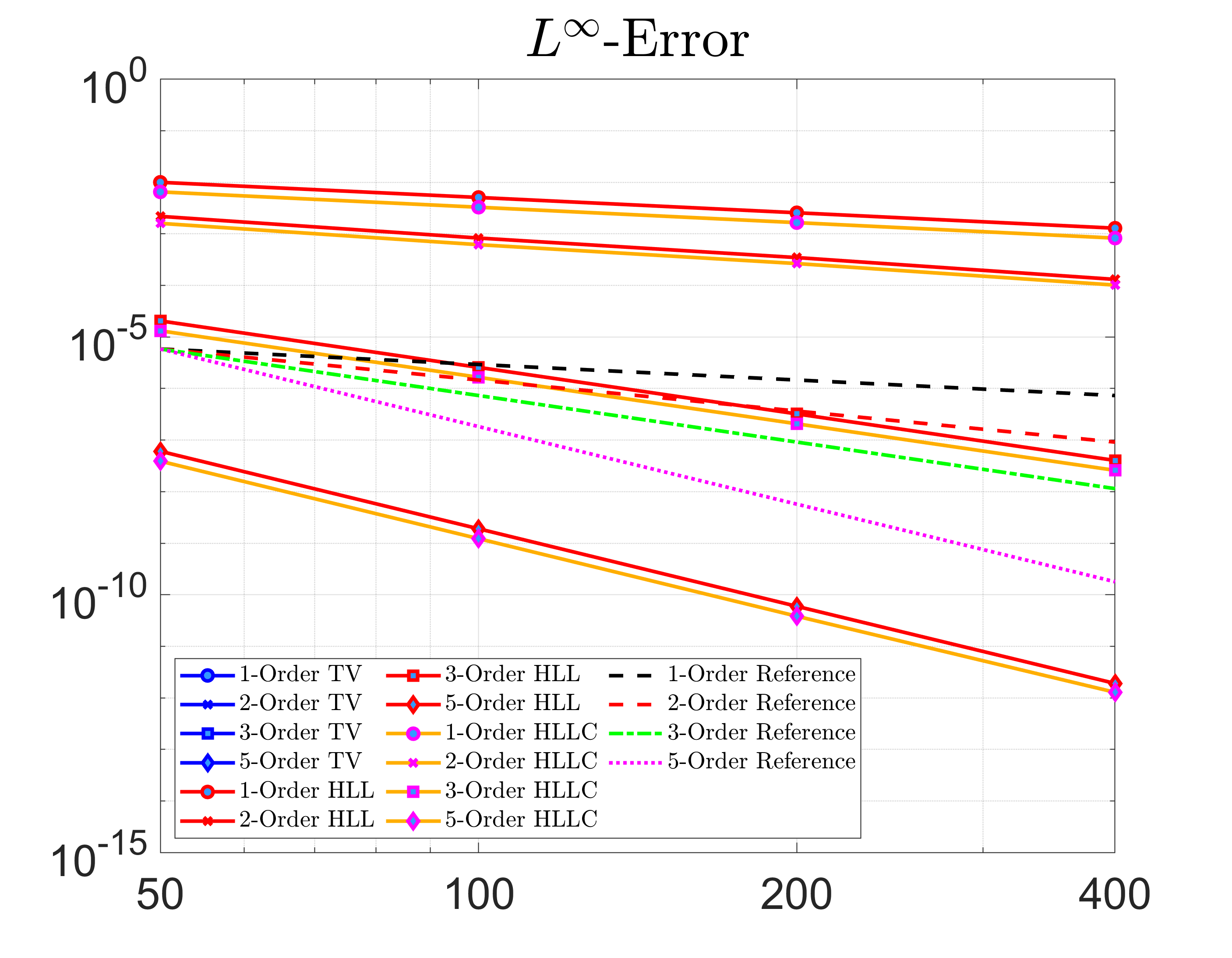}}
	\caption{\sf Example 6: $L^1$- (left) and $L^\infty$- (right) errors computed by the 1-, 2-, 3-, and 5-Orders TV splitting, CU (HLL), and HLLC schemes.}\label{fig44a}
\end{figure}

We also compute the numerical results using the CU (HLL) and HLLC schemes on the same meshes and present the obtained $L^1$- and $L^\infty$-errors in Figure~\ref{fig44a}. As observed, all the studied schemes obtain the expected order of accuracy. Moreover, the TV-splitting schemes demonstrate accuracy comparable to that of the HLLC schemes and exhibit higher accuracy than the corresponding CU (HLL) schemes.

\subsubsection*{Example 7---2-D Vortex Evolution Problem}
In this example taken from \cite{Titarev2005}; see also \cite{HS1999}, we consider the 2-D vortex evolution problem with the following initial conditions
\begin{equation*}
(\rho(x,y,0),u(x,y,0),v(x,y,0),p(x,y,0))=\left(T^{\frac{1}{\gamma-1}}, 1-\frac{\varepsilon}{2 \pi} e^{\frac{1}{2}(1-r^2)}y, 1+\frac{\varepsilon}{2 \pi} e^{\frac{1}{2}(1-r^2)}x, \rho^\gamma \right),
\end{equation*}
where $T=1-\frac{(\gamma-1)\varepsilon^2}{8\gamma \pi^2}e^{(1-r^2)}$,  $r^2=x^2+y^2$, and $\varepsilon=5$ is the vortex strength. The initial data, which is prescribed in the computational domain $[-5, 5]\times [-5,5]$ subject to the periodic boundary conditions, corresponds to a smooth vortex placed at the origin and is defined as the isentropic perturbation to the uniform flow of unit values of primitive variables and the exact solution is a vortex moving with a constant velocity at $45^{\circ}$ to the Cartesian mesh lines.  

\begin{figure}[ht!]
	\centerline{\includegraphics[trim=1.3cm 0.9cm 1.3cm 0.2cm, clip, width=6cm]{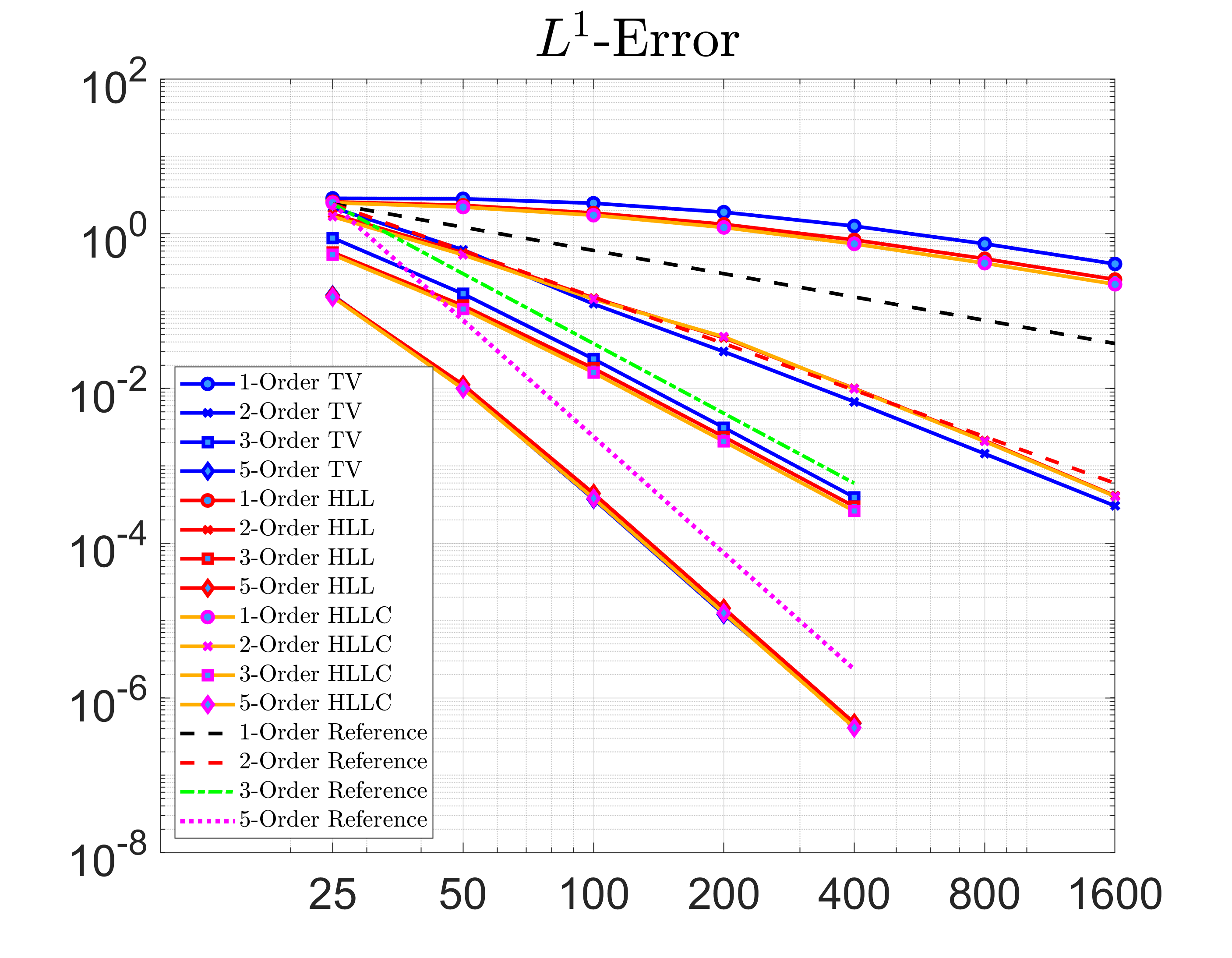}\hspace{1.cm}
		        \includegraphics[trim=1.3cm 0.9cm 1.3cm 0.2cm, clip, width=6cm]{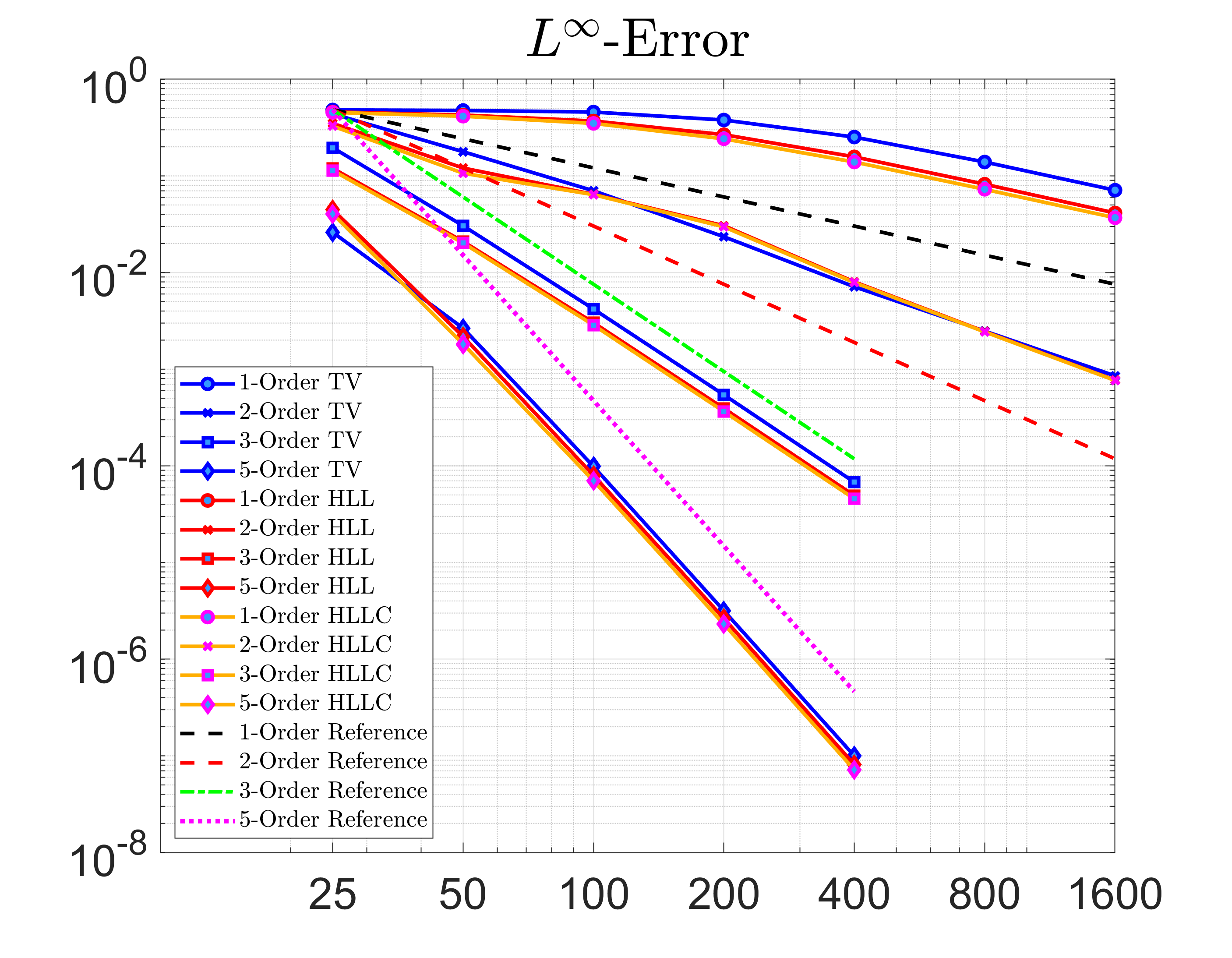}}
	\caption{\sf Example 7: $L^1$- (left) and $L^\infty$- (right) errors computed by the 1-, 2-, 3-, and 5-Orders TV splitting, HLL, and HLLC schemes.}\label{fig55a}
\end{figure}

We compute the numerical solution until the final time $t=10$ using the 1-, 2-, 3-, and 5-Order TV splitting schemes on a sequence of uniform meshes: $25 \times 25$,  $50\times 50$, $100\times 100$, $200\times 200$, and $400\times 400$, and then measure the $L^1$- and $L^\infty$-errors between the computed solutions and the exact solutions and the corresponding experimental convergence rates for the density. We also compare the performances of the TV splitting schemes with the CU (HLL) and HLLC schemes and present the obtained results in Tables \ref{tab52a} and \ref{tab52b}. One can see that for the 1- and 2-Order schemes, both first- and second-order convergence rates are observed only after significant mesh refinement, while the 3- and 5-Order TV splitting, CU (HLL), and HLLC schemes achieve the expected order of accuracy. In order to have a better view, we also show the $L^1$- and $L^\infty$-errors in Figure \ref{fig55a}, where we also show the results computed by the 1- and 2-Order schemes on finer meshes $800\times 800$ and $1600\times 1600$ to show that the 1- and 2-Order schemes achieve expected convergence rates after mesh refinement. Note that the $L^1$- and $L^\infty$-errors computed by the TV splitting schemes are slightly larger than the ones computed by the CU (HLL) and HLLC schemes, while HLLC schemes are the most accurate ones among the studied schemes. However, it can be observed that the TV-splitting schemes achieve accuracy comparable to that of the HLLC schemes and are more accurate than the corresponding CU (HLL) schemes; see Examples 4, 6, and 9. 

\begin{table}[ht!]
\centering
\resizebox{\linewidth}{!}{$
\begin{tabular}{|c|c|cc|cc|cc|cc|cc|cc|}
\hline
\multirow{2}{5em}{Method}&\multirow{2}{3em}{Mesh}&\multicolumn{2}{c|}{1-Order}&\multicolumn{2}{c|}{2-Order}&\multicolumn{2}{c|}{3-Order}&\multicolumn{2}{c|}{5-Order}\\
\cline{3-10}& &Error&Rate&Error&Rate&Error&Rate&Error&Rate\\
\hline
 \multirow{4}{5em}{TV \\[1.ex]splitting}& $25 \times 25$ &2.89 &---     &2.19e-00 &---  &8.89e-01&--- &1.62e-01 &---\\
                                        & $50 \times 50$ &2.84 &2.09e-02 &6.24e-01 &1.81 &1.68e-01&2.40&1.01e-02 &4.00\\
                                        & $100\times 100$&2.49 &1.90e-01 &1.24e-01 &2.33 &2.41e-02&2.80&3.73e-04 &4.76\\
                                        & $200\times 200$&1.91 &3.84e-01 &3.00e-02 &2.05 &3.12e-03&2.95&1.20e-05 &4.96\\
                                        & $400\times 400$&1.26 &5.97e-01 &6.75e-03 &2.15 &3.93e-04&2.99&4.65e-07 &4.70\\
\hline
 \multirow{4}{5em}{CU (HLL)}& $25 \times 25$ &2.60   &---      &1.81e-00 &---  &5.76e-01&--- &1.58e-01 &---\\
                            & $50 \times 50$ &2.34   &1.52e-01 &5.69e-01 &1.67 &1.18e-01&2.29&1.12e-02 &3.83\\
                            & $100\times 100$&1.86   &3.30e-01 &1.48e-01 &1.94 &1.83e-02&2.69&4.42e-04 &4.66\\
                            & $200\times 200$&1.33   &4.84e-01 &4.54e-02 &1.71 &2.37e-03&2.95&1.44e-05 &4.94\\
                            & $400\times 400$&0.84   &6.71e-01 &1.00e-02 &2.18 &2.98e-04&2.99&4.73e-07 &4.93\\
\hline
 \multirow{4}{5em}{HLLC}& $25 \times 25$ &2.53    &---       &1.68e-00 &---  &5.39e-01&--- &1.53e-01 &---\\
                        & $50 \times 50$ &2.22    &1.87e-01  &5.36e-01 &1.64 &1.07e-01&2.33&9.90e-03 &3.95\\
                        & $100\times 100$&1.74    &3.49e-01  &1.44e-01 &1.90 &1.61e-02&2.74&3.82e-04 &4.70\\
                        & $200\times 200$&1.21    &5.26e-01  &4.67e-02 &1.62 &2.07e-03&2.95&1.24e-05 &4.94\\
                        & $400\times 400$&0.74    &7.04e-01  &9.88e-03 &2.24 &2.60e-04&2.99&4.11e-07 &4.92\\
\hline
\end{tabular}
$}
\caption{\sf Example 7: $L^1$-errors and experimental convergence rates for the density $\rho$ computed by the 1-, 2-, 3-, and 5-Order TV splitting, CU (HLL), and HLLC schemes.\label{tab52a}}
\end{table}

\vskip 15pt

\begin{table}[ht!]
\centering
\resizebox{\linewidth}{!}{$
\begin{tabular}{|c|c|cc|cc|cc|cc|cc|cc|}
\hline
\multirow{2}{5em}{Method}&\multirow{2}{3em}{Mesh}&\multicolumn{2}{c|}{1-Order}&\multicolumn{2}{c|}{2-Order}&\multicolumn{2}{c|}{3-Order}&\multicolumn{2}{c|}{5-Order}\\
\cline{3-10}& &Error&Rate&Error&Rate&Error&Rate&Error&Rate\\
\hline
 \multirow{4}{5em}{TV \\[1.ex]splitting}& $25 \times 25$ &4.81e-01 &---      &4.38e-01 &---  &1.96e-01&--- &2.61e-02 &---\\
                                        & $50 \times 50$ &4.75e-01 &1.80e-02 &1.78e-01 &1.30 &3.05e-02&2.68&2.66e-03 &3.30\\
                                        & $100\times 100$&4.58e-01 &5.44e-02 &7.02e-02 &1.34 &4.20e-03&2.86&9.86e-05 &4.75\\
                                        & $200\times 200$&3.79e-01 &2.74e-01 &2.35e-02 &1.58 &5.43e-04&2.95&3.17e-06 &4.96\\
                                        & $400\times 400$&2.52e-01 &5.86e-01 &7.13e-03 &1.72 &6.78e-05&3.00&9.99e-08 &4.99\\
\hline
 \multirow{4}{5em}{CU (HLL)}& $25 \times 25$ &4.60e-01 &---      &3.50e-01 &---      &1.20e-01&--- &4.48e-02 &---\\
                            & $50 \times 50$ &4.26e-01 &1.08e-01 &1.21e-01 &1.53     &2.10e-02&2.51&2.18e-03 &4.36\\
                            & $100\times 100$&3.70e-01 &2.05e-01 &6.44e-02 &0.91     &3.03e-03&2.79&7.99e-05 &4.77\\
                            & $200\times 200$&2.67e-01 &4.69e-01 &3.05e-02 &1.08     &3.93e-04&2.95&2.64e-06 &4.92\\
                            & $400\times 400$&1.58e-01 &7.61e-01 &8.02e-03 &1.93     &4.93e-05&2.99&8.14e-08 &5.02\\
\hline
 \multirow{4}{5em}{HLLC}& $25 \times 25$ &4.54e-01 &---      &3.27e-01 &---      &1.14e-01&--- &4.04e-02 &---\\
                        & $50 \times 50$ &4.15e-01 &1.30e-01 &1.07e-01 &1.62     &2.02e-02&2.49&1.82e-03 &4.47\\
                        & $100\times 100$&3.49e-01 &2.47e-01 &6.39e-02 &0.74     &2.86e-03&2.82&7.02e-05 &4.70\\
                        & $200\times 200$&2.42e-01 &5.30e-01 &2.97e-02 &1.11     &3.65e-04&2.97&2.31e-06 &4.93\\
                        & $400\times 400$&1.39e-01 &7.97e-01 &7.88e-03 &1.91     &4.57e-05&3.00&7.18e-08 &5.00\\
\hline
\end{tabular}
$}
\caption{\sf Example 7: $L^\infty$-errors and experimental convergence rates for the density $\rho$ computed by the 1-, 2-, 3-, and 5-Order TV splitting, CU (HLL), and HLLC schemes. \label{tab52b}}
\end{table}

\subsubsection*{Example 8---Explosion Problem}
In this example, we consider the explosion problem studied in \cite{Toro2009}. We take the following initial conditions,
\begin{equation*}
(\rho(x,y,0),u(x,y,0),v(x,y,0),p(x,y,0))=\begin{cases}
(1,0,0,1),&x^2+y^2<0.16,\\
(0.125,0,0,0.1),&\mbox{otherwise},
\end{cases}
\end{equation*}
prescribed in the computational domain $[-1,1]\times[-1,1]$. subject to free boundary conditions at all the four sides.

\begin{figure}[ht!]
\centerline{\includegraphics[trim=5.5cm 8.2cm 2.7cm 7.5cm, clip, width=14.7cm]{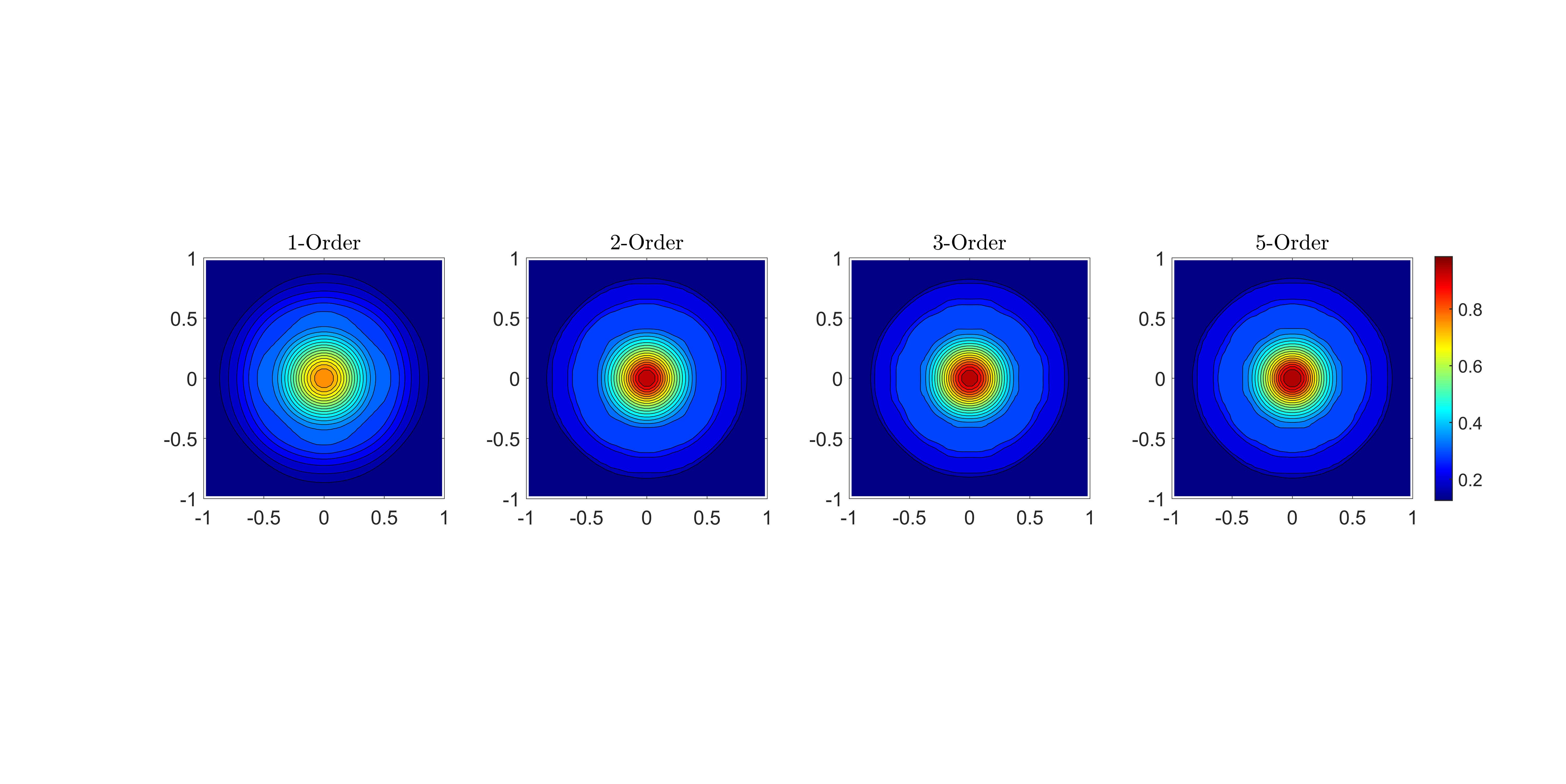}}
\caption{\sf{Example 8: Density $\rho$ computed by the 1-, 2-, 3-, and 5-Order schemes}\label{fig66}}
\end{figure}

\begin{figure}[ht!]
	\centerline{\includegraphics[trim=0.9cm 0.3cm 0.9cm 0.2cm, clip, width=6cm]{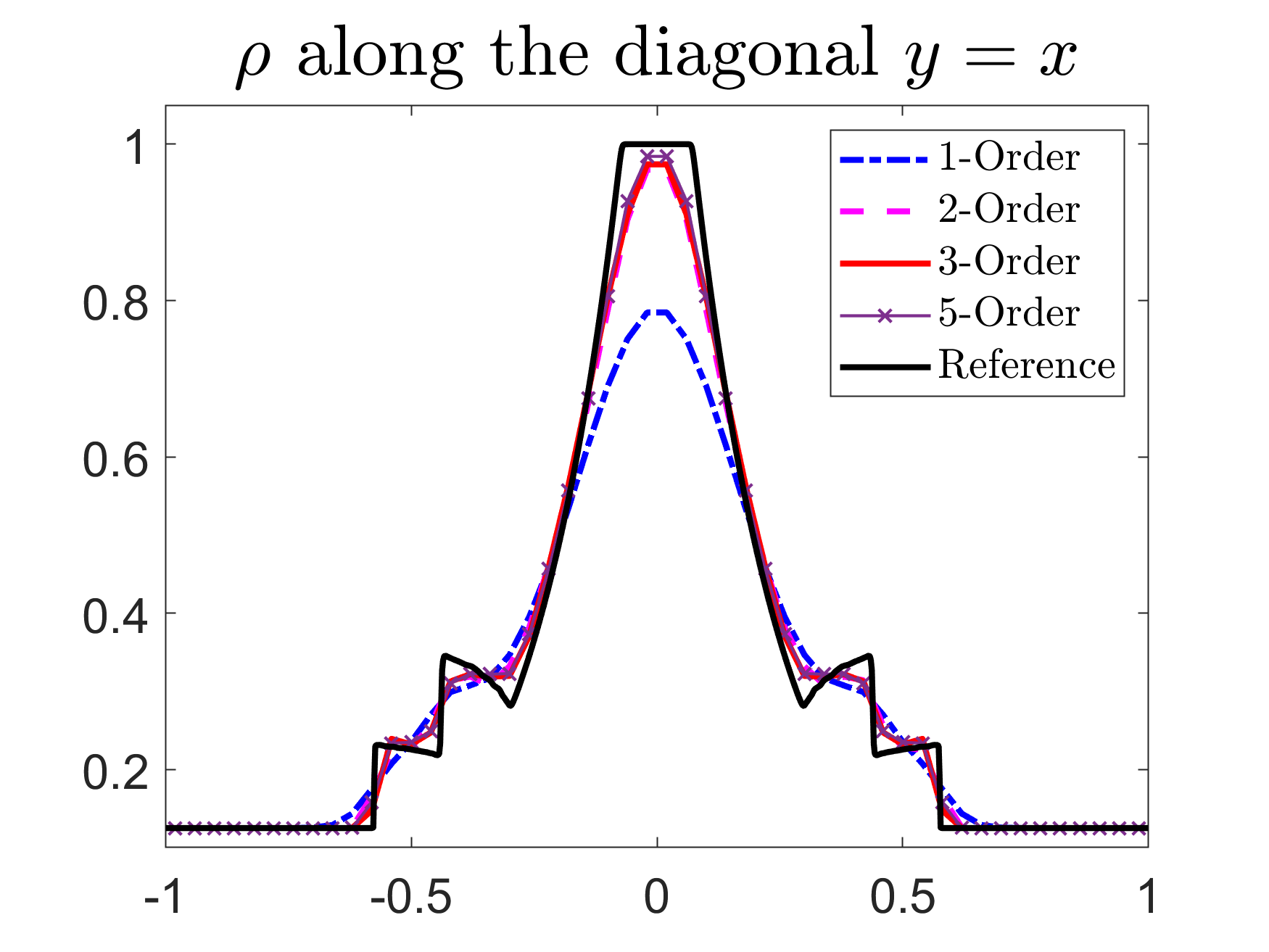}\hspace{0.5cm}
		        \includegraphics[trim=0.9cm 0.3cm 0.9cm 0.2cm, clip, width=6cm]{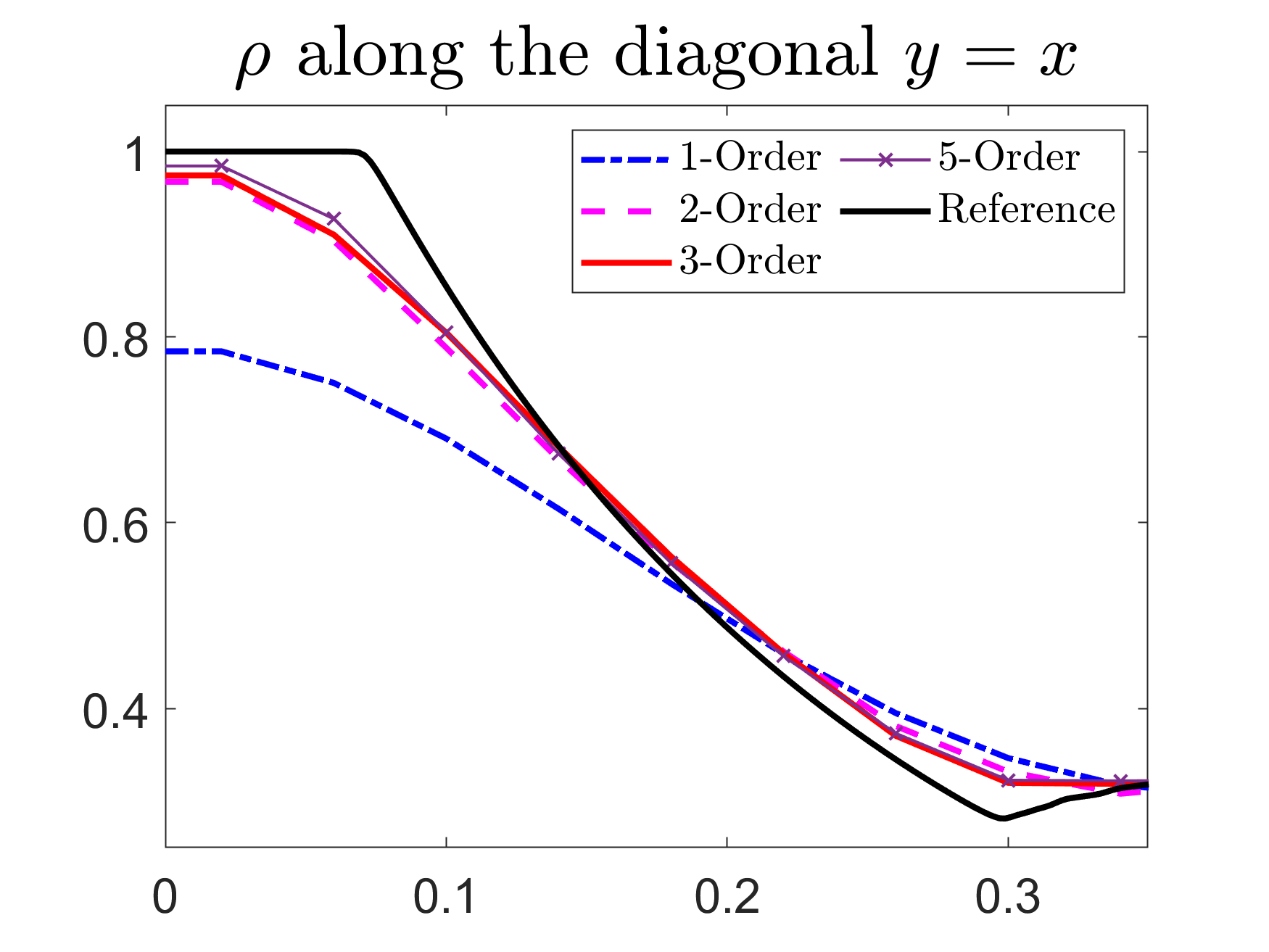}}
	\caption{\sf Example 8: Diagonal slices of the density $\rho$ computed by the 1-, 2-, 3-, and 5-Order schemes (left) and zoom at $x\in [0,0.35]$.}\label{fig6}
\end{figure}

We apply the studied 1-, 2-, 3-, and 5-Order schemes and compute the numerical solutions until the final time $t=0.25$ on a uniform mesh with $50 \times 50$ cells. The obtained results are presented in Figures \ref{fig66} and \ref{fig6}. In Figure \ref{fig6}, we show the slices of the densities along the diagonal $y = x$ for different schemes. As one can clearly see, using higher-order numerical schemes can achieve better resolution, with a particularly significant improvement observed when transiting from the 1-Order scheme to the 2-Order one.

\subsubsection*{Example 9---Implosion Problem}
In this example, we consider the implosion problem taken from \cite{Liska03} (see also \cite{CCHKL_22,Kurganov07,Liska03,CKX22}). The initial conditions,
\begin{equation*}
(\rho(x,y,0),u(x,y,0),v(x,y,0),p(x,y,0))=\begin{cases}
(0.125,0,0,0.14),&|x|+|y|<0.15,\\
(1,0,0,1),&\mbox{otherwise},
\end{cases}
\end{equation*}
are prescribed in the computational domain $[0,0.3]\times[0,0.3]$ with solid boundary conditions imposed at all four sides. This example was designed to test the amount of numerical diffusion present in different schemes as there is a jet forming near the origin and propagating along the diagonal $y=x$ direction, and schemes containing large numerical diffusion may not resolve the jet at all or the jet propagation velocity may be affected by the numerical diffusion.

\begin{figure}[ht!]
\centerline{\includegraphics[trim=5.5cm 8.2cm 2.7cm 7.5cm, clip, width=15.cm]{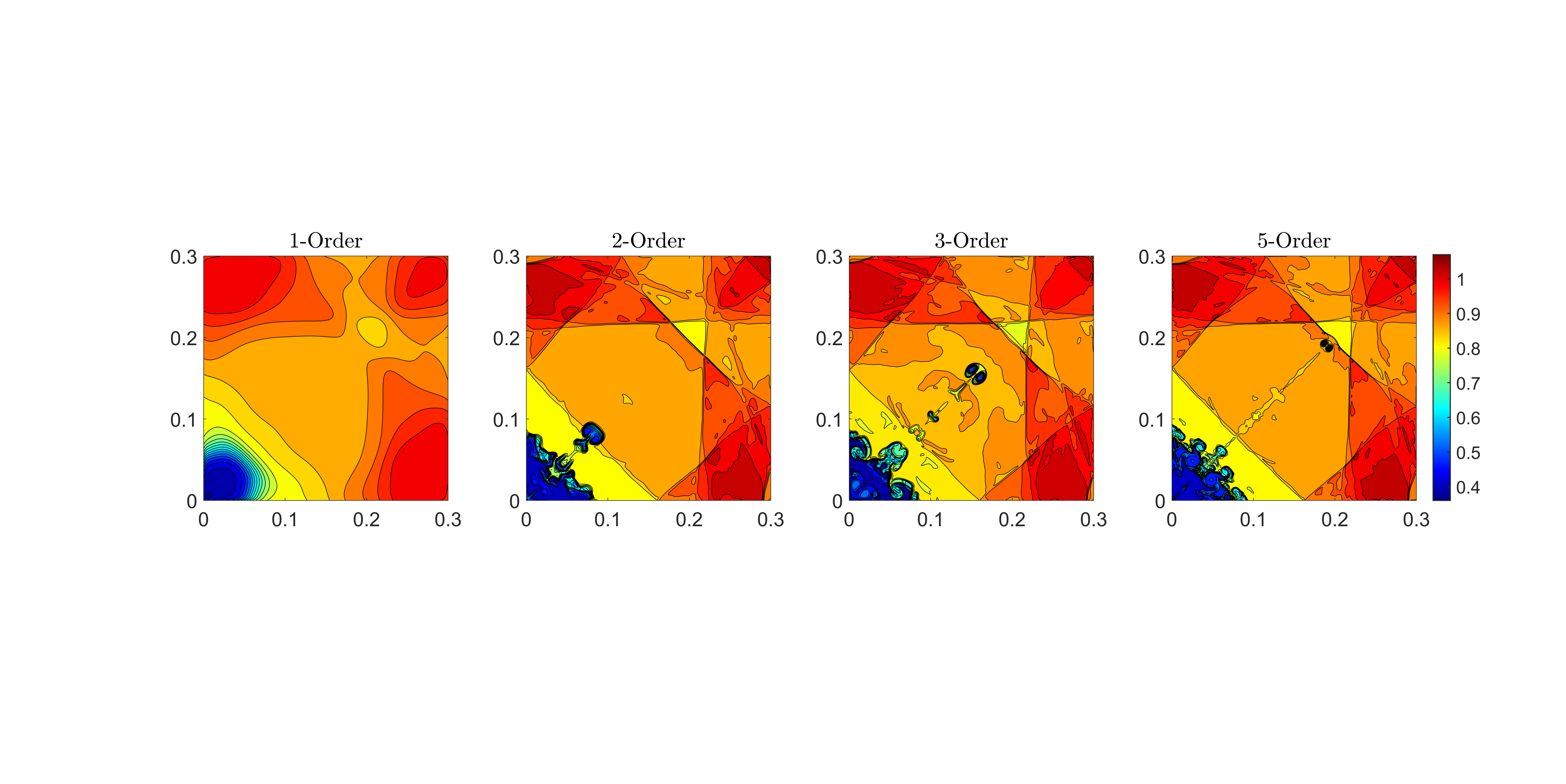}}
\caption{\sf{Example 9: Density ($\rho$) computed by the 1-, 2-, 3-, and 5-Order schemes}\label{fig7}}
\end{figure}

We compute the numerical solutions until the final time $t=2.5$ by the 1-, 2-, 3-, and 5-Order schemes on a uniform mesh of $400\times 400$ cells. The obtained results are depicted in Figure \ref{fig7}, where one can clearly observe that the jet propagates much further in the direction of $y=x$ when using higher-order schemes, clearly indicating that the high-order schemes are substantially less dissipative than the low-order ones.

\begin{figure}[ht!]
\centerline{\includegraphics[trim=5.5cm 8.2cm 2.7cm 7.5cm, clip, width=15.cm]{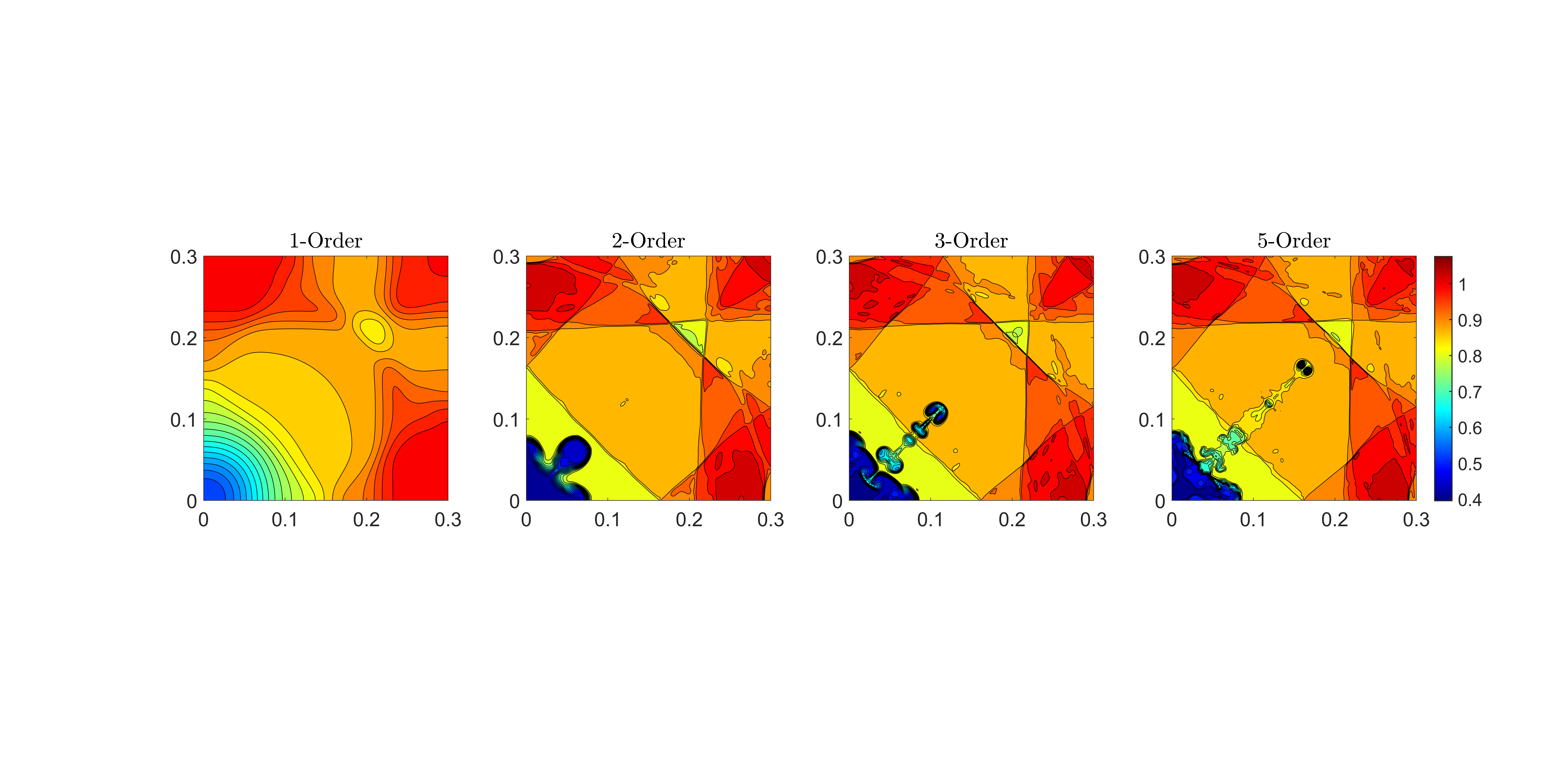}}
\vskip 12 pt
\centerline{\includegraphics[trim=5.5cm 8.2cm 2.7cm 7.5cm, clip, width=15.cm]{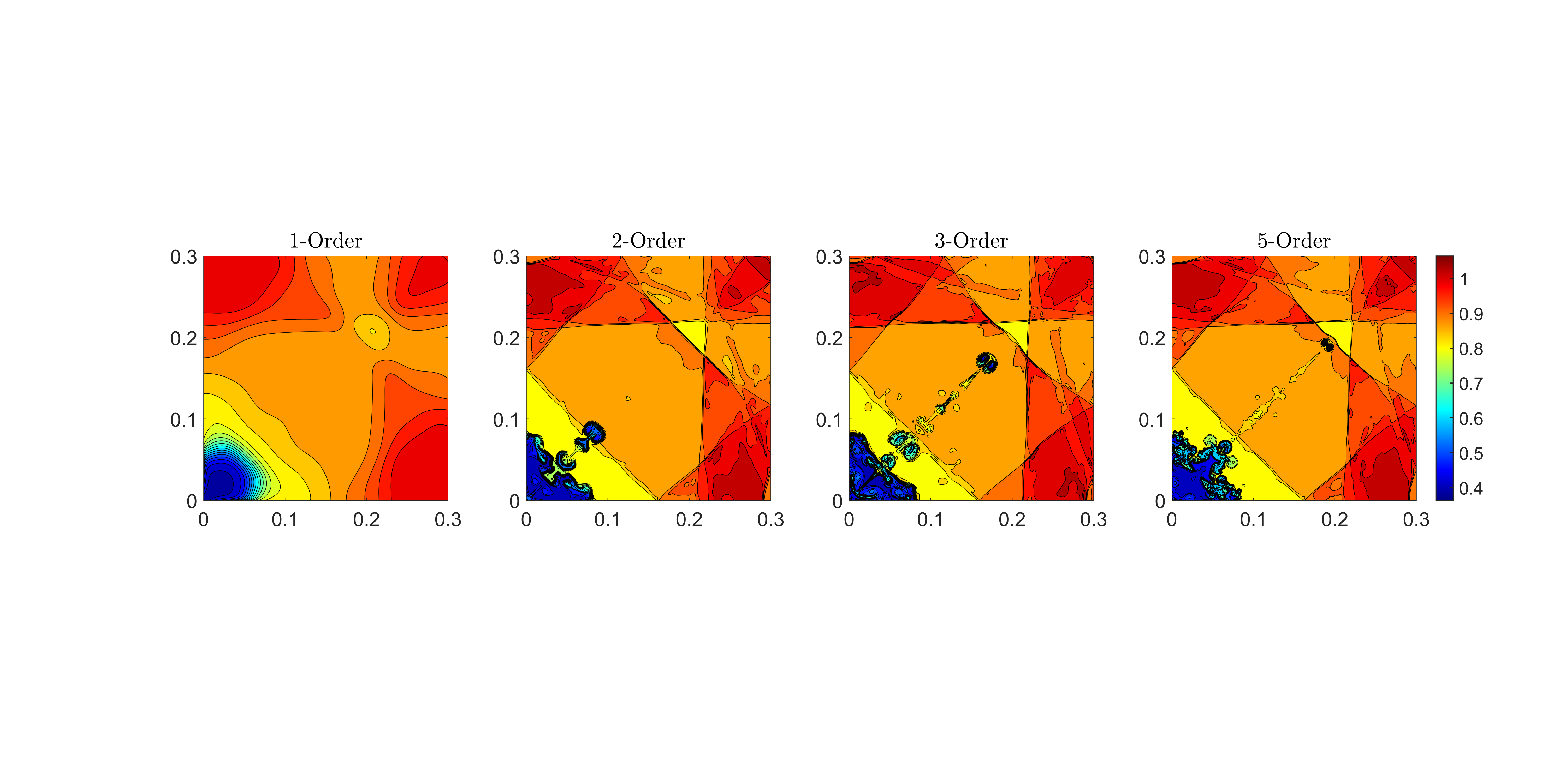}}
\caption{\sf{Example 9: Density $\rho$ computed by the 1-, 2-, 3-, and 5-Order CU (top row) and HLLC (bottom row) schemes}\label{fig7a}}
\end{figure}

In order to compare the studied TV splitting schemes with the CU, or HLL, and HLLC schemes, we compute the numerical results by the corresponding CU (HLL) and HLLC schemes and plot the corresponding results in Figure \ref{fig7a}. One can see that, the jets produced by the HLLC and TV splitting schemes move further than the corresponding CU (HLL) scheme. At the same time, the positions of the jets produced by the HLLC and TV splitting schemes are close, indicating that there is much less dissipation in the TV splitting schemes compared with the CU (HLL) ones.

\begin{figure}[ht!]
\centerline{\includegraphics[trim=5.5cm 6.8cm 2.6cm 6.cm, clip, width=13.cm]{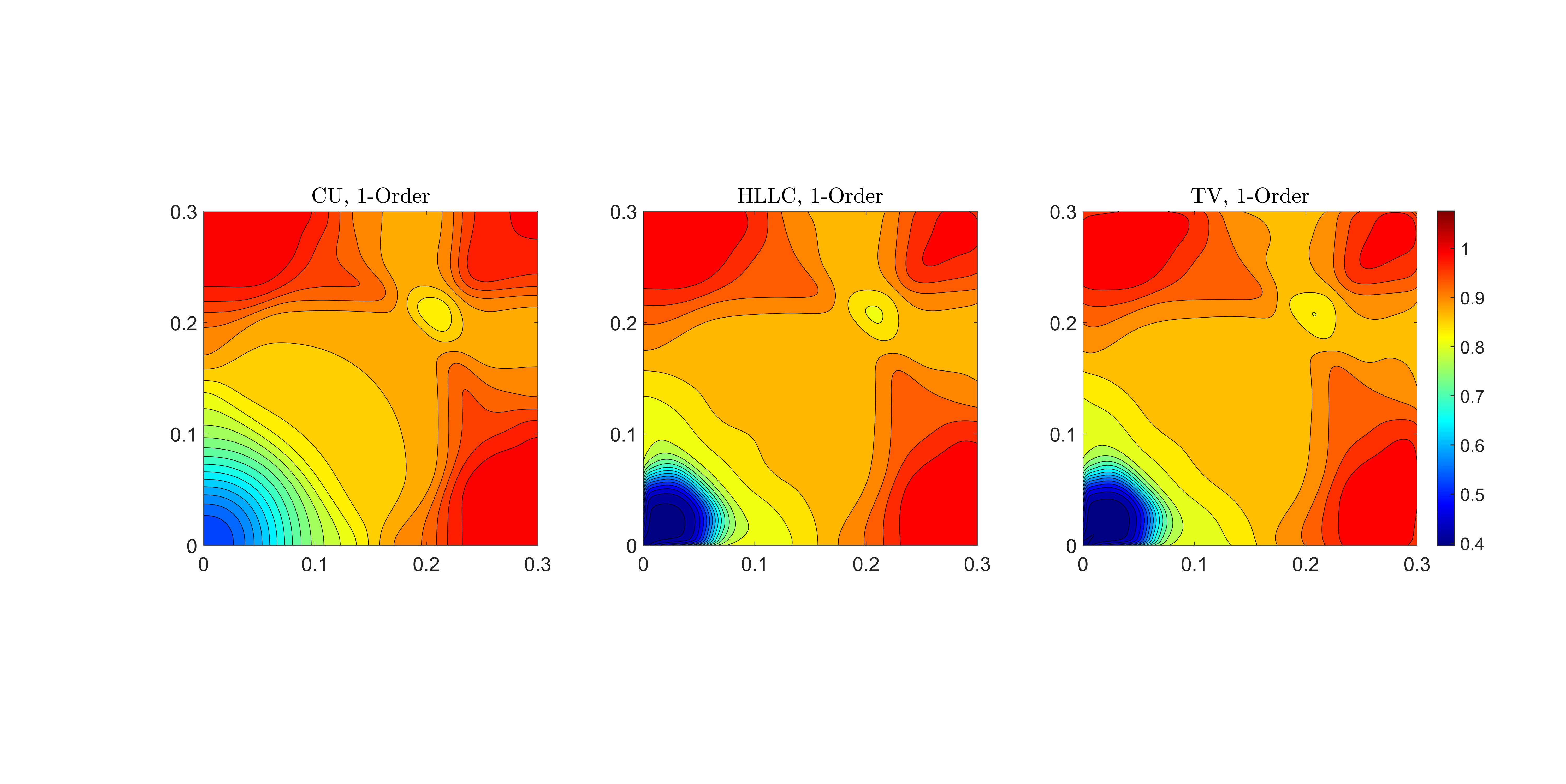}}
\vskip 12 pt
\centerline{\includegraphics[trim=5.5cm 6.8cm 2.6cm 6.cm, clip, width=13.cm]{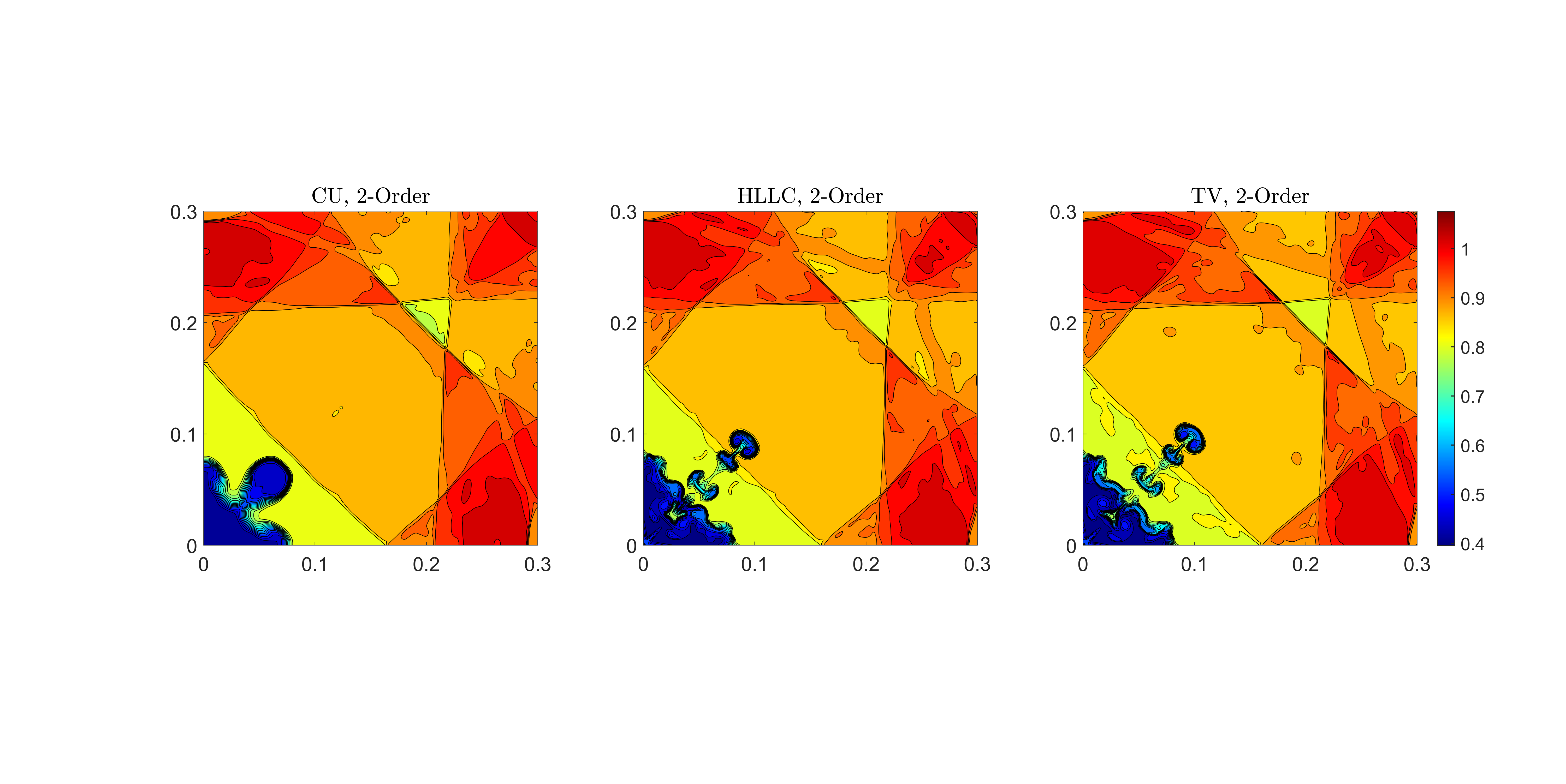}}
\vskip 12 pt
\centerline{\includegraphics[trim=5.5cm 6.8cm 2.6cm 6.cm, clip, width=13.cm]{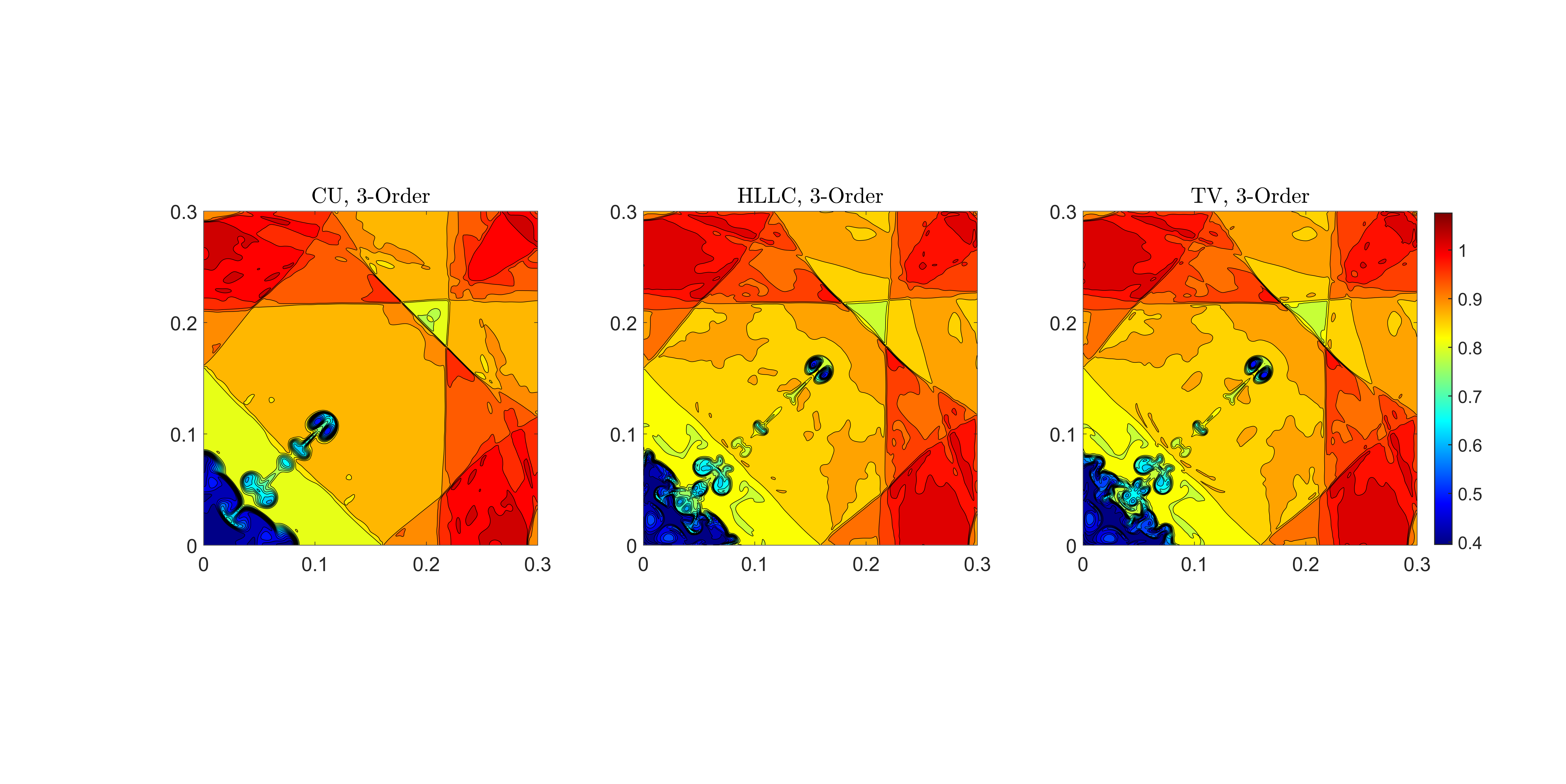}}
\vskip 12 pt
\centerline{\includegraphics[trim=5.5cm 6.8cm 2.6cm 6.cm, clip, width=13.cm]{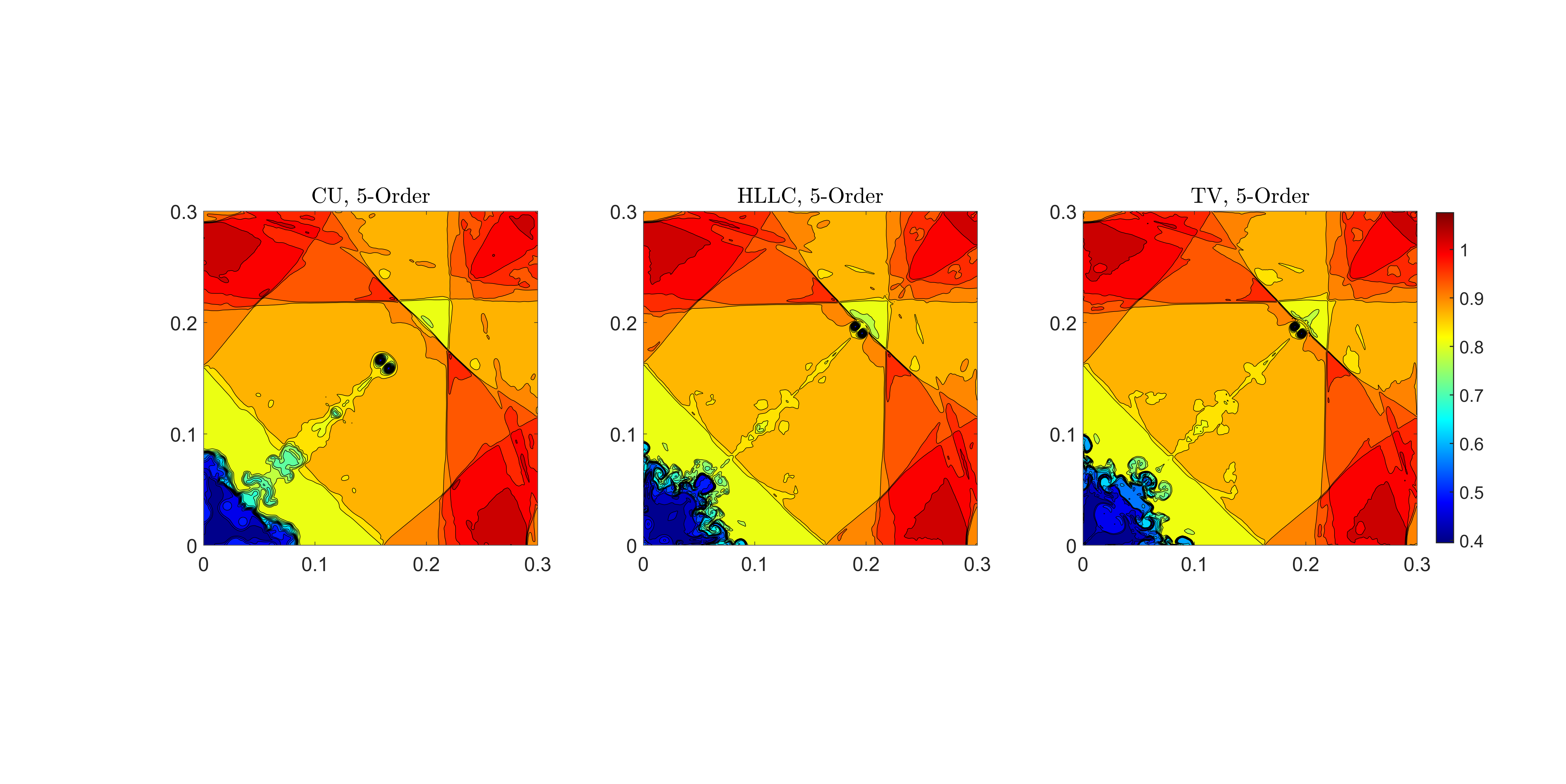}}
\caption{\sf{Example 9: Density $\rho$ computed by the CU (left column), HLLC (middle column), and TV splitting (right column) schemes.\label{fig7b}}}
\end{figure}

It is instructive to check whether the studied TV splitting schemes are more efficient than the CU (HLL) and HLLC schemes. To this end, we measure the CPU time
consumed during the CU (HLL) schemes and refine the mesh used by the HLLC and TV splitting scheme to the level that exactly the same CPU time is consumed to compute all three numerical solutions.  The corresponding meshes are $400 \times 400$ for the 1-, 2-, 3-, and 5-Order CU (HLL) schemes, $425\times 425$, $425\times 425$, $405\times 405$, and $408\times 408$ for the HLLC schemes, and $455\times 455$, $455\times 455$, $420\times 420$, and $418\times 418$ for the TV splitting scheme. The obtained numerical results, presented in Figure \ref{fig7b}, indicate that the TV splitting schemes still achieve a much higher resolution than the CU (HLL) schemes. At the same time, the TV splitting schemes are slightly less dissipative than the HLLC schemes; see, e.g., the 2-Order results. 

\subsubsection*{Example 10---Kelvin-Helmholtz (KH) Instability}
In this example, we study the KH instability taken from \cite{Fjordholm16,Panuelos20} (see also \cite{Feireisl21,CCHKL_22,CKX22}). We take the following initial data:
\begin{equation*}
\begin{aligned}
&(\rho(x,y,0),u(x,y,0))=\begin{cases}
(1,-0.5+0.5e^{(y+0.25)/L}),&y<-0.25,\\
(2,0.5-0.5e^{(-y-0.25)/L}),&-0.25<y<0,\\
(2,0.5-0.5e^{(y-0.25)/L}),&0<y<0.25,\\
(1,-0.5+0.5e^{(0.25-y)/L}),&y>0.25,
\end{cases}\\
&v(x,y,0)=0.01\sin(4\pi x),\qquad p(x,y,0)\equiv1.5,
\end{aligned}
\end{equation*}
where $L$ is a smoothing parameter (we take $L=0.00625$), which corresponds to a thin shear interface with a perturbed vertical velocity field $v$ in the conducted simulations. The periodic boundary conditions are imposed on all four sides of the computational domain $[-0.5,0.5]\times[-0.5,0.5]$.

\begin{figure}[ht!]
\centerline{
\includegraphics[trim=5.3cm 8.3cm 2.7cm 7.5cm, clip, width=14.7cm]{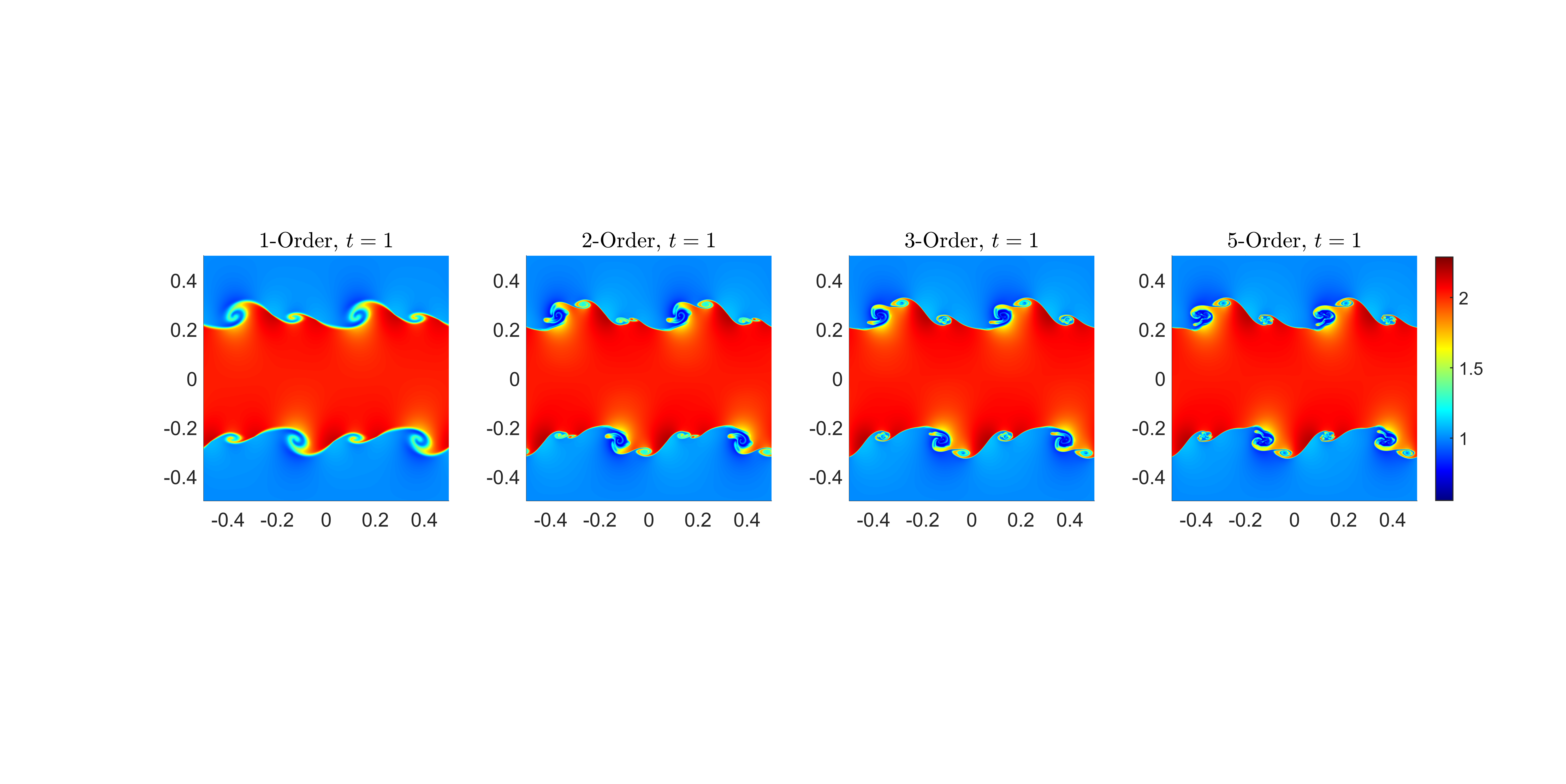}}
\vskip 16pt
\centerline{
\includegraphics[trim=5.3cm 8.3cm 2.7cm 7.5cm, clip, width=14.7cm]{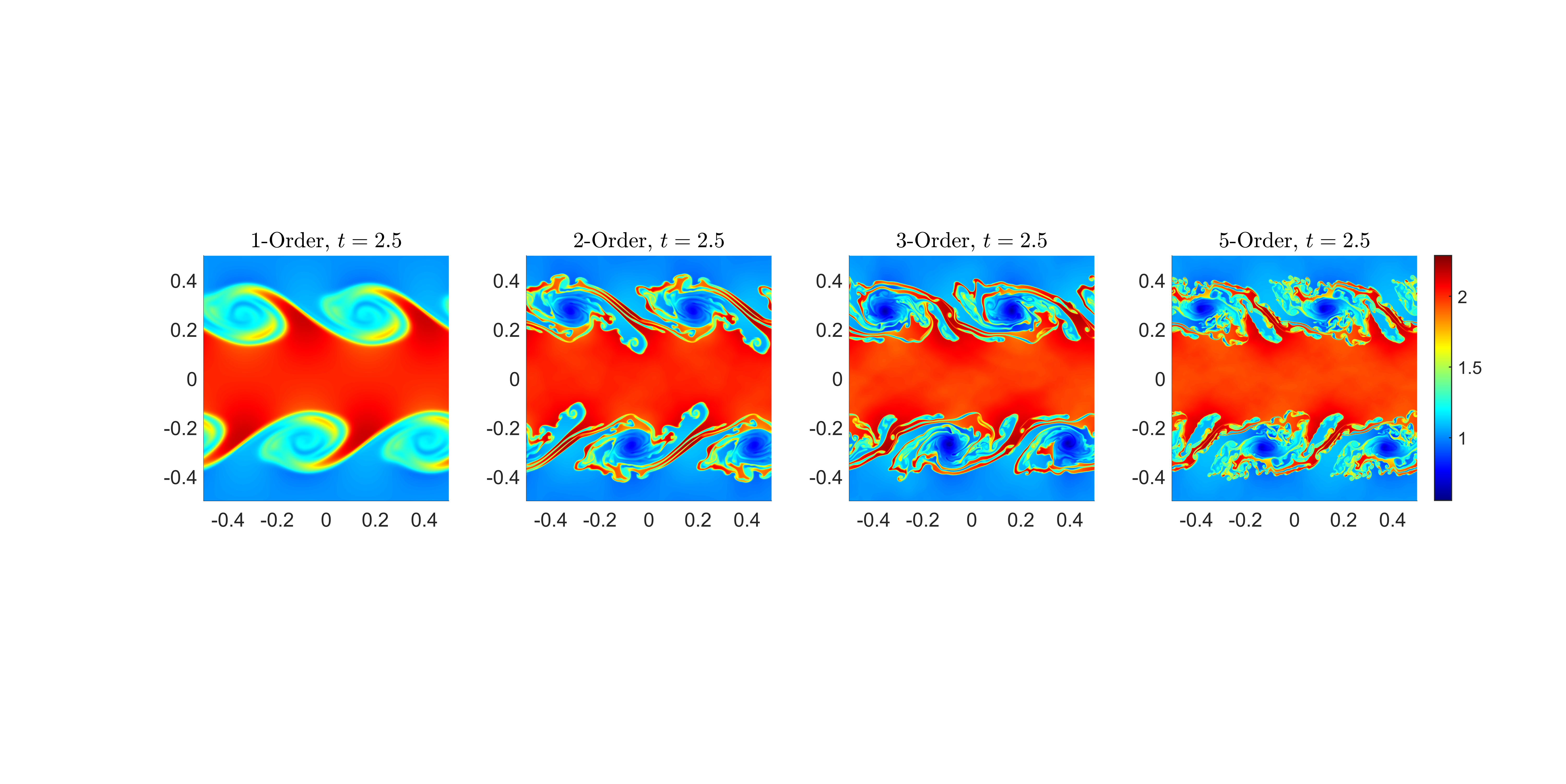}}
\vskip 16pt
\centerline{
\includegraphics[trim=5.3cm 8.3cm 2.7cm 7.5cm, clip, width=14.7cm]{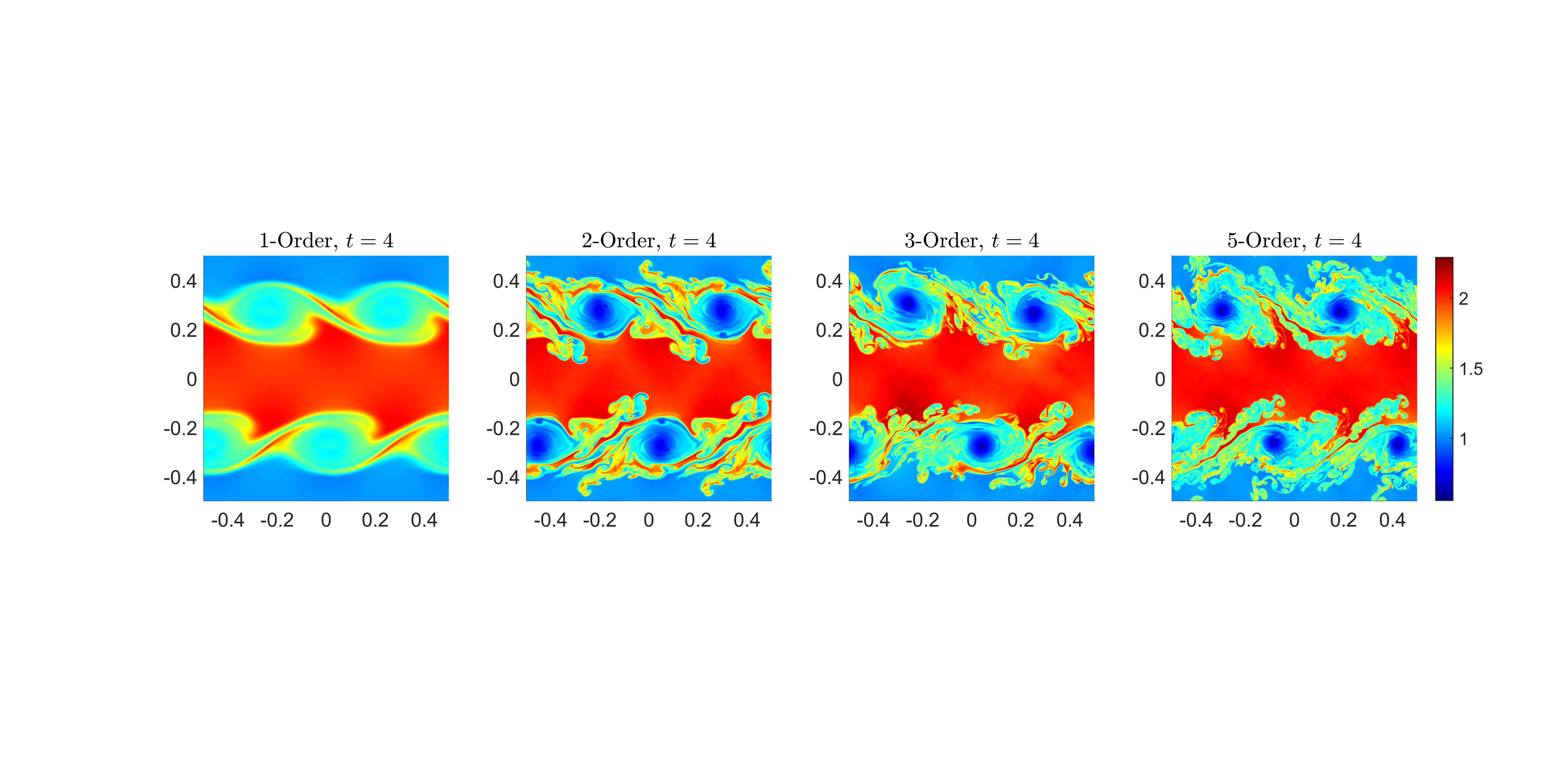}}
\caption{\sf Example 10: Time snapshots of density $\rho$ computed by the 1- (first column), 2- (second column), 3- (third column), and 5- (fourth column) Order schemes at $t=1$ (top row), $t=2.5$ (middle row), and $t=4$ (bottom row).\label{fig8}}
\end{figure}

We compute the numerical solutions until the final time $t=4$ by the 1-, 2-, 3-, and 5-Order schemes on a uniform mesh of $1024\times 1024$ cells, and plot the numerical results at times $t=1$, 2.5, and 4 in Figure \ref{fig8}. One can observe that at the early time $t=1$, the vortex streets generated by the high-order schemes are more pronounced. These structures grow exponentially over time, leading to increasingly complex turbulent mixing, particularly evident at later times $t=2.5$ and 4, clearly indicating that high-order schemes exhibit significantly less dissipation compared to low-order ones.

\subsubsection*{Example 11---Rayleigh-Taylor (RT) Instability}
In the last example taken from \cite{Shi03} (see also \cite{Wang20,CCHKL_22,CKX22,CKM25}), we investigate the RT instability, which is a physical phenomenon
occurring when a layer of heavier fluid is placed on top of a layer of lighter fluid. The model is governed by the 2-D Euler equations
\eref{3.1}--\eref{3.2} with added gravitational source terms and the modified system reads as
$$
\begin{aligned}
&\rho_t+(\rho u)_x+(\rho v)_y=0,\\
&(\rho u)_t+(\rho u^2 +p)_x+(\rho uv)_y=0,\\
&(\rho v)_t+(\rho uv)_x+(\rho v^2+p)_y=\rho,\\
&E_t+[u(E+p)]_x+[v(E+p)]_y=\rho v.
\end{aligned}
$$
We consider the following initial conditions:
\begin{equation*}
(\rho(x,y,0),u(x,y,0),v(x,y,0),p(x,y,0))=\begin{cases}
(2,0,-0.025\,c\cos(8\pi x),2y+1),&y<0.5,\\
(1,0,-0.025\,c\cos(8\pi x),y+1.5),&\mbox{otherwise},
\end{cases}
\end{equation*}
where $c:=\sqrt{\gamma p/\rho}$ is the speed of sound, prescribed in the computational domain $[0, 0.25]\times [0,1]$ with the solid wall boundary conditions imposed at $x=0$ and $x=0.25$, and the following Dirichlet boundary conditions imposed at the top and bottom boundaries:
$$
(\rho,u,v,p)|_{y=1}=(1,0,0,2.5),\qquad(\rho,u,v,p)|_{y=0}=(2,0,0,1).
$$

\begin{figure}[ht!]
\centerline{\includegraphics[trim=5.2cm 4.5cm 2.8cm 3.5cm, clip, width=14.7cm]{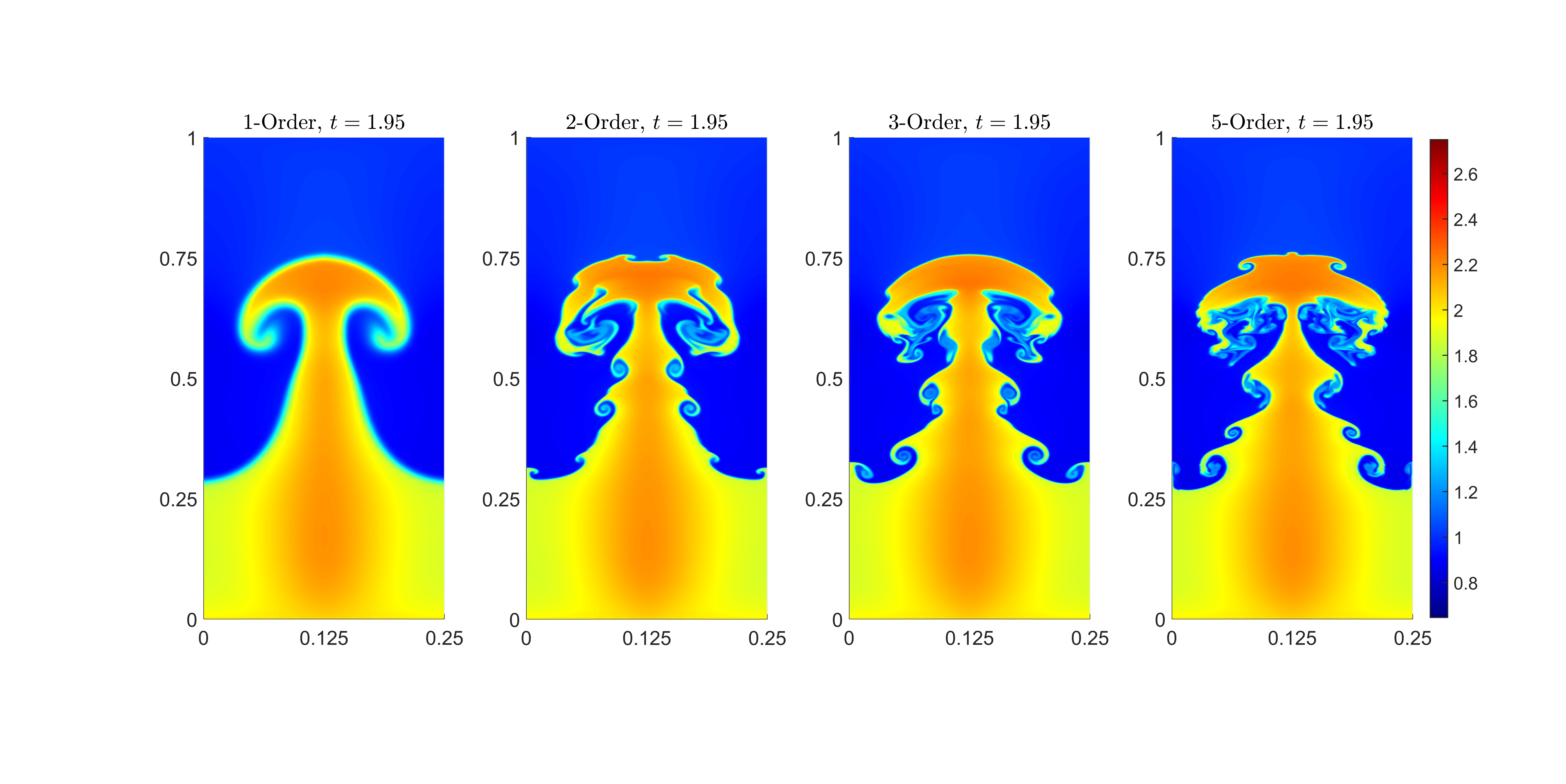}}
\vskip 12pt 
\centerline{\includegraphics[trim=5.2cm 4.5cm 2.8cm 3.5cm, clip, width=14.7cm]{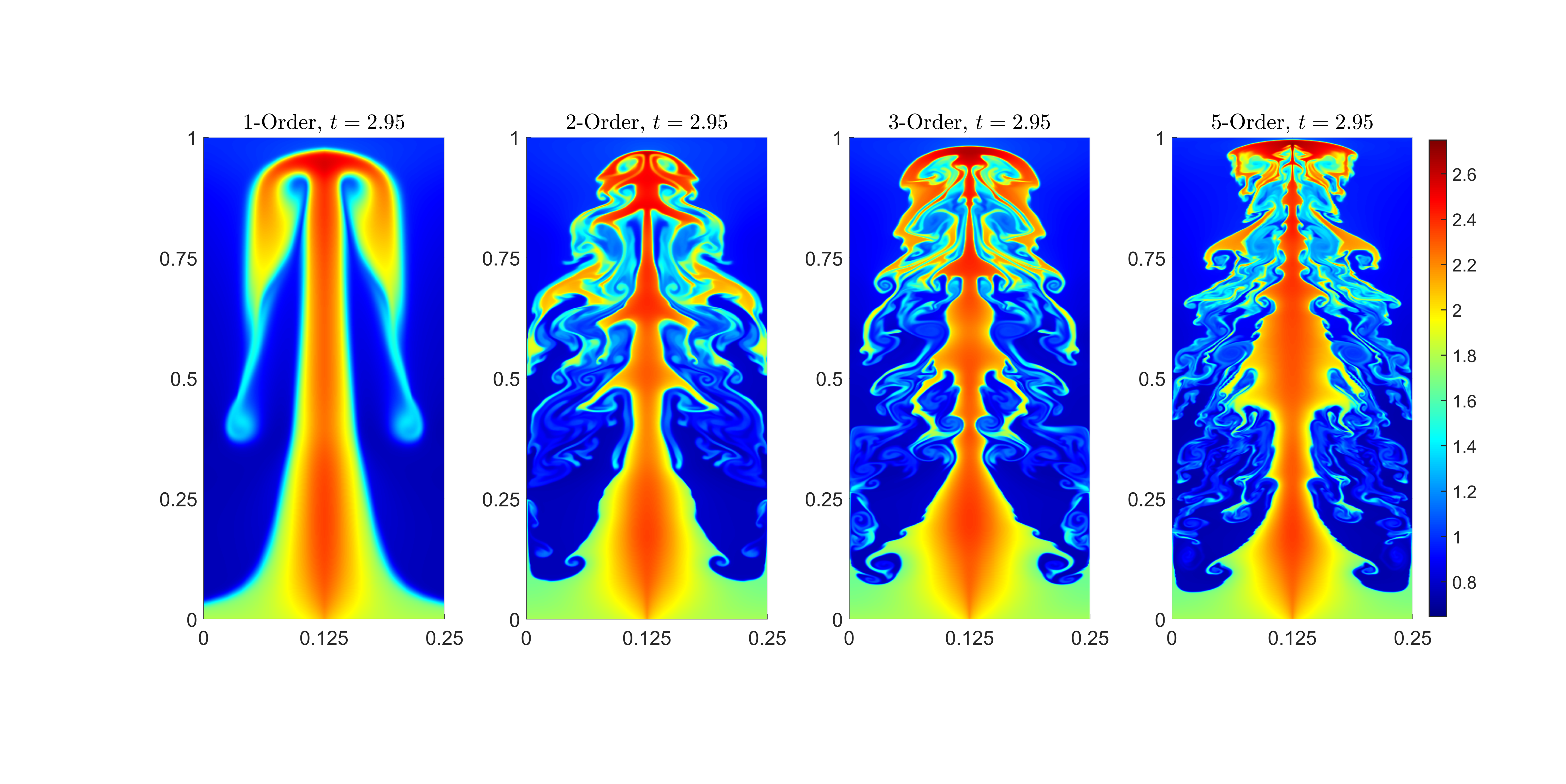}}
\caption{\sf Example 11: Density $\rho$ computed by the 1- (first column), 2- (second column), 3- (third column), and 5- (fourth column) Order schemes at $t=1.95$ (top row) and $t=2.95$ (bottom row).\label{fig9}}
\end{figure}

We compute the numerical solutions until the final time $t=2.95$ by 1-, 2-, 3-, and 5-Order schemes on the uniform mesh of $256\times 1024$ cells and present the numerical results at times $t=1.95$ and 2.95 in Figure \ref{fig9}. One can see that, there are pronounced differences between the solutions computed by different schemes and the structures captured by the high-order schemes are much more complicated, which again demonstrates that the high-order schemes can capture more details and are less dissipative than the low-order ones.

\subsection*{Acknowledgments}
The work of S. Chu was supported in part by the DFG (German Research Foundation) through HE5386/19-3, 27-1.
The work of M. Herty was funded by the DFG--SPP 2183: Eigenschaftsgeregelte Umformprozesse with the Project(s) HE5386/19-2,19-3 Entwicklung eines flexiblen isothermen Reckschmiedeprozesses f\"ur die eigenschaftsgeregelte Herstellung von Turbinenschaufeln aus Hochtemperaturwerkstoffen (424334423) and by the Deutsche
Forschungsgemeinschaft (DFG, German Research Foundation)--SPP 2410 Hyperbolic Balance Laws in Fluid Mechanics: Complexity, Scales,
Randomness (CoScaRa) within the Project(s) HE5386/26-1 (Numerische Verfahren f\"ur gekoppelte Mehrskalenprobleme,525842915) and
(Zuf\"allige kompressible Euler Gleichungen: Numerik und ihre Analysis, 525853336) HE5386/27-1. E. F. Toro acknowledges partial support from project number P130625149, entitled Development of a Numerical Multiphase Flow Tool for Applications to Petroleum Industry, funded by Repsol S.A. Spain.

\appendix
\section{One-Dimensional Local Characteristic Decomposition Based Third-Order WENO-Type Interpolant}\label{appa}
In this appendix, we briefly describe the 1-D local characteristic decomposition (LCD) based third-order WENO-Type interpolant.

Supposing the point values $\mU_j$ are given at uniform grid points $x=x_j$, we now demonstrate how to compute the 
interpolated left-sided value of $\mU$ at $x=x_\jph$, denoted as $\mU^-_\jph$. The corresponding right-sided value $\mU^+_\jph$, can be computed 
in a mirror-symmetric manner.

The value of $(U^{(i)}_\jph)^-$ is computed using a weighted average of the two linear interpolants $\bm{{\cal P}}_0(x)$ and $\bm{{\cal P}}_1(x)$, which are obtained using the stencils $[x_{j-1},x_j]$ and $[x_j,x_{j+1}]$, respectively:
\begin{equation}
(U^{(i)}_\jph)^-=\omega^{(i)}_0{\cal P}^{(i)}_0(x_\jph)+\omega^{(i)}_1{\cal P}^{(i)}_1(x_\jph),
\label{A1}
\end{equation}
where
\begin{equation}\label{A2}
{\cal P}^{(i)}_0(x_\jph)=-\frac{1}{2}U^{(i)}_{j-1}+\frac{3}{2}U^{(i)}_j,\quad {\rm and} \quad {\cal P}^{(i)}_1(x_\jph)= \frac{1}{2}U^{(i)}_j+\frac{1}{2}U^{(i)}_{j+1}.
\end{equation}
By performing a straightforward Taylor expansion, one can demonstrate that \eref{A1}--\eref{A2}  achieve third-order accuracy if the weights $\omega_k$ in \eref{A1} are chosen as:
\begin{equation}
\omega^{(i)}_k=\frac{\alpha^{(i)}_k}{\alpha^{(i)}_0+\alpha^{(i)}_1},\quad
\alpha^{(i)}_k=d_k\left[1+\frac{\tau^{(i)}_3}{\beta^{(i)}_k+\varepsilon}\right],
\label{A3}
\end{equation}
where $d_0=\frac{1}{4}$ and $d_1=\frac{3}{4}$, 
$$
\tau^{(i)}_3=\big|\beta^{(i)}_2-\beta^{(i)}_3\big|^p,
$$
and  the smoothness indicators $\beta^{(i)}_k$ are given by (see e.g., \cite{CH_third})
\begin{equation}
\beta^{(i)}_0=(U^{(i)}_{j-1}-U^{(i)}_j)^2, \quad \beta^{(i)}_1=(U^{(i)}_j-U^{(i)}_{j+1})^2, \quad 
\label{A4a}
\end{equation}
and
\begin{equation}\label{A4}
\begin{aligned}
 \beta_2&=\frac{13}{12}(U^{(i)}_{j-1}-2U^{(i)}_{j}+U^{(i)}_{j+1})^2+\frac{1}{4}(U^{(i)}_{j+1}-U^{(i)}_{j-1})^2,\\ \beta_3&=\frac{13}{12}(U^{(i)}_{j}-2U^{(i)}_{j+1}+U^{(i)}_{j+2})^2+\frac{1}{4}(3U^{(i)}_{j}-4U^{(i)}_{j+1}+U^{(i)}_{j+2})^2.
\end{aligned}
\end{equation}
Finally, in all of the numerical examples reported in \S\ref{sec4}, we have chosen $p=1.4$ and $\varepsilon=10^{-12}$.

Although the third-order interpolant \eref{A1}--\eref{A4} is essentially non-oscillatory, it is well-known that applying it to the conservative variables $\mU$ in a componentwise manner can result in spurious oscillations in the computed solution. To address this, we adopt the reconstruction procedure within the LCD framework.

To this end, we first introduce the matrix $\widehat A_\jph:=A(\widehat\mU_\jph)$, where $A=\frac{\partial \mF}{\partial \mU}$ and $\widehat\mU_\jph$ is either a simple average $(\mU_j+\mU_{j+1})/2$ or another type of average of the $\mU_j$ and $\mU_{j+1}$ states (in the numerical examples reported in \S\ref{sec4}, we have used the simple average). We then compute the matrices
$R_\jph$ and $R^{-1}_\jph$ such that $R^{-1}_\jph\widehat A_\jph R_\jph$ is a diagonal matrix and introduce the local characteristic
variables
in the neighborhood of $x=x_\jph$:
$$
\bm\Gamma_m=R^{-1}_\jph\mU_m,\quad m=j-1,\ldots,j+2.
$$
Equipped with the values $\bm\Gamma_{j-1}$,  $\bm\Gamma_j$, $\bm\Gamma_{j+1}$, and $\bm\Gamma_{j+2}$,  we apply the interpolation procedure \eref{A1}--\eref{A4} to each of the components $\Gamma^{(i)}$, $i=1,\ldots,d$ of $\bm\Gamma$ to obtain ${\bm\Gamma}_\jph^-$ (the values of ${\bm\Gamma}_\jph^+$ are computed, as mentioned above, in the mirror-symmetric way). Finally, the corresponding point values of $\mU$ are given by
\begin{equation*}
\mU^\pm_\jph=R_\jph{\bm\Gamma}^\pm_\jph.
\end{equation*}
\begin{remark}
A detailed explanation of how the average matrix $\widehat A_\jph$ and the corresponding matrices $R_\jph$ and $R^{-1}_\jph$ of the 1-D Euler equation of gas dynamics can be found in, e.g., \cite{CCHKL_22}.
\end{remark}

\section{1-D Local Characteristic Decomposition Based Fifth-Order WENO-Z Interpolant}\label{appb}
In this appendix, we briefly describe the 1-D LCD based fifth-order WENO-Z interpolant.

Given the point values $\mU_j$ at uniform grid points $x=x_j$, the value $(U^{(i)})^-_\jph$ is computed using a weighted average of the three parabolic interpolants $\bm{{\cal P}}_0(x)$, $\bm{{\cal P}}_1(x)$ and $\bm{{\cal P}}_2(x)$
obtained using the stencils $[x_{j-2},x_{j-1},x_j]$, $[x_{j-1},x_j,x_{j+1}]$, and $[x_j,x_{j+1},x_{j+2}]$, respectively:
\begin{equation}
(U^{(i)})^-_\jph=\sum_{k=0}^2\omega^{(i)}_k{\cal P}^{(i)}_k(x_\jph),
\label{C1}
\end{equation}
where
\begin{equation}
\begin{aligned}
&{\cal P}^{(i)}_0(x_\jph)=\frac{3}{8}U^{(i)}_{j-2}-\frac{5}{4}U^{(i)}_{j-1}+\frac{15}{8}U^{(i)}_j,\\
&{\cal P}^{(i)}_1(x_\jph)=-\frac{1}{8}U^{(i)}_{j-1}+\frac{3}{4}U^{(i)}_j+\frac{3}{8}U^{(i)}_{j+1},\\
&{\cal P}^{(i)}_2(x_\jph)=\frac{3}{8}U^{(i)}_j+\frac{3}{4}U^{(i)}_{j+1}-\frac{1}{8}U^{(i)}_{j+2}.
\end{aligned}
\label{C2}
\end{equation}
To ensure \eref{C1}--\eref{C2} is fifth-order accurate and nonoscillatory, one can take the weights $\omega_k$
in \eref{C1} to be
\begin{equation}
\omega^{(i)}_k:=\frac{\alpha^{(i)}_k}{\alpha^{(i)}_0+\alpha^{(i)}_1+\alpha^{(i)}_2},\quad
\alpha^{(i)}_k=d_k\left[1+\left(\frac{\tau^{(i)}_5}{\beta^{(i)}_k+\varepsilon}\right)^p\right],\quad\tau_5=|\beta^{(i)}_2-\beta^{(i)}_0|,
\label{C3}
\end{equation}
where $d_0=\frac{1}{16}$, $d_1=\frac{5}{8}$, and $d_2=\frac{5}{16}$, and the smoothness indicators $\beta^{(i)}_k$ are given by 
\begin{equation}
\begin{aligned}
&\beta^{(i)}_0=\frac{13}{12}\big(U^{(i)}_{j-2}-2U^{(i)}_{j-1}+U^{(i)}_j\big)^2+\frac{1}{4}\big(U^{(i)}_{j-2}-4U^{(i)}_{j-1}+3U^{(i)}_j\big)^2,\\
&\beta^{(i)}_1=\frac{13}{12}\big(U^{(i)}_{j-1}-2U^{(i)}_j+U^{(i)}_{j+1}\big)^2+\frac{1}{4}\big(U^{(i)}_{j-1}-U^{(i)}_{j+1}\big)^2,\\
&\beta^{(i)}_2=\frac{13}{12}\big(U^{(i)}_j-2U^{(i)}_{j+1}+U^{(i)}_{j+2}\big)^2+\frac{1}{4}\big(3U^{(i)}_{j}-4U^{(i)}_{j+1}+U^{(i)}_{j+2}\big)^2.
\end{aligned}
\label{C4}
\end{equation}
We have used $p=2$ and $\varepsilon=10^{-12}$ in all of the numerical examples reported in this paper. The corresponding right-sided value, $(U^{(i)}_\jph)^-$, can also be derived using a mirror-symmetric approach. 

As in Appendix \ref{appa}, to ensure the nonoscillatory nature of the reconstruction \eref{C1}--\eref{C4}, we also need to adopt the reconstruction procedure within the LCD framework. To this end, we first introduce the local characteristic variables in the neighborhood of $x=x_\jph$:
$$
\bm\Gamma_m=R^{-1}_\jph\mU_m,\quad m=j-2,\ldots,j+3.
$$
Equipped with the values $\bm\Gamma_{j-2}$, $\bm\Gamma_{j-1}$,  $\bm\Gamma_j$, $\bm\Gamma_{j+1}$, $\bm\Gamma_{j+2}$, and $\bm\Gamma_{j+3}$, we then apply the interpolation procedure \eref{C1}--\eref{C4} to each of the components $\Gamma^{(i)}$, $i=1,\ldots,d$ of $\bm\Gamma$ to obtain ${\bm\Gamma}_\jph^-$ (the values of ${\bm\Gamma}_\jph^+$ are computed in the mirror-symmetric way). Finally, the corresponding one-sided point values of $\mU$ are given by
\begin{equation*}
\mU^\pm_\jph=R_\jph{\bm\Gamma}^\pm_\jph.
\end{equation*}

\bibliographystyle{siamnodash}
\bibliography{ref}
\end{document}